\documentclass[12pt]{article}
\usepackage{latexsym, makeidx,amssymb,amsmath,stmaryrd}
\begin{document}
\baselineskip=17pt

\newcommand{\la}{\langle}
\newcommand{\ra}{\rangle}
\newcommand{\psp}{\vspace{0.4cm}}
\newcommand{\pse}{\vspace{0.2cm}}
\newcommand{\ptl}{\partial}
\newcommand{\dlt}{\delta}
\newcommand{\sgm}{\sigma}
\newcommand{\al}{\alpha}
\newcommand{\be}{\beta}
\newcommand{\G}{\Gamma}
\newcommand{\g}{\gamma}
\newcommand{\gm}{\gamma}
\newcommand{\vs}{\varsigma}
\newcommand{\Lmd}{\Lambda}
\newcommand{\lmd}{\lambda}
\newcommand{\td}{\tilde}
\newcommand{\vf}{\varphi}
\newcommand{\yt}{Y^{\nu}}
\newcommand{\wt}{\mbox{wt}\:}
\newcommand{\der}{\mbox{Der}\:}
\newcommand{\ad}{\mbox{ad}\:}
\newcommand{\stl}{\stackrel}
\newcommand{\ol}{\overline}
\newcommand{\ul}{\underline}
\newcommand{\es}{\epsilon}
\newcommand{\dmd}{\diamond}
\newcommand{\clt}{\clubsuit}
\newcommand{\vt}{\vartheta}
\newcommand{\ves}{\varepsilon}
\newcommand{\dg}{\dagger}
\newcommand{\tr}{\mbox{Tr}\:}
\newcommand{\ga}{{\cal G}({\cal A})}
\newcommand{\hga}{\hat{\cal G}({\cal A})}
\newcommand{\Edo}{\mbox{End}\:}
\newcommand{\for}{\mbox{for}}
\newcommand{\kn}{\mbox{ker}\:}
\newcommand{\Dlt}{\Delta}
\newcommand{\rad}{\mbox{Rad}}
\newcommand{\rta}{\rightarrow}
\newcommand{\mbb}{\mathbb}
\newcommand{\rd}{\mbox{Res}\:}
\newcommand{\stc}{\stackrel{\circ}}
\newcommand{\stdt}{\stackrel{\bullet}}
\newcommand{\lra}{\Longrightarrow}
\newcommand{\f}{\varphi}
\newcommand{\rw}{\rightarrow}
\newcommand{\op}{\oplus}
\newcommand{\co}{\Omega}
\newcommand{\dar}{\Longleftrightarrow}
\newcommand{\sta}{\theta}

\begin{center}{\Large \bf  Oscillator Variations  of the Classical
}\end{center}
\begin{center}{\Large \bf Theorem on Harmonic Polynomials} 
\end{center}
\vspace{0.2cm}

\begin{center}{\large Cuiling Luo$^a$ and Xiaoping
Xu$^b$}\end{center}

{\small \noindent a. College of Science, Hebei Polytechnic
University, Tangshan, Hebei 063009,\\ P. R. China.

\noindent b. Corresponding author, Hua Loo-Keng Key Mathematical
Laboratory, Institute of\\ Mathematics, Academy of  Mathematics and
Systems Sciences, Chinese Academy of \\ Sciences, Beijing, 100190,
P. R. China.}

\begin{abstract}{ We study two-parameter oscillator variations of the classical
theorem on harmonic polynomials, associated with noncanonical
oscillator representations of $sl(n,\mbb{F})$  and $o(n,\mbb{F})$.
We find the condition when the homogeneous solution spaces of the
variated Laplace equation are irreducible modules of the concerned
algebras and the homogeneous subspaces are direct sums of the images
of these solution subspaces under the powers of the dual
differential operator. This establishes a local
$(sl(2,\mbb{F}),sl(n,\mbb{F}))$ and $(sl(2,\mbb{F}),o(n,\mbb{F}))$
Howe duality, respectively. In generic case, the obtained
irreducible $o(n,\mbb{F})$-modules are infinite-dimensional
non-unitary modules without highest-weight vectors. As an
application, we determine the structure of noncanonical oscillator
representations of $sp(2n,\mbb{C})$. When both  parameters are equal
to the maximal allowed value, we obtain an infinite family of
explicit irreducible $({\cal G}, {\cal K})$-modules for
$o(n,\mbb{F})$ and $sp(2n,\mbb{C})$. Methodologically we have
extensively used partial differential equations to solve
representation problems.}\end{abstract}

\section{Introduction}

Harmonic polynomials are important objects in analysis, differential
geometry and physics. A fundamental theorem in classical harmonic
analysis says that the spaces of homogeneous harmonic polynomials
(solutions of Laplace equation) are irreducible modules of the
corresponding orthogonal Lie group (algebra) and the whole
polynomial algebra is a free module over the invariant polynomials
generated by harmonic polynomials. Bases of these irreducible
modules can be obtained easily (e.g., cf. [X]). The algebraic beauty
of the above theorem is that Laplace equation characterizes the
irreducible submodules of the polynomial algebra and the
corresponding quadratic invariant gives a decomposition of the
polynomial algebra into a direct sum of irreducible submodules. This
actually forms  an $(sl(2,\mbb{F}),o(n,\mbb{F}))$ Howe duality.

Lie algebras (Lie groups) serve as the symmetries in quantum physics
(e.g., cf. [FC, L, LF, G]). Their various representations provide
distinct concrete practical physical models. Many important physical
phenomena have been interpreted as the consequences of symmetry
breaking (e.g., cf. [LF]). Harmonic oscillators are basic objects in
quantum mechanics (e.g., cf. [FC, G]). Oscillator representations of
finite-dimensional simple Lie algebras are the most fundamental ones
in quantum physics. Their infinite-dimensional analogues are free
field representations of affine Kac-Moody algebras.

 The aim of this work is to establish certain two-parameter oscillator
variations of the classical theorem on harmonic polynomials,
associated with noncanonical oscillator representations of special
linear Lie algebras and orthogonal Lie algebras, which are obtained
by swapping differential operators and multiplication operators in
the canonical oscillator representations induced from the natural
representations. The Howe duality does not hold on the whole
polynomial algebras. But we find the condition when the homogeneous
solution spaces of the variated Laplace equation are irreducible
modules of the concerned algebras and the homogeneous subspaces are
direct sums of the images of these solution subspaces under the
powers of the dual differential operator. We may call this a {\it
local $(sl(2,\mbb{F}),sl(n,\mbb{F}))$ and
$(sl(2,\mbb{F}),o(n,\mbb{F}))$ Howe duality,} respectively. In
particular, we obtain explicit infinite-dimensional non-unitary
modules of orthogonal Lie algebras that are not of highest-weight
type. As an application of our results on special linear Lie
algebras, we prove that the homogeneous subspaces in noncanonical
oscillator representations of symplectic Lie algebras are
irreducible except some singular cases, in which the homogeneous
subspaces are direct sums of exactly two explicitly given
irreducible submodules. Explicit bases of all the above irreducible
modules in generic case are obtained.

Let ${\cal G}$ be a semisimple Lie algebra and let ${\cal K}$ be a
maximal proper reductive Lie subalgebra of ${\cal G}$. An
infinite-dimensional irreducible ${\cal G}$-module is said of
$({\cal G},{\cal K})$-{\it type} if it is a direct sum of
finite-dimensional irreducible ${\cal K}$-submodules. When both
parameters are equal to the maximal allowed value, we obtain an
infinite family of explicit irreducible $({\cal G}, {\cal
K})$-modules for orthogonal Lie algebras and symplectic Lie
algebras. Since our representations are not unitary, the concerned
modules are infinite-dimensional and we are dealing with pairs of
dual invariant differential operators, traditional methods fail to
solve our problems. In fact, we have extensively used the method of
solving flag partial differential equations developed in [X] by the
second author. Below we give a technical introduction.

 For convenience, we will use the notion
$\ol{i,i+j}=\{i,i+1,i+2,...,i+j\}$ for integers $i$ and $j$ with
$i\leq j$. Denote by $\mbb{N}$ the additive semigroup of nonnegative
integers.

Let $E_{r,s}$ be the square matrix with 1 as its $(r,s)$-entry and 0
as the others. The  compact orthogonal Lie algebra
$o(n,\mbb{R})=\sum_{1\leq r<s\leq n}\mbb{R}(E_{r,s}-E_{s,r}),$ whose
representation  on the polynomial algebra ${\cal
A}=\mbb{R}[x_1,...,x_n]$ is determined by $(E_{r,s}-E_{s,r})|_{\cal
A}=x_r\ptl_{x_s}-x_s\ptl_{x_r}$, which we call  the {\it canonical
oscillator representation of $o(n,\mbb{R})$} (e.g., cf. [FSS]).
Denote by ${\cal A}_k$ the subspace of homogeneous polynomials in
${\cal A}$ with degree $k$. Recall that the Laplace operator
$\Dlt=\ptl_{x_1}^2+\cdots+\ptl_{x_n}^2$ and its corresponding
invariant $\eta=x_1^2+x_2^2+\cdots+x_n^2$. When $n\geq 3$, it is
well known that the subspace of harmonic polynomials
$${\cal H}_k=\{f\in{\cal A}_k\mid
\Dlt(f)=0\}\eqno(1.1)$$ forms an irreducible $o(n,\mbb{R})$-module
and ${\cal A}_k={\cal H}_k\oplus\eta{\cal A}_{k-2},$  which is
equivalent to that ${\cal A}_k=\bigoplus_{i=1}^{\llbracket k/2
\rrbracket}\eta^i{\cal H}_{k-2i}$ is a direct sum of irreducible
submodules. Since the space
$\mbb{F}\Dlt+\mbb{F}[\Dlt,\eta]+\mbb{F}\eta$ forms an operator Lie
algebra isomorphic to $sl(2,\mbb{R})$, the above conclusion gives an
$(sl(2,\mbb{R}),o(n,\mbb{R}))$ Howe duality.

Below all the vector spaces are assumed over a field $\mbb{F}$ with
characteristic 0. Moreover, we always assume that $n\geq 2$ is an
integer. Let us reconsider the canonical oscillator representation
of $sl(n,\mbb{F})$:
$$E_{i,j}|_{\cal A}=x_i\ptl_j\qquad \for
\;\;i,j\in\ol{1,n}.\eqno(1.2)$$ Fix $1\leq n_1<n$. Note
$$[\ptl_{x_r},x_r]=1=[-x_r,\ptl_{x_r}].\eqno(1.3)$$
Changing operators $\ptl_{x_r}\mapsto -x_r$ and $x_r\mapsto
\ptl_{x_r}$ in (1.2) for $r\in\ol{1,n_1}$, we obtain the following
noncanonical oscillator representation of $sl(n,\mbb{F})$ determined
by:
$$E_{i,j}|_{\cal A}=\left\{\begin{array}{ll}-x_j\ptl_{x_i}-\delta_{i,j}&\mbox{if}\;
i,j\in\ol{1,n_1};\\ \ptl_{x_i}\ptl_{x_j}&\mbox{if}\;i\in\ol{1,n_1},\;j\in\ol{n_1+1,n};\\
-x_ix_j &\mbox{if}\;i\in\ol{n_1+1,n},\;j\in\ol{1,n_1};\\
x_i\partial_{x_j}&\mbox{if}\;i,j\in\ol{n_1+1,n}.
\end{array}\right.\eqno(1.4)$$
For any $k\in\mbb{Z}$, we denote
$${\cal A}_{\la
k\ra}=\mbox{Span}\:\{x^\al\mid\al\in\mbb{N}\:^n;\sum_{r=n_1+1}^n\al_r-\sum_{i=1}^{n_1}\al_i=k\}.
\eqno(1.5)$$ It was presented by Howe in his work [Ho] that for
$m_1,m_2\in\mbb{N}$ with $m_1>0$, ${\cal A}_{\la -m_1\ra}$ is an
irreducible highest-weight $sl(n,\mbb{F})$-submodule with highest
weight $m_1\lmd_{n_1-1}-(m_1+1)\lmd_{n_1}$ and ${\cal A}_{\la
m_2\ra}$ is an irreducible highest-weight $sl(n,\mbb{F})$-submodule
with highest weight
$-(m_2+1)\lmd_{n_1}+m_2(1-\dlt_{n_1,n-1})\lmd_{n_1+1}$.

 Denote ${\cal B}=\mbb{F}[x_1,...,x_n,y_1,...,y_n]$. Fix
 $n_1,n_2\in\ol{1,n}$ with $n_1\leq n_2$. Changing operators $\ptl_{x_r}\mapsto -x_r,\;
 x_r\mapsto
\ptl_{x_r}$  for $r\in\ol{1,n_1}$ and $\ptl_{y_s}\mapsto -y_s,\;
 y_s\mapsto\ptl_{y_s}$  for $s\in\ol{n_2+1,n}$, we get another noncanonical oscillator representation of $sl(n,\mbb{F})$ on
  ${\cal B}$ determined by
$$E_{i,j}|_{\cal
B}=E_{i,j}^x-E_{j,i}^y\qquad\for\;\;i,j\in\ol{1,n}\eqno(1.6)$$ with
$$E_{i,j}^x|_{\cal B}=\left\{\begin{array}{ll}-x_j\ptl_{x_i}-\delta_{i,j}&\mbox{if}\;
i,j\in\ol{1,n_1};\\ \ptl_{x_i}\ptl_{x_j}&\mbox{if}\;i\in\ol{1,n_1},\;j\in\ol{n_1+1,n};\\
-x_ix_j &\mbox{if}\;i\in\ol{n_1+1,n},\;j\in\ol{1,n_1};\\
x_i\partial_{x_j}&\mbox{if}\;i,j\in\ol{n_1+1,n}
\end{array}\right.\eqno(1.7)$$
and
$$E_{i,j}^y|_{\cal B}=\left\{\begin{array}{ll}y_i\ptl_{y_j}&\mbox{if}\;
i,j\in\ol{1,n_2};\\ -y_iy_j&\mbox{if}\;i\in\ol{1,n_2},\;j\in\ol{n_2+1,n};\\
\ptl_{y_i}\ptl_{y_j} &\mbox{if}\;i\in\ol{n_2+1,n},\;j\in\ol{1,n_2};\\
-y_j\partial_{y_i}-\delta_{i,j}&\mbox{if}\;i,j\in\ol{n_2+1,n}.
\end{array}\right.\eqno(1.8)$$
The related variated  Laplace operator becomes
$${\cal
D}=-\sum_{i=1}^{n_1}x_i\ptl_{y_i}+\sum_{r=n_1+1}^{n_2}\ptl_{x_r}\ptl_{y_r}-\sum_{s=n_2+1}^n
y_s\ptl_{x_s}\eqno(1.9)$$ and its dual
$$\eta=\sum_{i=1}^{n_1}y_i\ptl_{x_i}+\sum_{r=n_1+1}^{n_2}x_ry_r+\sum_{s=n_2+1}^n
x_s\ptl_{y_s}.\eqno(1.10)$$
 Set
$${\cal B}_{\la \ell_1,\ell_2\ra}=\mbox{Span}\{x^\al
y^\be\mid\al,\be\in\mbb{N}\:^n;\sum_{r=n_1+1}^n\al_r-\sum_{i=1}^{n_1}\al_i=\ell_1;
\sum_{i=1}^{n_2}\be_i-\sum_{r=n_2+1}^n\be_r=\ell_2\}\eqno(1.11)$$
for $\ell_1,\ell_2\in\mbb{Z}$. Define
$${\cal H}_{\la\ell_1,\ell_2\ra}=\{f\in {\cal B}_{\la
\ell_1,\ell_2\ra}\mid {\cal D}(f)=0\}.\eqno(1.12)$$  The following
is our first result:\psp

{\bf Theorem 1}. {\it For any $\ell_1,\ell_2\in\mbb{Z}$ such that
$\ell_1+\ell_2\leq n_1-n_2+1-\dlt_{n_1,n_2}$, ${\cal
H}_{\la\ell_1,\ell_2\ra}$ is an irreducible highest-weight
$sl(n,\mbb{F})$-module and ${\cal
B}_{\la\ell_1,\ell_2\ra}=\bigoplus_{m=0}^\infty\eta^m({\cal
H}_{\la\ell_1-m,\ell_2-m\ra})$ is a decomposition of irreducible
submodules. In particular, ${\cal B}_{\la\ell_1,\ell_2\ra}={\cal
H}_{\la\ell_1,\ell_2\ra}\oplus \eta({\cal
B}_{\la\ell_1-1,\ell_2-1\ra})$.}\psp

When $n_1+1<n_2<n$ and $\ell_1+\ell_2>n_1-n_2+1$, ${\cal
H}_{\la\ell_1,\ell_2\ra}$ is not irreducible and contains nonzero
elements in $\eta({\cal B}_{\la\ell_1-1,\ell_2-1\ra})$. Although the
space $\mbb{F}{\cal D}+\mbb{F}[{\cal D},\eta]+\mbb{F}\eta$ forms an
operator Lie algebra isomorphic to $sl(2,\mbb{R})$, we do not have
an $(sl(2,\mbb{F}),sl(n,\mbb{F}))$ Howe duality. We may call Theorem
1 an {\it local $(sl(2,\mbb{F}),sl(n,\mbb{F}))$ Howe duality.}

Consider the split
$$o(2n,\mathbb{F})=\sum_{i,j=1}^n\mathbb{F}(E_{i,j}-E_{n+j,n+i})+\sum_{1\leq
i<j\leq
n}[\mathbb{F}(E_{i,n+j}-E_{j,n+i})+\mathbb{F}(E_{n+j,i}-E_{n+i,j})]\eqno(1.13)$$
and  define a noncanonical oscillator representation of
$o(2n,\mbb{F})$ on ${\cal B}$ by
$$(E_{i,j}-E_{n+j,n+i})|_{\cal B}=E_{i,j}^x|_{\cal B}-E_{j,i}^y|_{\cal B},\eqno(1.14)$$
$$E_{i,n+j}|_{\cal B}=\left\{\begin{array}{ll}
\ptl_{x_i}\ptl_{y_j}&\mbox{if}\;i\in\ol{1,n_1},\;j\in\ol{1,n_2},\\
-y_j\ptl_{x_i}&\mbox{if}\;i\in\ol{1,n_1},\;j\in\ol{n_2+1,n},\\
x_i\ptl_{y_j}&\mbox{if}\;i\in\ol{n_1+1,n},\;j\in\ol{1,n_2},\\
-x_iy_j&\mbox{if}\;i\in\ol{n_1+1,n},\;j\in\ol{n_2+1,n}\end{array}\right.\eqno(1.15)$$
and
$$E_{n+i,j}|_{\cal B}=\left\{\begin{array}{ll}
-x_jy_i&\mbox{if}\;j\in\ol{1,n_1},\;i\in\ol{1,n_2},\\
-x_j\ptl_{y_i}&\mbox{if}\;j\in\ol{1,n_1},\;i\in\ol{n_2+1,n},\\
y_i\ptl_{x_j}&\mbox{if}\;j\in\ol{n_1+1,n},\;i\in\ol{1,n_2},\\
\ptl_{x_j}\ptl_{y_i}&\mbox{if}\;j\in\ol{n_1+1,n},\;i\in\ol{n_2+1,n}.\end{array}\right.\eqno(1.16)$$
Set
$${\cal B}_{\la
k\ra}=\bigoplus_{\ell_1,\ell_2\in\mbb{Z};\ell_1+\ell_2=k}{\cal
B}_{\la\ell_1,\ell_2\ra},\qquad {\cal H}_{\la k\ra}=\{f\in {\cal
B}_{\la k\ra}\mid {\cal D}(f)=0\}.\eqno(1.17)$$ Below we always take
${\cal K}=\sum_{i,j=1}^n\mbb{F}(E_{i,j}-E_{n+j,n+i})$. Our second
results is:\psp

{\bf Theorem 2}. {\it For any  $n_1-n_2+1-\dlt_{n_1,n_2}\geq
k\in\mbb{Z}$, ${\cal H}_{\la k\ra}$ is an irreducible
$o(2n,\mbb{F})$-submodule and ${\cal B}_{\la
k\ra}=\bigoplus_{i=0}^\infty\eta^i({\cal H}_{\la k-2i\ra})$ is a
decomposition of irreducible submodules. In particular, ${\cal
B}_{\la k\ra}={\cal H}_{\la k\ra}\oplus \eta({\cal B}_{\la
k-2\ra})$. The module ${\cal H}_{\la k\ra}$ under the assumption is
of highest-weight type only if $n_2=n$. When $n_1=n_2=n$, all the
irreducible modules ${\cal H}_{\la k\ra}$ with $0\geq k\in\mbb{Z}$
are of $({\cal G},{\cal K})$-type.}\psp

We may view Theorem 2 as an {\it local
$(sl(2,\mbb{F}),o(2n,\mbb{F}))$ Howe duality.}

 Note the split
$$o(2n+1,\mbb{F})=o(2n,\mbb{F})\oplus\bigoplus_{i=1}^n
[\mbb{F}(E_{0,i}-E_{n+i,0})+\mbb{F}(E_{0,n+i}-E_{i,0})].\eqno(1.18)$$
Let ${\cal B}'=\mbb{F}[x_0,x_1,...,x_n,y_1,...,y_n]$. We define a
noncanonical oscillator representation of $o(2n+1,\mbb{F})$ on
${\cal B}'$ by the differential operators in (1.14)-(1.16) and
$$ E_{0,i}|_{{\cal
B}'}=\left\{\begin{array}{ll}-x_0x_i&\mbox{if}\;i\in\ol{1,n_1},\\
x_0\ptl_{x_i}&\mbox{if}\;i\in\ol{n_1+1,n},\\
x_0\ptl_{y_i}&\mbox{if}\;i\in\ol{n+1,n+n_2},\\
-x_0y_i&\mbox{if}\;i\in\ol{n+n_2+1,2n}\end{array}\right.\eqno(1.19)$$
and
$$ E_{i,0}|_{{\cal
B}'}=\left\{\begin{array}{ll}\ptl_{x_0}\ptl_{x_i}&\mbox{if}\;i\in\ol{1,n_1},\\
x_i\ptl_{x_0}&\mbox{if}\;i\in\ol{n_1+1,n},\\
y_i\ptl_{x_0}&\mbox{if}\;i\in\ol{n+1,n+n_2},\\
\ptl_{x_0}\ptl_{y_i}&\mbox{if}\;i\in\ol{n+n_2+1,2n}.\end{array}\right.\eqno(1.20)$$
Now the variated Laplace operator becomes
$${\cal
D}'=\ptl_{x_0}^2-2\sum_{i=1}^{n_1}x_i\ptl_{y_i}+2\sum_{r=n_1+1}^{n_2}\ptl_{x_r}\ptl_{y_r}-2\sum_{s=n_2+1}^n
y_s\ptl_{x_s}\eqno(1.21)$$ and its dual operator
$$\eta'=x_0^2+2\sum_{i=1}^{n_1}y_i\ptl_{x_i}+2\sum_{r=n_1+1}^{n_2}x_ry_r+2\sum_{s=n_2+1}^n
x_s\ptl_{y_s}.\eqno(1.22)$$ Set
$${\cal B}'_{\la k\ra}=\sum_{i=0}^\infty {\cal B}_{\la
k-i\ra}x_0^i,\qquad {\cal H}'_{\la k\ra}=\{f\in {\cal B}'_{\la
k\ra}\mid {\cal D}'(f)=0\}.\eqno(1.23)$$ The following is our third
result.\psp

{\bf Theorem 3}. {\it  For any  $n_1-n_2+1-\dlt_{n_1,n_2}\geq
k\in\mbb{Z}$, ${\cal H}'_{\la k\ra}$ is an irreducible
$o(2n+1,\mbb{F})$-submodule and ${\cal B}'_{\la
k\ra}=\bigoplus_{i=0}^\infty(\eta')^i({\cal H}'_{\la k-2i\ra})$ is a
decomposition of irreducible submodules. In particular, ${\cal
B}'_{\la k\ra}={\cal H}'_{\la k\ra}\oplus \eta'({\cal B}'_{\la
k-2\ra})$. The module ${\cal H}_{\la k\ra}'$ under the assumption is
of highest-weight type only if $n_2=n$.  When $n_1=n_2=n$, all the
irreducible modules ${\cal H}_{\la k\ra}'$ with $0\geq k\in\mbb{Z}$
are of $({\cal G},{\cal K})$-type.} \psp

Again Theorem 2 can be viewed as an {\it local
$(sl(2,\mbb{F}),o(2n+1,\mbb{F}))$ Howe duality.}

Define a noncanonical oscillator representation of $sp(2n,\mbb{F})$
on ${\cal B}$ by (1.14)-(1.16). Using some results in the proof of
Theorem 1, we prove:\psp

{\bf Theorem 4}. {\it Let $k\in\mbb{Z}$. If $n_1<n_2$ or $k\neq 0$,
the subspace ${\cal B}_{\la k\ra}$ (cf. (1.17)) is an irreducible
$sp(2n,\mbb{F})$-submodule.  When $n_1=n_2$, the subspace ${\cal
B}_{\la 0\ra}$ is a direct sum of two irreducible
$sp(2n,\mbb{F})$-submodules. Moreover, each irreducible submodule is
of highest-weight module only if $n_2=n$.
 When $n_1=n_2=n$, all the
irreducible submodules are of $({\cal G},{\cal K})$-type.}\psp

 In addition,  the explicit expressions for
all the above irreducible modules are given. In the case of
highest-weight type,  the highest-weight vector and its weight of
the corresponding irreducible modules are also presented. Since the
representations with parameters $(n_1,n_2)$ are contragredient to
those with parameters $(n-n_2,n-n_1)$, the case $n_2<n_1$ has
virtually been handled.

In Section 2, we present some preparatory works, in particular, the
method of solving flag partial differential equations found in [X]
by the second author. In Section 3, we prove Theorem 1 when
$n_1<n_2$. Section 4 is devoted to the proof of Theorem 1 with
$n_1=n_2$. In Sections 5, 6 and 7, we prove Theorems 2, 3 and 4,
respectively.

\section{Preparation}

It is very often that Lie group theorists characterize certain
irreducible modules as kernels of a set of differential operators.
But how to solve the corresponding systems of partial differential
equations is in general unknown. It was realized by the second
author that these equations are of ``flag type" when the modules are
of highest-weight type. A linear transformation (operator) $T$ on a
vector space $V$ is called {\it locally nilpotent} if for any $v\in
V$, there exists a positive integer $k$ such that $T^k(v)=0$. A {\it
partial differential equation of flag type} is the linear
differential equation of the form:
$$(d_1+f_1d_2+f_2d_3+\cdots+f_{n-1}d_n)(u)=0,\eqno(2.1)$$
where $d_1,d_2,...,d_n$ are certain commuting locally nilpotent
differential operators on the polynomial algebra
$\mbb{F}[x_1,x_2,...,x_n]$ and $f_1,...,f_{n-1}$ are polynomials
satisfying $d_i(f_j)=0$ if $i>j.$ Many variable-coefficient
(generalized) Laplace equations, wave equations, Klein-Gordon
equations, Helmholtz equations are of this type. Solving such
equations is also important in finding invariant solutions of
nonlinear partial differential equations (e.g., cf. [I1, I2]). In
representation theory, we are more interested in polynomial
solutions of flag partial differential equations. The second author
[X] found an effective way of solving for them. The following lemma
is a slightly generalized form of Lemma 2.1 in [X].\psp

{\bf Lemma 2.1 (Xu [X])}. {\it Let ${\cal B}$ be a commutative
associative algebra and let ${\cal A}$ be a free ${\cal B}$-module
 generated  by a filtrated subspace $V=\bigcup_{r=0}^\infty V_r$
(i.e., $V_r\subset V_{r+1}$). Let $T_1$ be a linear operator on
${\cal B}\oplus {\cal A}$ with a right inverse $T_1^-$ such that
$$T_1({\cal B},{\cal A}),\;T_1^-({\cal B},{\cal A})\subset({\cal B},{\cal A}),\;\; T_1(\eta_1\eta_2)=T_1(\eta_1)\eta_2,\;\;
T_1^-(\eta_1\eta_2)=T_1^-(\eta_1)\eta_2 \eqno(2.2)$$ for $\eta_1 \in
{\cal B},\;\eta_2\in V$, and let $T_2$ be a linear operator on
${\cal A}$ such that
$$ T_2(V_{r+1})\subset {\cal B}V_r,\;\;
T_2(f\zeta)=fT_2(\zeta) \qquad\for\;\; r\in\mbb{N},\;\;f\in{\cal
B},\;\zeta\in{\cal A}.\eqno(2.3)$$ Then we
have \begin{eqnarray*}\hspace{1cm}&&\{g\in{\cal A}\mid (T_1+T_2)(g)=0\}\\
& =&\mbox{Span}\{ \sum_{i=0}^\infty(-T_1^-T_2)^i(hg)\mid g\in
V,\;h\in {\cal B};\;T_1(h)=0\}. \hspace{3.9cm}(2.4)\end{eqnarray*} }
\psp

Set
$$\es_i=(0,...,0,\stl{i}{1},0,...,0)\in \mbb{N}^{\:n}.\eqno(2.5)$$
 For each
$i\in\ol{1,n}$, we define the linear operator $\int_{(x_i)}$ on
${\cal A}$ by:
$$\int_{(x_i)}(x^\al)=\frac{x^{\al+\es_i}}{\al_i+1}\;\;\for\;\;\al\in
\mbb{N}^{\:n}.\eqno(2.6)$$ Furthermore, we let
$$\int_{(x_i)}^{(0)}=1,\qquad\int_{(x_i)}^{(m)}=\stl{m}{\overbrace{\int_{(x_i)}\cdots\int_{(x_i)}}}
\qquad\for\; \;0<m\in\mbb{Z}\eqno(2.7)$$ and denote
$$\ptl^{\al}=\ptl_{x_1}^{\al_1}\ptl_{x_2}^{\al_2}\cdots
\ptl_{x_n}^{\al_n},\;\;
\int^{(\al)}=\int_{(x_1)}^{(\al_1)}\int_{(x_2)}^{(\al_2)}\cdots
\int_{(x_n)}^{(\al_n)}\qquad\for\;\;\al\in
\mbb{N}^{\:n}.\eqno(2.8)$$ Obviously, $\int^{(\al)}$ is a right
inverse of $\ptl^\al$ for $\al\in \mbb{N}^{\:n}.$ We remark that
$\int^{(\al)}\ptl^\al\neq 1$ if $\al\neq 0$ due to $\ptl^\al(1)=0$.
In this paper,  our $T_1$'s are of the type $\ptl^\al$ and the right
inverse $T_1^-=\int^{(\al)}$.

Let $m_1,m_2,...,m_n$ be positive integers. Taking
$T_1=\ptl_{x_1}^{m_1},\;T_2=\ptl_{x_2}^{m_2}+\cdots+\ptl_{x_n}^{m_n}$
and $T_1^-=\int_{(x_1)}^{(m_1)}$, we find that the set
\begin{eqnarray*}\hspace{1.9cm}& &\{\sum_{k_2,...,k_n=0}^\infty(-1)^{k_2+\cdots+k_n}{k_2+\cdots+k_k\choose
k_2,...,k_n} \int_{(x_1)}^{((k_2+\cdots +k_n)m_1)}(x_1^{\ell_1})\\
& &\times\ptl_{x_2}^{k_2m_2}(x_2^{\ell_2})\cdots
\ptl_{x_n}^{k_nm_n}(x_n^{\ell_n})\mid
\ell_1\in\ol{0,m_1-1},\;\ell_2,...,\ell_n\in\mbb{N}\}\hspace{2cm}
(2.9)\end{eqnarray*} forms a basis of the space of polynomial
solutions for the equation
$$(\ptl_{x_1}^{m_1}+\ptl_{x_2}^{m_2}+\cdots+\ptl_{x_n}^{m_n})(u)=0.\eqno(2.10)$$
When all $m_i=2$, we get a basis of the space of harmonic
polynomials.

Cao [C] used Lemma 2.1 to prove that the subspaces of homogeneous
polynomial vector solutions of the $n$-dimensional Navier equations
in elasticity are exactly direct sums of three explicitly given
irreducible submodules when $n\neq 4$ and direct sums of four
explicitly given irreducible submodules if $n=4$ of the
corresponding orthogonal Lie group (algebra), and the whole
polynomial vector space is also a free module over the invariant
polynomials generated these solutions. The result can be viewed as a
vector generalization of the classical theorem on harmonic
polynomials.  Moreover, Cao solved the initial value problem for the
Navier equations based on the ideas in [X].

The idea of solving flag partial differential equation in [X] leads
the second author to find a family of special functions functions
$${\cal Y}_r(y_1,...,y_m)=\sum_{i_1,...,i_m=0}^\infty {i_1+\cdots+i_m\choose
i_1,...,i_m}\frac{y_1 ^{i_1}y_2 ^{i_2}\cdots y_m ^{i_m}}
{(r+\sum_{s=1}^msi_s)!},\eqno(2.11)$$by which we can solve the
initial value problem of the equation:
$$(\ptl_{x_1}^m-\sum_{r=1}^m\ptl_{x_1}^{m-i}f_i(\ptl_{x_2},...,\ptl_{x_n}))(u)=0,\eqno(2.12)$$
 where $f_i(\ptl_{x_2},...,\ptl_{x_n})\in\mbb{R}[\ptl_{x_2},...,\ptl_{x_n}].$

Let ${\cal A}=\mbb{F}[x_1,...,x_n]$ and let $gl(n,\mbb{F})$ act on
${\cal A}$ by (1.4). With the notion in (1.5), ${\cal
A}=\bigoplus_{k\in\mbb{Z}}{\cal A}_{\la k\ra}$ is a $\mbb{Z}$ graded
algebra and each homogeneous subspace ${\cal A}_{\la k\ra}$ is
infinite-dimensional. Set
$$\flat=\sum_{r=n_1+1}^nx_r\ptl_{x_r}-\sum_{i=1}^{n_1}x_i\ptl_{x_i}.\eqno(2.13)$$
Then
$${\cal A}_{\la k\ra}=\{f\in{\cal A}\mid\flat(f)=kf\}.\eqno(2.14)$$
Moreover, we have
$$\flat E_{i,j}=E_{j,i}\flat\;\;\mbox{on}\;\;{\cal
A}\qquad\for\;i,j\in\ol{1,n}.\eqno(2.15)$$ Thus ${\cal A}_{\la
k\ra}$ forms a ${\cal G}$-module for any subalgebra ${\cal G}$ of
$gl(n,\mbb{F})$.

For $\al\in\mbb{N}\:^n$, we denote $\alpha!=\prod_{i=1}^n\al_i!$ and
define a symmetric bilinear form $(\cdot|\cdot)$ on ${\cal A}$ by
$$(x^\alpha|x^\beta)=\dlt_{\alpha,\beta}(-1)^{\sum_{i=1}^{n_1}\alpha_i}\alpha!\qquad\for\;\;
\al,\be\in \mbb{N}\:^n.\eqno(2.16)$$ Then we have: \psp

{\bf Lemma 2.2}. {\it For any $A\in gl(n,\mbb{F})$ and
$f,g\in\mbb{\cal A}$, we have $(A(f)|g)=(f|A^t(g)),$ where $A^t$
denote the transpose of the matrix $A$.}

{\it Proof.} Let $\al,\be\in \mbb{N}\:^n$. For $i,j\in\ol{1,n_1}$,
$$(E_{i,j}(x^\al)|x^\be)=-\al_i(x^{\al+\es_j-\es_i}|x^\be)
-\dlt_{i,j}(x^\al|x^\be)\eqno(2.17)$$ and
$$(x^\al|E_{j,i}(x^\be))=-\be_j(x^\al|x^{\be+\es_i-\es_j})
-\dlt_{i,j}(x^\al|x^\be)\eqno(2.18)$$ by (1.4). Note
\begin{eqnarray*}\hspace{1cm}\al_i(x^{\al+\es_j-\es_i}|x^\be)&=&
\dlt_{\alpha+\es_j-\es_i,\beta}(-1)^{\sum_{i=1}^{n_1}\alpha_i}(\al_j+1)\alpha!
\\ &=&\be_j\dlt_{\alpha,\beta+\es_i-\es_j}(-1)^{\sum_{i=1}^{n_1}\alpha_i}\alpha!
=\be_j(x^\al|x^{\be+\es_i-\es_j})\hspace{2.7cm}(2.19)\end{eqnarray*}
by (2.16). Hence
$$(E_{i,j}(x^\al)|x^\be)=(x^\al|E_{j,i}(x^\be)).\eqno(2.20)$$
If $i,j\in\ol{n_1+1,n}$, then (2.19) holds and so does (2.20).

Consider $i\in\ol{1,n_1}$ and $j\in\ol{n_1+1,n}$.
$$(E_{i,j}(x^\al)|x^\be)=\al_i\al_j(x^{\al-\es_i-\es_j}|x^\be)=-
\dlt_{\alpha-\es_i-\es_j,\beta}(-1)^{\sum_{i=1}^{n_1}\alpha_i}\alpha!\eqno(2.21)$$
and
$$(x^\al|E_{j,i}(x^\be))=-(x^\al|x^{\be+\es_i+\es_j})=
-\dlt_{\alpha,\beta+\es_i+\es_j}(-1)^{\sum_{i=1}^{n_1}\alpha_i}\alpha!\eqno(2.22)$$
by (1.4) and (2.16). So (2.20) holds. Therefore, the lemma holds by
the symmetry of the form.$\qquad\Box$\psp

Let ${\cal G}$ be simple Lie subalgebra of $gl(n,\mbb{F})$ such that
$A^t\in{\cal G}$ if $A\in{\cal G}.$ Let $H$ be a Cartan subalgebra
of ${\cal G}$ and assume that ${\cal A}$ forms a weighted ${\cal
G}$-module with respect to $H$. Fix the positivity of roots and
denote by ${\cal G}_+$ the sum of positive root subspaces.  A {\it
singular vector} is a weight vector annihilated by positive root
vectors.

From now on,  we  count the number of singular vectors up to a
scalar multiple. Moreover, an element $g\in{\cal A}$ is called {\it
nilpotent with respect to} ${\cal G}_+$ if there exist a positive
integer $m$ such that
$$\xi_1\cdots \xi_m(g)=0\qquad\mbox{for
any}\;\xi_1,...,\xi_m\in{\cal G}_+.\eqno(2.23)$$ A subspace $V$ of
${\cal A}$ is called {\it nilpotent with respect to} ${\cal G}_+$ if
all its elements are nilpotent with respect to ${\cal G}_+$.  If the
elements of ${\cal G}_+|_{\cal A}$ are locally nilpotent and ${\cal
G}_+({\cal A}_i)\subset \sum_{r=0}^i{\cal A}_r$ for any
$i\in\mbb{N}$, then any element of ${\cal A}$ is nilpotent with
respect to ${\cal G}_+$ by Engel's Theorem.\psp

{\bf Lemma 2.3}. {\it If a submodule $N$ of ${\cal A}$ is nilpotent
with respect to ${\cal G}_+$, $N$ contains only one singular vector
$v$ and $(v|v)\neq0$, then $N$ is an irreducible summand of ${\cal
A}$}.

{\it Proof}. Under the nilpotent assumption, any nonzero submodule
of $N$ contains a singular vector. In particular,
$N_1=U(\mathcal{G})(v)$ is an irreducible submodule by the
uniqueness of singular vector.  Set
$$\bar N_1^\bot=\{u\in N|(u|w)=0\mid  w\in N_1\}.\eqno(2.24)$$ and
$${\cal R}=\{u\in N|(u|w)=0\mid w\in N\}.\eqno(2.25)$$ Note that
$\bar N_1^\bot$ and ${\cal R}$ are submodules of $N$ by Lemma 2.2.
If ${\cal R}\neq0$, it should contain a nonzero singular vector,
which is impossible according to the assumption $(v|v)\neq0$.
Therefore ${\cal R}=\{0\}$, and $N=N_1\bigoplus \bar N_1^\bot$. But
$\bar N_1^\bot=0$ by the same argument, and so $N=N_1$. The fact
${\cal R}=\{0\}$ implies that
$${\cal A}=N\oplus\{f\in{\cal A}\mid
(f|g)=0\;\for\;g\in N\}\eqno(2.26)$$ is a direct sum of ${\cal
G}$-submodules.$\qquad\Box$\psp

Let ${\cal Q}=\mbb{F}(x_1,...,x_n,y_1,...,y_n)$ be the space of
rational functions in $x_1,...,x_n,y_1,...,y_n$. Define a
representation of $sl(n,\mbb{F})$ on ${\cal Q}$ via
$$E_{i,j}|_{\cal
Q}=x_i\ptl_{x_j}-y_j\ptl_{y_i}\qquad\for\;\;i,j\in\ol{1,n}.\eqno(2.27)$$
Set $\zeta=\sum_{i=1}^n x_iy_i.$ Then
$$\xi(\zeta)=0\qquad\for\;\;\xi\in sl(n,\mbb{F}).\eqno(2.28)$$
Take
$$H=\sum_{i=1}^{n-1}\mbb{F}(E_{i,i}-E_{i+1,i+1})\eqno(2.29)$$
as a Cartan subalgebra of $sl(n,\mbb{F})$ and the subspace spanned
by positive root vectors:
$$sl(n,\mbb{F})_+=\sum_{1\leq i<j\leq}\mbb{F}E_{i,j}.\eqno(2.30)$$
The following lemma was proved in [X], which will be used in next
section.\psp

{\bf Lemma 4}. {\it Any singular function in ${\cal Q}$ is a
rational function in $x_1,y_n,\zeta$}.

\section{The $sl(n,\mbb{F})$-Variation I: $n_1<n_2$}

Fix $n_1,n_2\in\ol{1,n}$ such that $n_1\leq n_2$.  Recall that
${\cal Q}$ is the space of rational functions in
$x_1,...,x_n,y_1,...,y_n$. Define a representation of
$sl(n,\mbb{F})$ on ${\cal Q}$ determined by
$$E_{i,j}|_{\cal
Q}=E_{i,j}^x-E_{j,i}^y\qquad\for\;\;i,j\in\ol{1,n}\eqno(3.1)$$ with
$$E_{i,j}^x|_{\cal Q}=\left\{\begin{array}{ll}-x_j\ptl_{x_i}-\delta_{i,j}&\mbox{if}\;
i,j\in\ol{1,n_1};\\ \ptl_{x_i}\ptl_{x_j}&\mbox{if}\;i\in\ol{1,n_1},\;j\in\ol{n_1+1,n};\\
-x_ix_j &\mbox{if}\;i\in\ol{n_1+1,n},\;j\in\ol{1,n_1};\\
x_i\partial_{x_j}&\mbox{if}\;i,j\in\ol{n_1+1,n}
\end{array}\right.\eqno(3.2)$$
and
$$E_{i,j}^y|_{\cal Q}=\left\{\begin{array}{ll}y_i\ptl_{y_j}&\mbox{if}\;
i,j\in\ol{1,n_2};\\ -y_iy_j&\mbox{if}\;i\in\ol{1,n_2},\;j\in\ol{n_2+1,n};\\
\ptl_{y_i}\ptl_{y_j} &\mbox{if}\;i\in\ol{n_2+1,n},\;j\in\ol{1,n_2};\\
-y_j\partial_{y_i}-\delta_{i,j}&\mbox{if}\;i,j\in\ol{n_2+1,n}.
\end{array}\right.\eqno(3.3)$$
Recall $\flat$ in (2.13) and define
$$\flat'=\sum_{i=1}^{n_2}y_i\ptl{y_i}-\sum_{r=n_2+1}^ny_r\ptl{y_r}.\eqno(3.4)$$
 Moreover, the deformed Laplace operator
${\cal D}$ in (1.9) and its dual $\eta$ in (1.10). Then
$$TE_{i,j}|_{\cal Q}=E_{i,j}|_{\cal Q}T\qquad\for\;\;T=\flat,\flat',{\cal
D},\eta;\;i,j\in\ol{1,n}.\eqno(3.5)$$ In addition,
$$[\flat,{\cal D}]=[\flat',{\cal D}]=-{\cal D},\qquad
[\flat,\eta]=[\flat',\eta]=\eta.\eqno(3.6)$$

By (3.1)-(3.3), we find
$$E_{i,r}|_{\cal
Q}=-x_r\ptl_{x_i}-y_r\ptl_{y_i}\qquad\for\;\;1\leq i<r\leq
n_1,\eqno(3.7)$$
$$E_{i,n_1+s}|_{\cal
Q}=\ptl_{x_i}\ptl_{x_{n_1+s}}-y_{n_1+s}\ptl_{y_i}\qquad\for\;\;i\in\ol{1,n_1},\;s\in\ol{1,n_2-n_1},
\eqno(3.8)$$
$$E_{r,s}|_{\cal Q}=x_r\ptl_{x_s}-y_s\ptl_{y_r}\qquad\for\;\;n_1<
r<s\leq n_2,\eqno(3.9)$$
$$E_{n_2,n_2+1}=x_{n_2}\ptl_{x_{n_2+1}}-\ptl_{y_{n_2}}\ptl_{y_{n_2+1}},\eqno(3.10)$$
$$E_{i,r}|_{\cal
Q}=x_i\ptl_{x_r}+y_i\ptl_{y_r}\qquad\for\;\;n_2+1\leq i<r\leq
n.\eqno(3.11)$$ The subalgebra $sl(n,\mbb{F})_+$ in (2.30) is
generated by the above $E_{i,j}$.

 Denote
$$\zeta_1=x_{n_1-1}y_{n_1}-x_{n_1}y_{n_1-1},\;\;
\zeta=\sum_{r=n_1+1}^{n_2}x_ry_r,\;\;\zeta_2=x_{n_2+1}y_{n_2+2}-x_{n_2+2}y_{n_2+1}.
\eqno(3.12)$$ We will process according to three cases.\psp

{\it Case 1}. $n_1+1<n_2$ \psp

 Assume $n_1+1<n_2<n$. Suppose that $f\in{\cal Q}$ is a singular vector. By Lemma 2.4,
$f$ can be written as a rational function in
$$\{x_i,y_r,\zeta_1,\zeta,\zeta_2\mid
n_2+2\neq i\in\ol{1,n_1+1}\bigcup\ol{n_2+1,n},\;n_1-1\neq
r\in\ol{1,n_1}\bigcup\ol{n_2,n}\}.\eqno(3.13)$$ Note
$$E_{n_1-1,n_1}(f)=-x_{n_1}\ptl_{x_{n_1-1}}(f)=0\eqno(3.14)$$
by (3.7) and
$$E_{n_2+1,n_2+2}(f)=y_{n_2+1}\ptl_{y_{n_2+2}}(f)=0\eqno(3.15)$$
by (3.11). So $f$ is independent of $x_{n_1-1}$ and $y_{n_2+2}$. For
$i\in\ol{1,n_1-2}$, we have
\begin{eqnarray*}\qquad E_{i,n_1-1}(f)&=&-x_{n_1-1}\ptl_{x_i}(f)-y_{n_1-1}\ptl_{y_i}(f)\\
&=& -x_{n_1-1}(\ptl_{x_i}(f)+x_{n_1}^{-1}y_{n_1}\ptl_{y_i}(f))
+x_{n_1}^{-1}\zeta_1\ptl_{y_i}(f)=0\hspace{2.5cm}(3.16)\end{eqnarray*}
 by (3.7). Since both
 $\ptl_{x_i}(f)+x_{n_1}^{-1}y_{n_1}\ptl_{y_i}(f)$ and
 $x_{n_1}^{-1}\zeta_1\ptl_{y_i}(f)$ are independent of $x_{n_1-1}$,
 we have $\ptl_{y_i}(f)=0$, which implies $\ptl_{x_i}(f)$=0 by (3.16). Thus $f$ is
 independent of $\{x_i,y_i\mid i\in\ol{1,n_1-1}$. Similarly, we can
 prove that $f$ is
 independent of $\{x_i,y_i\mid i\in\ol{n_2+1,n}$. Therefore, $f$
 only depends on
 $$\{x_{n_1},x_{n_1+1},x_{n_2+1},y_{n_1},y_{n_2},y_{n_2+1},\zeta_1,\zeta,\zeta_2\}.
 \eqno(3.17)$$

According to (3.8) and (3.12), $E_{n_1,n_1+1}|_{\cal
Q}=\ptl_{x_{n_1}}\ptl_{x_{n_1+1}}-y_{n_1+1}\ptl_{y_{n_1}}$ and
$$E_{n_1,n_1+1}(f)=f_{x_{n_1}x_{n_1+1}}+y_{n_1+1}(f_{x_{n_1}\zeta}-y_{n_1-1}f_{\zeta_1\zeta}-f_{y_{n_1}}-x_{n_1-1}f_{\zeta_1})=0. \eqno(3.18)$$ Applying
$E_{n_1+1,n_2}|_{\cal
Q}=x_{n_1+1}\ptl_{x_{n_2}}-y_{n_2}\ptl_{y_{n_1+1}}$ to the above
equation, we get
$$-f_{x_{n_1}\zeta}+y_{n_1-1}f_{\zeta_1\zeta}+f_{y_{n_1}}+x_{n_1-1}f_{\zeta_1}=0\eqno(3.19)$$
by (3.9). According to (3.12),
$$x_{n_1-1}=y_{n_1}^{-1}\zeta_1+x_{n_1}y_{n_1}^{-1}y_{n_1-1}.\eqno(3.20)$$
Substituting it into (3.19), we get
$$y_{n_1-1}(f_{\zeta_1\zeta}+y_{n_1}^{-1}x_{n_1}f_{\zeta_1})+f_{y_{n_1}}
+y_{n_1}^{-1}\zeta_1f_{\zeta_1}-f_{x_{n_1}\zeta}=0.\eqno(3.21)$$
 Since $f$ is independent of $y_{n_1-1}$, we have
$$f_{\zeta_1\zeta}+y_{n_1}^{-1}x_{n_1}f_{\zeta_1}=0.\eqno(3.22)$$
Thus
$$f_{\zeta_1}=e^{-y_{n_1}^{-1}x_{n_1}\zeta}g\eqno(3.23)$$
for some function $g$ in the variables of (3.17) except $\zeta$,
i.e., $g_{\zeta}=0$. But $f$ is a rational function in the variables
of (3.17) and so is $f_{\zeta_1}$. Hence (3.23) forces
$f_{\zeta_1}=0$, that is, $f$ is independent of $\zeta_1$.
Similarly, we can prove that $f$ is independent of $\zeta_2$. Now
$f$ only depends on
$$\{x_{n_1},x_{n_1+1},x_{n_2+1},y_{n_1},y_{n_2},y_{n_2+1},\zeta\}.\eqno(3.24)$$
Since $\zeta=\sum_{i=n_1+1}^{n_2}x_iy_i$, $f\in{\cal
B}=\mbb{F}[x_1,...,x_n,y_1,...,y_n]$ if and only if $f$ is a
polynomial in the variables (3.24).  Now (3.18) and (3.19) are
equivalent to
$$f_{x_{n_1}x_{n_1+1}}=0,\qquad
f_{x_{n_1}\zeta}-f_{y_{n_1}}=0.\eqno(3.25)$$ Similarly, we can prove
$$f_{y_{n_2}y_{n_2+1}}=0,\qquad
f_{y_{n_2+1}\zeta}-f_{x_{n_2+1}}=0.\eqno(3.26)$$ Set
$$\phi(m_1,m_2)=\sum_{i=0}^{\infty}\frac{y_{n_1}^i(\ptl_{x_{n_1}}\ptl_{\zeta})^i(x_{n_1}^{m_1}\zeta^{m_2})}{i!}\qquad\for
\;\;m_1,m_2\in\mbb{N}.\eqno(3.27)$$ By Lemma 2.1 with
$T_1=\ptl_{y_{n_1}},\; T_1^-=\int_{(y_{n_1})}$ (cf. (2.6)) and
$T_2=-\ptl_{x_{n_1}}\ptl_{\zeta}$, the polynomial solution space of
(3.25) is
$$[\mbb{F}[x_{n_1+1},\zeta]+\sum_{m_1=1}^\infty\sum_{m_2=0}^\infty\mbb{F}\phi(m_1,m_2)]
[\mbb{F}[x_{n_2+1},y_{n_2},y_{n_2+1}]].\eqno(3.28)$$ Denote
$$\psi(m_1,m_2)=\sum_{i=0}^{\infty}\frac{x_{n_2+1}^i(\ptl_{y_{n_2+1}}\ptl_{\zeta})^i(y_{n_2+1}^{m_1}\zeta^{m_2})}{i!}\qquad\for
\;\;m_1,m_2\in\mbb{N},\eqno(3.29)$$
\begin{eqnarray*}\qquad\phi(m_1,m_2,m_3)&=&\sum_{r=0}^\infty\frac{x_{n_2+1}^r(\ptl_{y_{n_2+1}}\ptl_{\zeta})^r(\phi(m_1,m_2)y_{n_2+1}^{m_3})}{r!}
\\ &=&\sum_{i,r=0}^\infty\frac{y_{n_i}^ix_{n_2+1}^r\ptl_{x_{n_1}}^i\ptl_{y_{n_2+1}}^r\ptl_{\zeta}^{i+r}
(x_{n_1}^{m_1}\zeta^{m_2}y_{n_2+1}^{m_3})}{i!r!}.\hspace{3cm}(3.30)\end{eqnarray*}

Solving (3.26) by Lemma 2.1 with $T_1=\ptl_{x_{n_2+1}}\;
T_1^-=\int_{(x_{n_2+1})}$ (cf. (2.6)) and
$T_2=-\ptl_{y_{n_2+1}}\ptl_{\zeta}$, we find the polynomial solution
space of the system (3.25) and (3.26) is
\begin{eqnarray*}\qquad& &\mbb{F}[x_{n_1+1},y_{n_2},\zeta]+\sum_{m_1,m_3=1}^\infty\sum_{m_2=0}^\infty \mbb{F}\phi(m_1,m_2,m_3)
\\ &&+\sum_{m_1=1}^\infty\sum_{m_2=0}^\infty(\mbb{F}[y_{n_2}]\phi(m_1,m_2)+
\mbb{F}[x_{n_1+1}]\psi(m_1,m_2)).\hspace{4.3cm}(3.31)\end{eqnarray*}
According to (1.10),
$$x_{n_1+1}^{m_1}y_{n_2}^{m_2}\zeta^{m_3}=\eta^{m_3}(x_{n_1+1}^{m_1}y_{n_2}^{m_2}),\eqno(3.32)$$
$$\eta^{m_2}(x_{n_1}^{m_1}y_{n_2}^{m_3})=
(\zeta+y_{n_1}\ptl_{x_{n_1}})^{m_2}(x_{n_1}^{m_1}y_{n_2}^{m_3})
=\phi(m_1,m_2)y_{n_2}^{m_3},\eqno(3.33)$$
$$\eta^{m_2}(y_{n_2+1}^{m_1}x_{n_1+1}^{m_3})=(\zeta+x_{n_2+1}\ptl_{y_{n_2+1}})^{m_2}
(y_{n_2+1}^{m_1}x_{n_1+1}^{m_3})
=\psi(m_1,m_2)x_{n_1+1}^{m_3},\eqno(3.34)$$
$$\eta^{m_2}(x_{n_1}^{m_1}y_{n_2+1}^{m_3})=(\zeta+y_{n_1}\ptl_{x_{n_1}}+x_{n_2+1}\ptl_{y_{n_2+1}})^{m_2}(x_{n_1}^{m_1}
y_{n_2+1}^{m_3}) =\phi(m_1,m_2,m_3).\eqno(3.35)$$ It can be verified
that $\{\eta^{m_1}(x_i^{m_2}y_j^{m_3})\mid
m_1,m_2,m_3\in\mbb{N};i=n_1,n_1+1;j=n_2,n_2+1\}$ are singular
vectors. By (3.31)-(3.35), the  nonzero vectors in
$$\{\mbb{F}[\eta](x_i^{m_1}y_j^{m_2})\mid
m_1,m_2\in\mbb{N};i=n_1,n_1+1;j=n_2,n_2+1\}\eqno(3.36)$$ are all the
singular vectors of $sl(n,\mbb{F})$ in ${\cal
 B}=\mbb{F}[x_1,...,x_{n_1},y_1,...,y_{n_2}]$.

 Similarly, when $n_2=n$ and $n_1\leq n-2$,
 the  nonzero vectors in
$$\{\mbb{F}[\eta](x_i^{m_1}y_n^{m_2})\mid
m_1,m_2\in\mbb{N};i=n_1,n_1+1\}\eqno(3.37)$$ are all the singular
vectors of $sl(n,\mbb{F})$ in ${\cal B}$.

Denote
$${\cal H}=\{f\in{\cal B}\mid{\cal D}(f)=0\}.\eqno(3.38)$$
By (3.5), ${\cal H}$ forms an $sl(n,\mbb{F})$-submodule.  Recall
${\cal B}_{\la \ell_1,\ell_2\ra}$ defined in (1.11). Then
$${\cal B}_{\la \ell_1,\ell_2\ra}=\{f\in{\cal B}\mid\flat(f)=\ell_1f;\flat'(f)=\ell_2f\}
\eqno(3.39)$$ by (2.13) and (3.4). Moreover, ${\cal
B}=\bigoplus_{\ell_1,\ell_2\in\mbb{Z}}{\cal B}_{\la
\ell_1,\ell_2\ra}$ becomes a $\mbb{Z}^2$-graded algebra.
 According to (3.5), ${\cal B}_{\la \ell_1,\ell_2\ra}$ forms an
$sl(n,\mbb{F})$-submodule, and so does
$${\cal H}_{\la \ell_1,\ell_2\ra}={\cal B}_{\la
\ell_1,\ell_2\ra}\bigcap {\cal H}.\eqno(3.40)$$

Next (1.9) and (1.10) imply
$$[{\cal D},\eta]=n_2-n_1+\flat+\flat',\;\;{\cal
D}(x_i^{m_1}y_j^{m_2})=0\eqno(3.41)$$ for
$m_1,m_2\in\mbb{N},\;i=n_1,n_1+1$ and $j=n_2,n_2+1$. Thus
$$x_{n_1+1}^{m_1}y_{n_2}^{m_2}\in {\cal H}_{\la m_1,m_2\ra},
\qquad x_{n_1+1}^{m_1}y_{n_2+1}^{m_2}\in {\cal H}_{\la
m_1,-m_2\ra},\eqno(3.42)$$
$$x_{n_1}^{m_1}y_{n_2}^{m_2}\in {\cal H}_{\la -m_1,m_2\ra},
\qquad x_{n_1}^{m_1}y_{n_2+1}^{m_2}\in {\cal H}_{\la
-m_1,-m_2\ra}.\eqno(3.43)$$ For any $g\in {\cal
H}_{\la\ell_1,\ell_2\ra}$and $0<m\in\mbb{N}$, we have
$\eta^m(g)\in{\cal B}_{\ell_1+m,\ell_2+m}$ and
$${\cal
D}(\eta^m(g))=m(n_2-n_1+\ell_1+\ell_2+m-1)\eta^{m-1}(g).\eqno(3.44)$$
Thus
$${\cal
D}(\eta^m(g))=0\;\;\mbox{if and only if}\;\;\ell_1+\ell_2\leq
n_1-n_2\;\;\mbox{and}\;\;m=n_1-n_2-\ell_1-\ell_2+1.\eqno(3.45)$$ If
so,
$$\eta^m(g)\in{\cal
H}_{n_1-n_2-\ell_2+1,n_1-n_2-\ell_1+1}.\eqno(3.46)$$ Note
$$(n_1-n_2-\ell_2+1)+(n_1-n_2-\ell_1+1)=n_1-n_2+2+(n_1-n_2-\ell_1-\ell_2)\geq
n_1-n_2+2.\eqno(3.47)$$

Let $f_{\la\ell_1,\ell_2\ra}\in {\cal H}_{\la\ell_1,\ell_2\ra}$ be a
singular vector in (3.42) and (3.43). Then the singular vectors in
${\cal H}$ are nonzero weight vectors in
$$\mbox{Span}\{f_{\la\ell_1,\ell_2\ra},\eta^{n_1-n_2+1-r_1-r_2}(f_{\la
r_1,r_2\ra})\mid\ell_1,\ell_2,r_1,r_2\in\mbb{Z};r_1+r_2\leq
n_1-n_2\}\eqno(3.48)$$ by (3.36), where
$$\eta^{n_1-n_2+1-r_1-r_2}(f_{\la r_1,r_2\ra})\in {\cal H}_{\la
n_1-n_2+1-r_2,n_1-n_2+1-r_1\ra}.\eqno(3.49)$$  Thus when
$n_1+1<n_2<n,$
  we have
$${\cal H}_{\la \ell_1,\ell_2\ra}\;\mbox{has a unique singular
vector if}\;\;\ell_1+\ell_2\leq n_1-n_2+1\eqno(3.50)$$ and
$${\cal H}_{\la \ell_1,\ell_2\ra}\;\mbox{has exactly two singular
vectors if}\;\;\ell_1+\ell_2> n_1-n_2+1.\eqno(3.51)$$ In the case
$n_1+1<n_2=n,$ ${\cal B}_{\la\ell_1,\ell_2\ra}=0$ if $\ell_2<0$, and
for $\ell_1\in\mbb{Z}$, $\ell_2\in\mbb{N}$, $${\cal H}_{\la
\ell_1,\ell_2\ra}\;\mbox{has a unique singular vector
if}\;\;\ell_1\geq n_1-n+2\;\mbox{or}\;\ell_1+\ell_2\leq
n_1-n+1.\eqno(3.52)$$
$${\cal H}_{\la
\ell_1,\ell_2+1\ra}\;\mbox{has exactly two singular vector
if}\;\mbox{and}\;\;n_1-n+1-\ell_2\leq\ell_1\leq
n_1-n+1.\eqno(3.53)$$

Observe that the symmetric bilinear form $(\cdot|\cdot)$ on ${\cal
B}$ is determined by
$$(x^{\al}y^\be|x^{\al_1}y^{\be_1})=\dlt_{\al,\al_1}
\dlt_{\be,\be_1}(-1)^{\sum_{i=1}^{n_1}\al_i+\sum_{r=n_2+1}^n\be_r}\al!\be!\qquad\for
\;\;\al,\be,\al_1,\be_1\in\mbb{N}\:^n.\eqno(3.54)$$ When
$n_1+1<n_2<n$,  Lemma 2.3 tells us that ${\cal H}_{\la
\ell_1,\ell_2\ra}$ for $\ell_1,\ell_2\in\mbb{Z}$ is an irreducible
summand of ${\cal B}_{\ell_1,\ell_2}$ if and only if
$\ell_1+\ell_2\leq n_1-n_2+1$. It can be verified that
$$({\cal D}(x^{\al}y^\be)|x^{\al_1}y^{\be_1})=
(x^{\al}y^\be|\eta(x^{\al_1}y^{\be_1})).\eqno(3.55)$$  Recall that
$f_{\la\ell_1,\ell_2\ra}\in {\cal H}_{\la\ell_1,\ell_2\ra}$ is a
singular vector in (3.42) and (3.43). Thus
$$(f_{\la\ell_1,\ell_2\ra}|f_{\la\ell_1,\ell_2\ra})\neq
0\eqno(3.56)$$ and
$$(f_{\la\ell_1,\ell_2\ra}|f_{\la\ell_1',\ell_2'\ra})=0\qquad\mbox{if}\;\;(\ell_1,\ell_2)\neq
(\ell_1',\ell_2').\eqno(3.57)$$ Recall $sl(n,\mbb{F})_+$ in (2.30)
and let $sl(n,\mbb{F})_-=\sum_{1\leq i<j\leq n}\mbb{F}E_{j,i}$ be
the subalgebra spanned by the negative root vectors. Moreover,
$(sl(n,\mbb{F})_-)^t=sl(n,\mbb{F})_+.$ According to (3.7)-(3.11),
${\cal B}$ is nilpotent with respect to $sl(n,\mbb{F})_+$. Thus all
${\cal H}_{\la\ell_1,\ell_2\ra}$ with $\ell_1+\ell_2\leq n_1-n_2+1$
are irreducible $sl(n,\mbb{F})$-submodules by Lemma 2.3 and (3.50),
and so are $\eta^m({\cal H}_{\la\ell_1,\ell_2\ra})$ for any
$m\in\mbb{N}$ by (3.5).

We extend the transpose to an algebraic anti-isomorphism on
$U(sl(n,\mbb{F}))$ by $1^t=1$ and
$$(A_1A_2\cdots A_r)^t=A_r^t\cdots A_2^tA_1^t\qquad \for\;\; A_i\in sl(n,\mbb{F}).\eqno(3.58)$$  By the irreducibility,
$$ {\cal H}_{\la\ell_1,\ell_2\ra}=U(sl(n,\mbb{F})_-)(f_{\la\ell_1,\ell_2\ra})\qquad\mbox{if}\;\;\ell_1+\ell_2\leq
n_1-n_2+1.\eqno(3.59)$$

Let $\ell_1,\ell_2,\ell_1',\ell_2'\in\mbb{Z}$ such that
$\ell_1+\ell_2,\ell_1'+\ell_2'\leq n_1-n_2+1$ and
$(\ell_1,\ell_2)\neq (\ell_1',\ell_2')$. Then
$$(w(f_{\la\ell_1,\ell_2\ra})|f_{\la\ell_1',\ell_2'\ra})=
(f_{\la\ell_1,\ell_2\ra}|w^t(f_{\la\ell_1',\ell_2'\ra}))=0\qquad\for\;\;w\in
U(sl(n,\mbb{F})_-)sl(n,\mbb{F})_-\eqno(3.60)$$ by Lemma 2.2. Since
$f_{\la\ell_1,\ell_2\ra}$ is a weight vector, we have
$U(H)(f_{\la\ell_1,\ell_2\ra})\subset\mbb{F}f_{\la\ell_1,\ell_2\ra}$
(cf. (2.29)). Thus for any $w_1,w_2\in U(sl(n,\mbb{F})_-)$,
$$(w_1(f_{\la\ell_1,\ell_2\ra})|w_1(f_{\la\ell_1',\ell_2'\ra}))=
(w_2^tw_1(f_{\la\ell_1,\ell_2\ra})|f_{\la\ell_1',\ell_2'\ra})=c(f_{\la\ell_1,\ell_2\ra})|f_{\la\ell_1',\ell_2'\ra})
\eqno(3.61)$$ for some $c\in\mbb{F}$ by (3.60). Hence (3.59) implies
$$({\cal H}_{\la\ell_1,\ell_2\ra}|{\cal
H}_{\la\ell_1',\ell_2'\ra})=\{0\}.\eqno(3.62)$$

For $f\in {\cal H}_{\la\ell_1,\ell_2\ra},\;g\in {\cal B}$ and
$m,m_1\in\mbb{N}$ such that $m\leq m_1$, we find
$$(\eta^m(f)|\eta^{m_1}(g))=({\cal D}^{m_1}\eta^m(f)|g)
=\dlt_{m_1,m}m![\prod_{r=0}^{m-1}
(\ell_1+\ell_2+n_2-n_1+r)](f|g)\eqno(3.63)$$by (3.44) and (3.55). In
particular, the singular vectors $\eta^{n_1-n_2+1-r_1-r_2}(f_{\la
r_1,r_2\ra})$ for $r_1,r_2\in\mbb{Z}$ with $r_1+r_2\leq n_1-n_2$ are
isotropic polynomials. Moreover, for $m,m_1\in\mbb{N}$ and
$\ell_1,\ell_2,\ell_1',\ell_2'\in\mbb{Z}$ such that
$\ell_1+\ell_2,\ell_1'+\ell_2'\leq n_1-n_2+1$,
$$(\eta^m({\cal H}_{\la\ell_1,\ell_2\ra})|\eta^{m_1}({\cal
H}_{\la\ell_1',\ell_2'\ra}))=\{0\}\qquad\mbox{if}\;\;(m,\ell_1,\ell_1)\neq
(m_1,\ell_1',\ell_1')\eqno(3.64)$$ by (3.62) and (3.63). On the
other hand, $$(\eta^m(f_{\la\ell_1,\ell_2\ra})|\eta^m
(f_{\la\ell_1,\ell_2\ra}))=m![\prod_{r=0}^{m-1}
(\ell_1+\ell_2+n_2-n_1+r)](f_{\la\ell_1,\ell_2\ra}|
f_{\la\ell_1,\ell_2\ra})\neq 0\eqno(3.65)$$ by (3.63). Since the
radical of $(\cdot|\cdot)$ on $\eta^m({\cal
H}_{\la\ell_1,\ell_2\ra})$ is a proper submodule by Lemma 2.2, the
irreducibility of $\eta^m({\cal H}_{\la\ell_1,\ell_2\ra})$  implies
that
$$(\cdot|\cdot)\;\;\mbox{is nondegenerate rewtricted to}\;\;\eta^m({\cal
H}_{\la\ell_1,\ell_2\ra}).\eqno(3.66)$$

Fix $\ell_1,\ell_2\in\mbb{Z}$ with $\ell_1+\ell_2\leq n_1-n_2+1$.
Set
$$\hat{\cal B}_{\la\ell_1,\ell_2\ra}=\sum_{m=0}^\infty\eta^m({\cal
H}_{\la\ell_1-m,\ell_2-m\ra}).\eqno(3.67)$$ By (3.64) and (3.66),
the above sum is a direct sum and $(\cdot|\cdot)$ is nondegenerate
restricted to $\hat{\cal B}_{\la\ell_1,\ell_2\ra}$. Hence
$${\cal B}_{\la\ell_1,\ell_2\ra}=\hat{\cal
B}_{\la\ell_1,\ell_2\ra}\oplus (\hat{\cal
B}_{\la\ell_1,\ell_2\ra}^\perp\bigcap{\cal
B}_{\la\ell_1,\ell_2\ra}).\eqno(3.68)$$ If $\hat{\cal
B}_{\la\ell_1,\ell_2\ra}^\perp\bigcap{\cal
B}_{\la\ell_1,\ell_2\ra}\neq \{0\}$, then it contains a singular
vector, which must be of the form
$\eta^{m_1}(f_{\la\ell_1-m_1,\ell_2-m_1\ra})$ for some
$m_1\in\mbb{N}$ by (3.36). This contradicts (3.65). Thus $\hat{\cal
B}_{\la\ell_1,\ell_2\ra}^\perp\bigcap{\cal B}_{\la\ell_1,\ell_2\ra}=
\{0\}$, equivalently
$${\cal B}_{\la\ell_1,\ell_2\ra}=\bigoplus_{m=0}^\infty\eta^m({\cal
H}_{\la\ell_1-m,\ell_2-m\ra})\eqno(3.69)$$ is completely reducible.
Applying (3.69) to ${\cal B}_{\la\ell_1-1,\ell_2-1\ra}$, we have
$${\cal B}_{\ell_1,\ell_2}={\cal H}_{\la\ell_1,\ell_2\ra}\oplus
\eta({\cal B}_{\la\ell_1-1,\ell_2-1\ra})\qquad\mbox{if}\;\;
\ell_1+\ell_2\leq n_1-n_2+1.\eqno(3.70)$$

Assume $n_1+1<n_2=n$. For $\ell_1\in\mbb{Z}$ and $\ell_2\in\mbb{N}$
such that $\ell_1\geq n_1-n+2\;\mbox{or}\;\ell_1+\ell_2\leq
n_1-n+1$, all  ${\cal H}_{\ell_1,\ell_2}$  are irreducible
submodules of ${\cal B}_{\ell_1,\ell_2}$ by Lemma 2.3, (3.52) and
(3.54). When $\ell_1+\ell_2\leq n_1-n_2+1$, (3.64), (3.66), (3.69)
and (3.70) also hold by the same arguments as in the above. In
summary, we have:\psp

{\bf Theorem 3.1}. {\it Suppose $n_1+1<n_2$. For
$\ell_1,\ell_2\in\mbb{Z}$ with $\ell_1+\ell_2\leq n_1-n_2+1$ ,
${\cal H}_{\la\ell_1,\ell_2\ra}$ is an irreducible highest-weight
$sl(n,\mbb{F})$-module and
$${\cal B}_{\la\ell_1,\ell_2\ra}=\bigoplus_{m=0}^\infty\eta^m({\cal
H}_{\la\ell_1-m,\ell_2-m\ra})\eqno(3.71)$$ is an orthogonal
decomposition of irreducible submodules. In particular, ${\cal
B}_{\la\ell_1,\ell_2\ra}={\cal H}_{\la\ell_1,\ell_2\ra}\oplus
\eta({\cal B}_{\la\ell_1-1,\ell_2-1\ra})$. The symmetric bilinear
form $(\cdot|\cdot)$ restricted to $\eta^m({\cal
H}_{\la\ell_1-m,\ell_2-m\ra})$ is nondengerate. If $n_2<n$, all
${\cal H}_{\la\ell_1,\ell_2\ra}$ for $\ell_1,\ell_2\in\mbb{Z}$ with
$\ell_1+\ell_2> n_1-n_2+1$ have exactly two singular vectors.

 Assume $n_2=n$. Then ${\cal
B}_{\la\ell,0\ra}={\cal H}_{\la\ell,0\ra}$ with $\ell\in\mbb{Z}$ are
irreducible highest-weight $sl(n,\mbb{F})$-modules.  All ${\cal
H}_{\la\ell_1,\ell_2\ra}$ for $\ell_1\in\mbb{Z}$ and
$\ell_2\in\mbb{N}$ such that $\ell_1\geq n_1-n+2$ are also
irreducible highest-weight $sl(n,\mbb{F})$-modules. Moreover, for
$\ell_2\in 1+\mbb{N}$ and $n_1-n_2+1+\ell_2\leq \ell_1\in\mbb{Z}$,
the orthogonal decompositions in (3.71) also holds. Furthermore,
${\cal H}_{\la\ell_1,\ell_2+1\ra}$ for $\ell_1\in\mbb{Z}$ and
$\ell_2\in\mbb{N}$ such that $n_1-n+1-\ell_2\leq \ell_1\leq n_1-n+1$
have exactly two singular vectors.}\psp

Indeed, we have more detailed information.  Suppose $n_1+1<n_2<n$.
For $m_1,m_2\in\mbb{N}$ with $m_1+m_2\geq n_2-n_1-1$, ${\cal
H}_{\la-m_1,-m_2\ra}$ has  a highest-weight vector
$x_{n_1}^{m_1}y_{n_2+1}^{m_2}$ of weight
$m_1\lmd_{n_1-1}-(m_1+1)\lmd_{n_1}-(m_2+1)\lmd_{n_2}+m_2(1-\dlt_{n_2,n-1})\lmd_{n_2+1}$.
When $m_1,m_2\in\mbb{N}$ with $m_2-m_1\geq n_2-n_1-1$, ${\cal
H}_{\la m_1,-m_2\ra}$ has a highest-weight vector
$x_{n_1+1}^{m_1}y_{n_2+1}^{m_2}$ of weight
$-(m_1+1)\lmd_{n_1}+m_1\lmd_{n_1+1}-(m_2+1)\lmd_{n_2}+m_2(1-\dlt_{n_2,n-1})\lmd_{n_2+1}$.
If $m_1,m_2\in\mbb{N}$ with $m_1-m_2\geq n_2-n_1-1$, ${\cal H}_{\la
-m_1,m_2\ra}$ is has a highest-weight vector
$x_{n_1}^{m_1}y_{n_2}^{m_2}$ of weight
$m_1\lmd_{n_1-1}-(m_1+1)\lmd_{n_1}+m_2\lmd_{n_2-1}-(m_2+1)\lmd_{n_2}$.

Assume $n_1+1<n_2=n$. When $m_1,m_2\in\mbb{N}$,  ${\cal H}_{\la
m_1,m_2\ra}$  has a highest-weight vector $x_{n_1+1}^{m_1}y_n^{m_2}$
of weight $-(m_1+1)\lmd_{n_1}+m_1\lmd_{n_1+1}+m_2\lmd_{n-1}$. If
$m_1,m_2\in\mbb{N}$ with $m_1\leq n-n_1-2$ or $m_2-m_1\leq n_1-n+1$,
${\cal H}_{\la -m_1,m_2\ra}$ has a highest-weight vector
$x_{n_1}^{m_1}y_n^{m_2}$ of weight
$m_1\lmd_{n_1-1}-(m_1+1)\lmd_{n_1}+m_2\lmd_{n-1}$.

By Lemma 2.1 with
$T_1=\ptl_{x_{n_1+1}}\ptl_{y_{n_1+1}},\;T_1^-=\int_{(x_{n_1+1})}\int_{(y_{n_1+1})}$
and $T_2={\cal D}-\ptl_{x_{n_1}+1}\ptl_{y_{n_1}+1}$, ${\cal
H}_{\la\ell_1,\ell_2\ra}$ has a basis
\begin{eqnarray*}\qquad& &\big\{\sum_{i=0}^\infty\frac{(x_{n_1+1}y_{n_1+1})^i({\cal
D}-\ptl_{x_{n_1}+1}\ptl_{y_{n_1}+1})^i(x^\al
y^\be)}{\prod_{r=1}^i(\al_{n_1+1}+r)(\be_{n_1+1}+r)}
\mid\al,\be\in\mbb{N}\:^n;\\ &
&\al_{n_1+1}\be_{n_1+1}=0;\sum_{r=n_1+1}^n\al_r-\sum_{i=1}^{n_1}\al_i=\ell_1;
\sum_{i=1}^{n_2}\be_i-\sum_{r=n_2+1}^n\be_r=\ell_2
\big\}.\hspace{1.8cm}(3.72)\end{eqnarray*} \psp

{\it Case 2}. $n_1+1=n_2$ \psp

In this case, $\zeta=x_{n_1+1}y_{n_1+1}$. First we consider the
subcase $n_2<n$. Suppose that $f\in{\cal Q}$ is a singular vector.
According to the arguments in (3.13)-(3.17), $f$ is a rational
function in
$$\{x_{n_1},x_{n_1+1},x_{n_1+2},y_{n_1},y_{n_1+2},\zeta,\zeta_1,\zeta_2\}.\eqno(3.73)$$
Moreover, (3.18) holds. Substituting (3.20) and
$y_{n_1+1}=x_{n_1+1}^{-1}\zeta$ into (3.18), we still get (3.22),
which implies $f_{\zeta_1}=0$. Symmetrically, $f_{\zeta_2}=0$. Hence
we can rewrite $f$ as a rational function in
$$\{x_{n_1},x_{n_1+1},x_{n_1+2},y_{n_1},y_{n_1+1},y_{n_1+2}\}.
 \eqno(3.74)$$
Now $f$ is a singular vector if and only if it is a weight vector
satisfying the following system of differential equations
$$(\ptl_{x_{n_1}}\ptl_{x_{n_1+1}}-y_{n_1+1}\ptl_{y_{n_1}})(f)=0,\eqno(3.75)$$
$$(x_{n_1+1}\ptl_{x_{n_1+2}}-\ptl_{y_{n_1+1}}\ptl_{y_{n_1+2}})(f)=0
\eqno(3.76)$$ by (3.8) and (3.10) with $n_2=n_1+1$. Note
$$E_{n_1,n_1+2}|_{\cal Q}=[E_{n_1,n_1+1}|_{\cal Q},E_{n_1+1,n_1+2}|_{\cal
Q}]=
\ptl_{x_{n_1}}\ptl_{x_{n_1+2}}-\ptl_{y_{n_1}}\ptl_{y_{n_1+2}}\eqno(3.77)$$
by (3.8) and (3.10) with $n_2=n_1+1$. So
$$(\ptl_{x_{n_1}}\ptl_{x_{n_1+2}}-\ptl_{y_{n_1}}\ptl_{y_{n_1+2}})(f)=0.\eqno(3.78)$$

For our purpose of representation, we only consider $f$ is a
polynomial in $\{x_i,y_i\mid i=n_1,n_1+1,n_1+2\}$. Set
\begin{eqnarray*}\qquad\phi_{m_1,m_2,m_3}&=&[\prod_{s=1}^{m_2}(m_1+s)]\sum_{i=0}^\infty
\frac{x_{n_1}^{m_1+i}x_{n_1+2}^i(\ptl_{y_{n_1}}\ptl_{y_{n_1+2}})^i
(y_{n_1}^{m_2}y_{n_1+2}^{m_3})} {i!\prod_{r=1}^i(m_1+r)}\\
&=&(y_{n_1}\ptl_{x_{n_1}}+x_{n_1+2}\ptl_{y_{n_1+2}})^{m_2}
(x_{n_1}^{m_1+m_2}y_{n_1+2}^{m_3})\hspace{4.3cm}(3.79)\end{eqnarray*}
and
\begin{eqnarray*}\qquad\psi_{m_1,m_2,m_3}&=&\frac{(m_1+m_2)!\prod_{s=1}^{m_1}(m_3+s)}{m_1!}
\sum_{i=0}^\infty
\frac{x_{n_1}^ix_{n_1+2}^{m_1+i}(\ptl_{y_{n_1}}\ptl_{y_{n_1+2}})^i
(y_{n_1}^{m_2}y_{n_1+2}^{m_3})} {i!\prod_{r=1}^i(m_1+r)}
\\&=&\sum_{i=0}^{m_2}{m_2\choose i}
\frac{(m_1+m_2)!x_{n_1}^ix_{n_1+2}^{m_1+i}
y_{n_1}^{m_2-i}\ptl_{y_{n_1+2}}^{m_1+i}(y_{n_1+2}^{m_1+m_3})}
{(m_1+i)!}
\\&=&\sum_{i=0}^{m_2}{m_2\choose m_2-i}
\frac{(m_1+m_2)!x_{n_1}^i
y_{n_1}^{m_2-i}(x_{n_1+2}\ptl{y_{n_1+2}})^{m_1+i}(y_{n_1+2}^{m_1+m_3})}
{(m_1+i)!}
\\&=&\sum_{i=0}^{m_2}
\frac{(m_1+m_2)!(y_{n_1}\ptl_{x_{n_1}})^{m_2-i}(x_{n_1}^{m_2})
(x_{n_1+2}\ptl{y_{n_1+2}})^{m_1+i}(y_{n_1+2}^{m_1+m_3})}
{(m_2-i)!(m_1+i)!}\\&=&\sum_{r=0}^\infty
\frac{(m_1+m_2)!(y_{n_1}\ptl_{x_{n_1}})^r(x_{n_1}^{m_2})
(x_{n_1+2}\ptl{y_{n_1+2}})^{m_1+m_2-r}(y_{n_1+2}^{m_1+m_3})}
{r!(m_1+m_2-r)!}\\
&=&(y_{n_1}\ptl_{x_{n_1}}+x_{n_1+2}\ptl_{y_{n_1+2}})^{m_1+m_2}
(x_{n_1}^{m_2}y_{n_1+2}^{m_1+m_3}).\hspace{3.8cm}(3.80)\end{eqnarray*}
By Lemma 2.1 with
$T_1=\ptl_{x_{n_1}}\ptl_{x_{n_1+2}},\;T_1^-=\int_{(x_{n_1})}\int_{(x_{n_1+2})}$
(cf. (2.6)) and $T_2=\ptl_{y_{n_1}}\ptl_{y_{n_1+2}}$, the polynomial
solution space of (3.78) is
$$\mbox{Span}\{\phi_{m_1,m_2,m_3}x_{n_1+1}^{m_4}y_{n_1+1}^{m_5},\psi_{m_1+1,m_2,m_3}x_{n_1+1}^{m_4}y_{n_1+1}^{m_5}\mid
m_i\in\mbb{N}\}.\eqno(3.81)$$

Note
$$\phi_{m_1,m_2,0}x_{n_1+1}^{m_3}y_{n_1+1}^{m_4}=[\prod_{r=1}^{m_2}(m_1+r)]x_{n_1}^{m_1}y_{n_1}^{m_2}
x_{n_1+1}^{m_3}y_{n_1+1}^{m_4}\qquad\for\;\;m_i\in\mbb{N}\eqno(3.82)$$
and
$$x_{n_1+1}^{m_1}y_{n_1+1}^{m_2}\psi_{m_3,0,m_4}
=[\prod_{i=1}^{m_3}(m_4+i)]x_{n_1+1}^{m_1}y_{n_1+1}^{m_2}x_{n_1+2}^{m_3}y_{n_1+2}^{m_4}\qquad\for\;\;m_i\in\mbb{N}.
\eqno(3.83)$$ In particular, all the polynomials in (3.83) are
solutions of the equation (3.75). Now
$$\eta=\sum_{i=1}^{n_1}y_i\ptl_{x_i}+x_{n_1+1}y_{n_1+1}+\sum_{s=n_1+2}^n
x_s\ptl_{y_s}.\eqno(3.84)$$
 Write
\begin{eqnarray*} \qquad h_{m_1,m_2,m_3}&=&
\frac{(m_1+m_2)!}{m_1!}
\sum_{i=0}^\infty\frac{x_{n_1}^ix_{n_1+1}^{m_1+i}y_{n_1+1}^{m_3+i}\ptl_{y_{n_1}}^i(y_{n_1}^{m_2})}
{i!\prod_{r=1}^i(m_1+r)}\\&=&\sum_{i=0}^{m_2}{m_2\choose
i}\frac{(m_1+m_2)!x_{n_1}^ix_{n_1+1}^{m_1+i}y_{n_1+1}^{m_3+i}y_{n_1}^{m_2-i}}
{(m_1+i)!}
\\
&=&\sum_{i=0}^{m_2}{m_2\choose
m_2-i}\frac{(m_1+m_2)!x_{n_1}^ix_{n_1+1}^{m_1+i}y_{n_1+1}^{m_3+i}y_{n_1}^{m_2-i}}
{(m_1+i)!}
\\ &=&\sum_{i=0}^{m_2}\frac{(m_1+m_2)!(y_{n_1}\ptl_{x_{n_1}})^{m_2-i}(x_{n_1}^{m_2})x_{n_1+1}^{m_1+i}y_{n_1+1}^{m_3+i}}
{(m_2-i)!(m_1+i)!}\\ &=&
\eta^{m_1+m_2}(x_{n_1}^{m_2}y_{n_1+1}^{m_3-m_1})\hspace{7.8cm}(3.85)\end{eqnarray*}
and calculate
\begin{eqnarray*}\qquad& &\eta^{m_2}(x_{n_1}^{m_1+m_2})y_{n_1+1}^{m_3}
=\eta^{m_2}(x_{n_1}^{m_1+m_2}y_{n_1+1}^{m_3})
\\ &=&\sum_{i=0}^{m_2}{m_2\choose
m_2-i}[\prod_{r=1}^{m_2-i}(m_1+i+r)]y_{n_1}^{m_2-i}x_{n_1}^{m_1+i}x_{n_1+1}^i
y_{n_1+1}^{m_3+i}\\&=&\sum_{i=0}^{m_2}{m_2\choose
i}\frac{[\prod_{r=1}^{m_2}(m_1+r)]y_{n_1}^{m_2-i}x_{n_1}^{m_1+i}x_{n_1+1}^i
y_{n_1+1}^{m_3+i}}{\prod_{s=1}^i(m_1+s)}\\&=&[\prod_{r=1}^{m_2}(m_1+r)]
\sum_{i=0}^{m_2}\frac{x_{n_1}^{m_1+i}x_{n_1+1}^i
y_{n_1+1}^{m_3+i}\ptl_{y_{n_1}}^i(y_{n_1}^{m_2})}{i![\prod_{s=1}^i(m_1+s)]}.
\hspace{5.2cm}(3.86)
\end{eqnarray*}
Lemma 2.1 with
$T_1=\ptl_{x_{n_1}}\ptl_{x_{n_1+1}},\;T_1^-=\int_{(x_{n_1})}\int_{(x_{n_1+1})}$
(cf. (2.6)) and $T_2=-y_{n_1+1}\ptl_{y_{n_1}}$ tells us that
$\mbox{Span}\{h_{m_1,m_2,m_3},\eta^{m_2}(x_{n_1}^{m_1+m_2})y_{n_1+1}^{m_3}\mid
m_i\in\mbb{N}\}$ is the solution space of (3.75) in
$\mbox{Span}\{x_{n_1}^{m_1}y_{n_1}^{m_2}x_{n_1+1}^{m_3}y_{n_1+1}^{m_4}\mid
m_i\in\mbb{N}\}.$ In particular, (3.85) and (3.86) can be viewed as
algorithms of solving the equation (3.75).

On the other hand,
$$\ptl_{x_{n_1}}(\phi_{0,m_2,m_3})=m_2\psi_{1,m_2,m_3-1},\eqno(3.87)$$
$$\ptl_{x_{n_1}}(\phi_{m_1,m_2,m_3})=(m_1+m_2)\phi_{m_1-1,m_2,m_3}\qquad\mbox{if}\;m_1>0,\eqno(3.88)$$
$$
\ptl_{x_{n_1}}(\psi_{m_1,m_2,m_3})=m_2\psi_{m_1+1,m_2-1,m_3-1},\eqno(3.89)$$
$$\ptl_{y_{n_1}}(\phi_{m_1,m_2,m_3})=m_2(m_1+m_2)\phi_{m_1,m_2-1,m_3},\eqno(3.90)$$
$$\ptl_{y_{n_1}}(\psi_{m_1,m_2,m_3})=m_2(m_1+m_2)\psi_{m_1,m_2-1,m_3}.\eqno(3.91)$$
Applying the algorithm (3.85) to (3.81), we get that
$\hat\phi_{m_1,m_2,m_3}y_{n+1}^{m_4}$ and
$\hat\psi_{m_1,m_2,m_3}y_{n+1}^{m_4}$ are the solutions of (3.75) by
(3.87)-(3.91), where
$$\hat\phi_{m_1,m_2,m_3}=
\sum_{i=0}^\infty{m_2\choose
i}\phi_{m_1+i,m_2-i,m_3}(x_{n_1+1}y_{n_1+1})^i=\eta^{m_2}(x_{n_1}^{m_1+m_2}y_{n_1+2}^{m_3}),
\eqno(3.92)$$
\begin{eqnarray*}\qquad\hat\psi_{m_1,m_2,m_3}&=&
\sum_{i=0}^{m_1}{m_1+m_2\choose
i}\psi_{m_1-i,m_2,m_3+i}(x_{n_1+1}y_{n_1+1})^i\\
& &+\sum_{r=1}^{m_2}{m_1+m_2\choose
m_1+r}\phi_{m_2-r,r,m_1+m_3}(x_{n_1+1}y_{n_1+1})^{m_1+r}\\&=&\eta^{m_1+m_2}(x_{n_1}^{m_2}y_{n_1+2}^{m_1+m_3}).
\hspace{7.7cm}(3.93)\end{eqnarray*} Using the algorithm (3.86), we
find that the solution space of (3.75) in (3.81) is
\begin{eqnarray*}\qquad \;\;\;& &\mbox{Span}\{x_{n_1+1}^{m_1}y_{n_1+1}^{m_2}x^{m_3}_{n_1+2}y^{m_4}_{n_1+2},
\hat\phi_{m_1,m_2,m_3}y_{n_1+1}^{m_4},h_{m_1,m_2,m_3},\\
& &\qquad\;\;\hat\psi_{m_1+1,m_2+1,m_3}y_{n_1+1}^{m_4}\mid
m_i\in\mbb{N}\}.\hspace{6.6cm}(3.94)\end{eqnarray*}

 According to (3.84), (3.92) and (3.93),
$$\ptl_{x_{n_1+2}}(\hat\phi_{m_1,m_2,m_3})=m_2m_3\hat\phi_{m_1,m_2-1,m_3-1},\eqno(3.95)$$
$$\ptl_{y_{n_1+1}}(\hat\phi_{m_1,m_2,m_3})=m_2x_{n_1+1}\hat\phi_{m_1,m_2-1,m_3},\eqno(3.96)$$
$$\ptl_{y_{n_1+2}}(\hat\phi_{m_1,m_2,m_3})=m_3\hat\phi_{m_1,m_2,m_3-1},\eqno(3.97)$$
$$\ptl_{x_{n_1+2}}(\hat\psi_{m_1,m_2,m_3})=(m_1+m_2)(m_1+m_3)\hat\psi_{m_1-1,m_2,m_3},
\eqno(3.98)$$
$$\ptl_{y_{n_1+1}}(\hat\psi_{m_1,m_2,m_3})=(m_1+m_2)x_{n_1+1}\hat\psi_{m_1-1,m_2,m_3+1},
\eqno(3.99)$$
$$\ptl_{y_{n_1+2}}(\hat\psi_{m_1,m_2,m_3})=(m_1+m_3)\hat\psi_{m_1,m_2,m_3-1}.
\eqno(3.100)$$

Put \begin{eqnarray*} \qquad g_{m_1,m_2,m_3}&=&
\frac{(m_1+m_2)!}{m_1!}
\sum_{i=0}^\infty\frac{y_{n_1+2}^iy_{n_1+1}^{m_1+i}x_{n_1+1}^{m_3+i}\ptl_{x_{n_1+2}}^i(x_{n_1+2}^{m_2})}
{i!\prod_{r=1}^i(m_1+r)}\\ &=&
\eta^{m_1+m_2}(y_{n_1+2}^{m_2}x_{n_1+1}^{m_3-m_1})\hspace{7.4cm}(3.101)\end{eqnarray*}
and
\begin{eqnarray*}\qquad
g'_{m_1,m_2,m_3}&=&
\sum_{i=0}^\infty\frac{[\prod_{s=1}^{m_2}(m_1+s)]y_{n_1+2}^{m_1+i}y_{n_1+1}^ix_{n_1+1}^{m_3+i}\ptl_{x_{n_1+2}}^i(x_{n_1+2}^{m_2})}
{i!\prod_{r=1}^i(m_1+r)}\\
&=&\eta^{m_2}(y_{n_1+2}^{m_1+m_2}x_{n_1+1}^{m_3}).\hspace{7.8cm}(3.102)\end{eqnarray*}
Symmetrically, $\mbox{ Span}\{g_{m_1,m_2,m_3},g'_{m_1,m_2,m_3}\mid
m_i\in \mbb{N}\}$ is the solution space of (3.76) in
Span$\{x_{n_1+1}^{m_1}y_{n_1+1}^{m_2}x^{m_3}_{n_1+2}y^{m_4}_{n_1+2}\mid
m_i\in\mbb{N}\}$ by Lemma 2.1 with
$T_1=\ptl_{y_{n_1+1}}\ptl_{y_{n_1+2}},\;T_1^-=\int_{(y_{n_1+1})}\int_{(y_{n_1+2})}$
(cf. (2.6)) and $T_2=-x_{n_1+1}\ptl_{x_{n_1+2}}$. Observe that
$\{\hat\phi_{m_1,m_2,m_3},\hat\psi_{m_1,m_2,m_3},\\
h_{m_1,m_2,m_3}\mid m_i\in\mbb{N}\}$ are solutions of (3.76). Thus
the solution space of (3.76) in (3.94) is
\begin{eqnarray*}\qquad \;\;\;&
&\mbox{Span}\{g_{m_1,m_2,m_3},g'_{m_1,m_2,m_3},h_{m_1,m_2,m_3},
\hat\phi_{m_1,m_2,m_3},\\
&&\qquad\;\;\hat\phi_{m_1,m_2,0}y_{n_1+1}^{m_3},\hat\psi_{m_1+1,m_2+1,m_3}\mid
m_i\in\mbb{N}\}\hspace{4.8cm}(3.103)\end{eqnarray*} by (3.97) and
(3.100).

Expressions (3.85), (3.92), (3.93), (3.101) and (3.102) imply that
 the solution space of the
singular vectors in ${\cal B}$ is
\begin{eqnarray*} & &\mbox{Span}\{\eta^{m_2}(x_i^{m_1}y_j^{m_3}),
x_{n_1+1}^{m_1}y_{n_1+1}^{m_2},\eta^{m_1+m_2}(x_{n_1}^{m_2}y_{n_1+1}^{m_3-m_1}),\eta^{m_1+m_2}(y_{n_1+2}^{m_2}x_{n_1+1}^{m_3-m_1})
\\ & &\qquad\;\;\mid m_r\in\mbb{N};
(i,j)=(n_1,n_1+1),(n_1,n_1+2),(n_1+1,n_1+2)\}.\hspace{2.1cm}(3.104)\end{eqnarray*}
Remind that in this case, $${\cal
D}=-\sum_{i=1}^{n_1}x_i\ptl_{y_i}+\ptl_{x_{n_1+1}}\ptl_{y_{n_1+1}}-\sum_{s=n_1+2}^n
y_s\ptl_{x_s}.\eqno(3.105)$$ We have
$${\cal
D}[\eta^{m_1+m_2}(x_{n_1}^{m_2}y_{n_1+1}^{m_3-m_1})]=(m_1+m_2)m_3
\eta^{m_1+m_2-1}(x_{n_1}^{m_2}y_{n_1+1}^{m_3-m_1})\eqno(3.106)$$ by
(3.44). Thus we find a singular
$$\eta^{m_1+m_2}(x_{n_1}^{m_2}y_{n_1+1}^{-m_1})\in {\cal H}_{\la
m_1,m_2\ra}\eqno(3.107)$$ of new type if $m_1,m_2\geq 1$.
Symmetrically, $\eta^{m_1+m_2}(y_{n_1+2}^{m_2}x_{n_1+1}^{-m_1})\in
{\cal H}_{\la m_2,m_1\ra}$ is a singular vector.

Recall the singular vectors
$$f_{\la-m_1,-m_2\ra}=x_{n_1}^{m_1}y_{n_1+2}^{m_2}\in{\cal
H}_{\la -m_1,-m_2\ra},\qquad
f_{\la-m_1,m_2\ra}=x_{n_1}^{m_1}y_{n_1+1}^{m_2}\in{\cal H}_{\la
-m_1,m_2\ra},\eqno(3.108)$$
$$f_{\la m_1,-m_2\ra}=x_{n_1+1}^{m_1}y_{n_1+2}^{m_2}\in{\cal
H}_{\la m_1,-m_2\ra}.\eqno(3.109)$$ Moreover, we have the singular
vectors
$$\eta^{-\ell_1-\ell_2}(f_{\la\ell_1,\ell_2\ra})\in {\cal
H}_{\la -\ell_2,-\ell_1\ra}\qquad
\for\;\;\ell_1,\ell_2\in\mbb{Z}\;\mbox{with}\;\ell_1+\ell_2\leq
-1.\eqno(3.110)$$ Therefore, any singular vector in ${\cal H}$ (cf.
(3.38)) is a nonzero weight vector in
\begin{eqnarray*}& & \mbox{Span}\{f_{\la\ell_1,\ell_2\ra},\eta^{-\ell_1'-\ell_2'}
(f_{\la\ell_1',\ell_2'\ra}),\eta^{m_1+m_2}(x_{n_1}^{m_2}y_{n_1+1}^{-m_1}),\eta^{m_1+m_2}(y_{n_1+2}^{m_2}x_{n_1+1}^{-m_1})\\
&
&\qquad\;\;\mid\ell_1,\ell_2,\ell_1',\ell_2'\in\mbb{Z},\;m_1,m_2\in\mbb{N}+1;\ell_1\leq
0\;\mbox{or}\;\ell_2\leq 0;\ell_1'+\ell_2'\leq
-1\}.\hspace{1.2cm}(3.111)\end{eqnarray*}

Assume $n_2=n$. We similarly find that the solution space of the
singular vectors in ${\cal B}$ is
$$\mbox{Span}\{\eta^{m_2}(x_{n-1}^{m_1}y_n^{m_3}),
x_n^{m_1}y_n^{m_2},\eta^{m_1+m_2}(x_{n-1}^{m_2}y_n^{m_3-m_1})\mid
m_i\in\mbb{N}\}.\eqno(3.112)$$ In particular, any singular vector in
${\cal H}$ (cf. (3.38)) is a nonzero weight vector in
\begin{eqnarray*}\qquad& &
\mbox{Span}\{x_{n-1}^{m_1}y_n^{m_2},x_n^{m_1},\eta^{m_1+1}
(x_{n-1}^{m_1+m_2+1}y_n^{m_2}),\\
&&\qquad\;\;\eta^{m_1+m_2+2}(x_{n-1}^{m_2+1}y_n^{-m_1-1})\mid
m_1,m_2\in\mbb{N}\}.\hspace{5.2cm}(3.113)\end{eqnarray*} By the
arguments of (3.55)-(3.70), we have:\psp

{\bf Theorem 3.2}. {\it Suppose $n_1+1=n_2$. For
$\ell_1,\ell_2\in\mbb{Z}$ with $\ell_1+\ell_2\leq 0$ or $n_2=n$ and
$0\leq\ell_2\leq\ell_1$,  ${\cal H}_{\la\ell_1,\ell_2\ra}$ is an
irreducible highest-weight $sl(n,\mbb{F})$-module and
$${\cal B}_{\la\ell_1,\ell_2\ra}=\bigoplus_{m=0}^\infty\eta^m({\cal
H}_{\la\ell_1-m,\ell_2-m\ra})\eqno(3.114)$$ is an orthogonal
decomposition of irreducible submodules. In particular, ${\cal
B}_{\la\ell_1,\ell_2\ra}={\cal H}_{\la\ell_1,\ell_2\ra}\oplus
\eta({\cal B}_{\la\ell_1-1,\ell_2-1\ra})$. The symmetric bilinear
form $(\cdot|\cdot)$ restricted to $\eta^m({\cal
H}_{\la\ell_1-m,\ell_2-m\ra})$ is nondegenerate.

Assume $n_2<n$. For $m_1,m_2\in\mbb{N}+1$, ${\cal H}_{\la
m_1,m_2\ra}$ has exactly three singular vectors. All the submodules
${\cal H}_{\la\ell_1,\ell_2\ra}$ for $\ell_1,\ell_2\in\mbb{Z}$ such
$\ell_1+\ell_2>0$ and $\ell_1\ell_2\leq 0$ have two singular
vectors. Consider $n_2=n$. For $m_1,m_2\in\mbb{N}$ with $m_1<m_2$,
${\cal H}_{\la m_1,m_2\ra}$ is also an irreducible highest-weight
$sl(n,\mbb{F})$-module. All submodules ${\cal H}_{\la
-m_1,m_1+m_2+1\ra}$ with $m_1,m_2\in\mbb{N}$ have have exactly two
singular vectors}.\psp

Indeed, we have more detailed information.  Suppose $n_2<n$. For
$m_1,m_2\in\mbb{N}$, ${\cal H}_{\la-m_1,-m_2\ra}$ has  a
highest-weight vector $x_{n_1}^{m_1}y_{n_1+2}^{m_2}$ of weight
$m_1\lmd_{n_1-1}-(m_1+1)\lmd_{n_1}-(m_2+1)\lmd_{n_1+1}+m_2(1-\dlt_{n_1,n-2})\lmd_{n_1+2}$.
When $m_1,m_2\in\mbb{N}$ with $m_2-m_1\geq 0$, ${\cal H}_{\la
m_1,-m_2\ra}$ has a highest-weight vector
$x_{n_1+1}^{m_1}y_{n_1+2}^{m_2}$ of weight
$-(m_1+1)\lmd_{n_1}+(m_1-m_2-1)\lmd_{n_1+1}+m_2(1-\dlt_{n_1,n-2})\lmd_{n_1+2}$.
If $m_1,m_2\in\mbb{N}$ with $m_1-m_2\geq 0$, ${\cal H}_{\la
-m_1,m_2\ra}$ is has a highest-weight vector
$x_{n_1}^{m_1}y_{n_1+1}^{m_2}$ of weight
$m_1\lmd_{n_1-1}+(m_2-m_1-1)\lmd_{n_1}-(m_2+1)\lmd_{n_1+1}$.

Assume $n_2=n$. For $m_1,m_2\in\mbb{N}$ with $m_2\leq m_1$, ${\cal
H}_{\la -m_1,m_2\ra}$ has a highest-weight vector
$x_{n-1}^{m_1}y_n^{m_2}$ of weight
$m_1\lmd_{n-2}+(m_2-m_1-1)\lmd_{n-1}$. Moreover,
 ${\cal H}_{\la
m,0\ra}$ has a highest-weight vector $x_{n-1}^m$ of weight
$m\lmd_{n-2}-(m+1)\lmd_{n-1}$ for $m\in\mbb{Z}$. For
$m_1,m_2\in\mbb{N}+1$, ${\cal H}_{\la m_1,m_2\ra}$  has a
highest-weight vector $\eta^{m_1+m_2}(x_{n-1}^{m_2}y_n^{-m_1})$ of
weight $m_2\lmd_{n-2}+(m_1-m_2-1)\lmd_{n-1}$. Again ${\cal
H}_{\la\ell_1,\ell_2\ra}$ has a basis of the form (3.72).

\section{The $sl(n,\mbb{F})$-Variation II: $n_1=n_2$}

In this section, we continue the discussion from last section.
Recall $n\geq 2$.\pse

{\it Case 3}. $n_1=n_2$. \psp

In this case, the variated Laplace operator
$${\cal
D}=-\sum_{i=1}^{n_1}x_i\ptl_{y_i}-\sum_{s=n_1+1}^n
y_s\ptl_{x_s}\eqno(4.1)$$ and its dual
$$\eta=\sum_{i=1}^{n_1}y_i\ptl_{x_i}+\sum_{s=n_2+1}^n
x_s\ptl_{y_s}.\eqno(4.2)$$

First we consider the subcase $1<n_1<n-1$. Suppose that $f\in{\cal
Q}$ is a singular vector. According to the arguments in
(3.13)-(3.17), $f$ is a rational function in
$$\{x_{n_1},x_{n_1+1},y_{n_1},y_{n_1+1},\zeta_1,\zeta_2\}
 \eqno(4.3)$$
(cf. (3.12)). Note
$$E_{n_1,n_1+1}|_{\cal
Q}=\ptl_{x_{n_1}}\ptl_{x_{n_1+1}}-\ptl_{y_{n_1}}\ptl_{y_{n_1+1}}\eqno(4.4)$$
by (3.1)-(3.3). Now $E_{n_1,n_1+1}(f)=0$ implies
$$(\ptl_{x_{n_1}}\ptl_{x_{n_1+1}}-\ptl_{y_{n_1}}\ptl_{y_{n_1+1}})(f)=0,\eqno(4.5)$$
equivalently,
\begin{eqnarray*}\qquad& &
(x_{n_1-1}x_{n_1+2}-y_{n_1-1}y_{n_1+2})f_{\zeta_1\zeta_2}-y_{n_1-1}f_{\zeta_1x_{n_1+1}}
-x_{n_1-1}f_{\zeta_1y_{n_1+1}}\\ &
&+y_{n_1+2}f_{\zeta_2x_{n_1}}+x_{n_1+2}f_{\zeta_2y_{n_1}}+f_{x_{n_1}x_{n_1+1}}-
f_{y_{n_1}y_{n_1+1}}=0.\hspace{4cm}(4.6)\end{eqnarray*} According to
(3.12),
$$y_{n_1-1}=x_{n_1}^{-1}y_{n_1}x_{n_1-1}-x_{n_1}^{-1}\zeta_1,\qquad
y_{n_1+2}=x_{n_1+1}^{-1}\zeta_2+x_{n_1+1}^{-1}y_{n_1+1}x_{n_1+2}.\eqno(4.7)$$
Substituting (4.7) into (4.6), the coefficient of
$x_{n_1-1}x_{n_1+2}$ implies $f_{\zeta_1\zeta_2}=0$. Thus
$$f=g+h\qquad\mbox{with}\;\;g_{\zeta_2}=h_{\zeta_1}=0.\eqno(4.8)$$
Now (4.6) becomes
\begin{eqnarray*}& &
x_{n_1}^{-1}\zeta_1g_{\zeta_1x_{n_1+1}}-(x_{n_1}^{-1}y_{n_1}g_{\zeta_1x_{n_1+1}}+g_{\zeta_1y_{n_1+1}})x_{n_1-1}
+(x_{n_1+1}^{-1}y_{n_1+1}h_{\zeta_2x_{n_1}}+h_{\zeta_2y_{n_1}})x_{n_1+2}
\\ & &+x_{n_1+1}^{-1}\zeta_2h_{\zeta_2x_{n_1}}+g_{x_{n_1}x_{n_1+1}}-
g_{y_{n_1}y_{n_1+1}}+h_{x_{n_1}x_{n_1+1}}-
h_{y_{n_1}y_{n_1+1}}=0,\hspace{2.7cm}(4.9)\end{eqnarray*} which
implies
$$x_{n_1}^{-1}y_{n_1}g_{\zeta_1x_{n_1+1}}+g_{\zeta_1y_{n_1+1}}=0,\qquad
x_{n_1+1}^{-1}y_{n_1+1}h_{\zeta_2x_{n_1}}+h_{\zeta_2y_{n_1}}=0.\eqno(4.10)$$

For the representation purpose, we assume that $g$ is polynomial in
$\zeta_1$ with $g|_{\zeta_1=0}=0$ and $h$ is polynomial in
$\zeta_2$. Set
$$\zeta_3=x_{n_1}y_{n_1+1}-x_{n_1+1}y_{n_1}.\eqno(4.11)$$
By (4.10),
$$g\;\;\mbox{is a function
in}\;\;x_{n_1},y_{n_1},\zeta_1,\zeta_3.\eqno(4.12)$$ Moreover, (4.9)
says
$$-x_{n_1}^{-1}y_{n_1}\zeta_1g_{\zeta_1\zeta_3}-y_{n_1}g_{x_{n_1}\zeta_3}-x_{n_1}g_{y_{n_1}\zeta_3}=0.\eqno(4.13)$$
Again we can assume that $g=\hat g+\td g$ is polynomial in
$x_{n_1},y_{n_1},\zeta_3$ with $\hat g|_{\zeta_3=0}=0$ and $\td
g_{\zeta_3}=0$. Then (4.13) is equivalent to
$$y_{n_1}\zeta_1\hat g_{\zeta_1}+x_{n_1}y_{n_1}\hat g_{x_{n_1}}+x_{n_1}^2\hat g_{y_{n_1}}=0.\eqno(4.14)$$
This shows that $\hat g$ is a function in
$\zeta_1/x_{n_1},x_{n_1}^2-y_{n_1}^2,\zeta_3$. If $\hat g$ is a
polynomial, then $\hat g=0$. So the polynomial solution of $g$ must
be a polynomial in $x_{n_1},y_{n_2},\zeta_1$ with $g_{\zeta_1}\neq
0$. Similarly, if $h_{\zeta_2}\neq 0$ and $h|_{\zeta_2=0}=0$, the
polynomial solution of $h$ must be a polynomial in
$x_{n+1},y_{n+1},\zeta_2$. Assume $h_{\zeta_2}=0$. Then
$$h_{x_{n_1}x_{n_1+1}}-h_{y_{n_1}y_{n_1+1}}=0\eqno(4.15)$$
by (4.9).

 By Lemma 2.1, (3.78)-(3.80) and (4.2),
the polynomial solution of $h$ must be in
$$\mbox{Span}\{\eta^{m_3}(x_{n_1}^{m_1}y_{n_1+1}^{m_2})\mid m_1,m_2,m_3\in\mbb{N}\}.\eqno(4.16)$$
Therefore, a singular vector in ${\cal B}$ must be a nonzero weight
vector in
$$\mbox{Span}\{x_{n_1}^{m_1}y_{n_1}^{m_2}\zeta_1^{m_3+1},
x_{n_1+1}^{m_1}y_{n_1+1}^{m_2}\zeta_2^{m_3+1},
\eta^{m_3}(x_{n_1}^{m_1}y_{n_1+1}^{m_2})\mid
m_i\in\mbb{N}\}.\eqno(4.17)$$
 Note
$$x_{n_1}^{m_1}y_{n_1}^{m_2}\zeta_1^{m_3+1}\in{\cal
B}_{\la-m_1-m_3-1,m_2+m_3+1\ra},\eqno(4.18)$$
$$x_{n_1+1}^{m_1}y_{n_1+1}^{m_2}\zeta_2^{m_3+1}\in {\cal B}_{\la
m_1+m_3+1,-m_2-m_3-1\ra}.\eqno(4.19)$$ Moreover,
$${\cal D}(x_{n_1}^{m_1}y_{n_1}^{m_2}\zeta_1^{m_3+1})=-
m_2x_{n_1}^{m_1+1}y_{n_1}^{m_2-1}\zeta_1^{m_3+1}=0\dar
m_2=0\eqno(4.20)$$ and
$${\cal D}(x_{n_1+1}^{m_1}y_{n_1+1}^{m_2}\zeta_2^{m_3+1})
=-m_1x_{n_1+1}^{m_1-1}y_{n_1+1}^{m_2+1}\zeta_2^{m_3+1}=0\dar
m_1=0\eqno(4.21)$$ by (3.12) and (4.1). Furthermore,
$$x_{n_1}^{m_1}y_{n_1}^{m_2}\zeta_1^{m_3+1}=
\frac{\eta^{m_2}(x_{n_1}^{m_1+m_2}\zeta_1^{m_3+1})}{\prod_{r=1}^{m_2}(m_1+r)},\;\;
x_{n_1+}^{m_1}y_{n_1+1}^{m_2}\zeta_2^{m_3+1}=
\frac{\eta^{m_1}(y_{n_1+1}^{m_1+m_2}\zeta_2^{m_3+1})}{\prod_{r=1}^{m_1}(m_2+r)}
\eqno(4.21)$$ by (4.2). Indeed,
$$\eta^{m_1+1}(x_{n_1}^{m_1}\zeta_1^{m_2})=\eta^{m_1+1}(y_{n_1+1}^{m_1}\zeta_2^{m_2})=0
\qquad\for\;\;m_1,m_2\in\mbb{N}.\eqno(4.22)$$

Since $x_{n_1}^{m_1}y_{n_1+1}^{m_2}\in{\cal H}_{\la-m_1,-m_2\ra}$,
 (3.45) says that $\eta^m(x_{n_1}^{m_1}y_{n_1+1}^{m_2})$ with
$m>0$ is a singular vector only if $m=m_1+m_2+1$. But
$\eta^{m_1+m_2+1}(x_{n_1}^{m_1}y_{n_1+1}^{m_2})=0$ by (4.2). Thus
any singular vector in ${\cal H}$ (cf. (3.38)) is a nonzero weight
vector in
$$\mbox{Span}\{x_{n_1}^{m_1}\zeta_1^{m_2+1},y_{n_1+1}^{m_1}\zeta_2^{m_2+1},
x_{n_1}^{m_1}y_{n_1+1}^{m_2}\mid m_1,m_2\in\mbb{N}\}.\eqno(4.23)$$
Since ${\cal B}$ is nilpotent with respect to $sl(n,\mbb{F})_+$ (cf.
(2.30)), any nonzero submodule of ${\cal B}$ has a singular vector.
The above fact implies ${\cal H}_{\la\ell_1,\ell_2\ra}=\{0\}$ for
$\ell_1,\ell_2\in\mbb{Z}$ such that $\ell_1+\ell_2>0.$ Observe that
\begin{eqnarray*}& &(x_{n_1}^{m_1}\zeta_1^{m_2}|x_{n_1}^{m_1}\zeta_1^{m_2})\\&=&(\sum_{i=0}^{m_2}{m_2\choose
i}(-1)^ix_{n_1-1}^{m_2-i}x_{n_1}^{m_1+i}y_{n_1-1}^iy_{n_1}^{m_2-i}|\sum_{i=0}^{m_2}{m_2\choose
i}(-1)^ix_{n_1-1}^{m_2-i}x_{n_1}^{m_1+i}y_{n_1-1}^iy_{n_1}^{m_2-i})\\
&=&(-1)^{m_1+m_2}m_2!\sum_{i=0}^{m_2}{m_2\choose
i}(m_1+i)!(m_2-i))!\neq 0\hspace{5.2cm}(4.24)\end{eqnarray*} by
(3.55). Similarly,
$(y_{n_1+1}^{m_1}\zeta_2^{m_2}|y_{n_1}^{m_1}\zeta_2^{m_2})\neq 0$.

Next we assume $n_1=n_2=1$ and $n\geq 3$. By the arguments in the
above, a singular vector in ${\cal B}$ must be a nonzero weight
vector in
$$\mbox{Span}\{x_{n_1+1}^{m_1}y_{n_1+1}^{m_2}\zeta_2^{m_3+1},
\eta^{m_3}(x_{n_1}^{m_1}y_{n_1+1}^{m_2})\mid
m_i\in\mbb{N}\}.\eqno(4.25)$$ Thus any singular vector in ${\cal H}$
(cf. (3.38)) is a nonzero weight vector in
$$\mbox{Span}\{y_{n_1+1}^{m_1}\zeta_2^{m_2+1},
x_{n_1}^{m_1}y_{n_1+1}^{m_2}\mid m_1,m_2\in\mbb{N}\}.\eqno(4.26)$$
The above fact implies ${\cal H}_{\la\ell_1,\ell_2\ra}=\{0\} $ for
$\ell_1,\ell_2\in\mbb{Z}$ such that $\ell_1+\ell_2>0 $ or
$\ell_2>0.$

Consider the subcase $n_1=n_2=n-1$ and $n\geq 3$. A singular vector
in ${\cal B}$ must be a nonzero weight vector in
$$\mbox{Span}\{x_{n_1}^{m_1}y_{n_1}^{m_2}\zeta_1^{m_3+1},
\eta^{m_3}(x_{n_1}^{m_1}y_{n_1+1}^{m_2})\mid
m_i\in\mbb{N}\}.\eqno(4.27)$$ Thus any singular vector in ${\cal H}$
(cf. (3.38)) is a nonzero weight vector in
$$\mbox{Span}\{x_{n_1}^{m_1}\zeta_1^{m_2+1},
x_{n_1}^{m_1}y_{n_1+1}^{m_2}\mid m_1,m_2\in\mbb{N}\}.\eqno(4.28)$$
The above fact implies ${\cal H}_{\la\ell_1,\ell_2\ra}=\{0\} $ for
$\ell_1,\ell_2\in\mbb{Z}$ such that $\ell_1+\ell_2>0 $ or
$\ell_1>0.$

Suppose $n_1=n_2=1$ and $n=2$. A singular vector in ${\cal B}$ must
be a nonzero weight vector in
$$\mbox{Span}\{\eta^{m_3}(x_1^{m_1}y_2^{m_2})\mid
m_i\in\mbb{N}\}.\eqno(4.29)$$ Thus any singular vector in ${\cal H}$
(cf. (3.38)) is a nonzero weight vector in
$$\mbox{Span}\{x_1^{m_1}y_2^{m_2}\mid m_1,m_2\in\mbb{N}\}.\eqno(4.30)$$
The above fact implies ${\cal H}_{\la\ell_1,\ell_2\ra}=\{0\}$ for
$\ell_1,\ell_2\in\mbb{Z}$ such that $\ell_1>0 $ or $\ell_2>0.$

Finally, we assume $n_1=n_2=n$. A singular vector in ${\cal B}$ must
be a nonzero weight vector in
$$\mbox{Span}\{x_{n_1}^{m_1}y_{n_1}^{m_2}\zeta_1^{m_3}\mid m_i\in\mbb{N}\}.\eqno(4.31)$$
 Thus any singular vector in
${\cal H}$ (cf. (3.38)) is a nonzero weight vector in
$$\mbox{Span}\{x_{n_1}^{m_1}\zeta_1^{m_2}\mid m_1,m_2\in\mbb{N}\}.\eqno(4.32)$$
The above fact implies ${\cal H}_{\la\ell_1,\ell_2\ra}=\{0\}$ for
$\ell_1,\ell_2\in\mbb{Z}$ such that $\ell_1+\ell_2>0.$ Indeed, all
${\cal B}_{\la -m_1, m_2\ra}$ with $m_1,m_2\in\mbb{N}$ are
finite-dimensional and completely reducible by Weyl's Theorem of
complete reducibility. Moreover, their irreducible summands are
completely determined by (4.31).

By the arguments of (3.55)-(3.70), we obtain:\psp

{\bf Theorem 4.1}. {\it Suppose $n_1=n_2$. Let
 $\ell_1,\ell_2\in\mbb{Z}$ such that $\ell_2\geq 0$ when $n_1=n$.  Assume $\ell_1+\ell_2\leq 0$ and: (a)
 $\ell_2\leq 0 $  if $n_1=1$ and $n\geq 3$;
 (b) $\ell_1\leq 0$ if $n_1=n-1$ and $n\geq 3$; (c) $\ell_1,\ell_2\leq
 0$ when $n_1=1$ and $n=2$. Then
   ${\cal H}_{\la\ell_1,\ell_2\ra}$ is an
irreducible highest-weight $sl(n,\mbb{F})$-module and
$${\cal B}_{\la\ell_1,\ell_2\ra}=\bigoplus_{m=0}^\infty\eta^m({\cal
H}_{\la\ell_1-m,\ell_2-m\ra})\eqno(4.33)$$ is an orthogonal
decomposition of irreducible submodules. The symmetric bilinear form
restricted to $\eta^m({\cal H}_{\la\ell_1-m,\ell_2-m\ra})$.
 In particular, ${\cal
B}_{\la\ell_1,\ell_2\ra}={\cal H}_{\la\ell_1,\ell_2\ra}\oplus
\eta({\cal B}_{\la\ell_1-1,\ell_2-1\ra})$. If the conditions fails,
${\cal H}_{\la\ell_1,\ell_2\ra}=\{0\}$. When $n_1=n_2=n$, all the
above irreducible modules are of finite-dimensional}.\psp

Suppose $n_1<n-1$. Let $m_1,m_2\in\mbb{N}$. The subspace ${\cal
H}_{\la-m_1,-m_2\ra}$ has  a highest-weight vector
$x_{n_1}^{m_1}y_{n_1+1}^{m_2}$ of weight
$m_1(1-\dlt_{1,n_1})\lmd_{n_1-1}-(m_1+m_2+2)\lmd_{n_1}+m_2\lmd_{n_1+1}$.
If $n_1\geq 2$, the subspace ${\cal H}_{\la-m_1-m_2-1,m_2+1\ra}$ has
a highest-weight vector $x_{n_1}^{m_1}\zeta_1^{m_2+1}$ of weight
$(m_2+1)\lmd_{n_1-2}-m_1\lmd_{n_1-1}-(m_1+m_2+3)\lmd_{n_1}$. The
subspace ${\cal H}_{\la m_1+1,-m_2-m_1-1\ra}$ has  a highest-weight
vector $y_{n_1+1}^{m_2}\zeta_2^{m_1+1}$ of weight
$-(m_1+m_2+3)\lmd_{n_1}+m_2\lmd_{n_1+1}-(m_1+1)(1-\dlt_{n_1,n-2})\lmd_{n_1+2}$.

Consider $n_1=n-1$. The subspace ${\cal H}_{\la-m_1,-m_2\ra}$ has a
highest-weight vector $x_{n_1}^{m_1}y_{n_1+1}^{m_2}$ of weight
$m_1(1-\dlt_{n,2})\lmd_{n-2}-(m_1+m_2+2)\lmd_{n-1}$. If $n\geq 3$,
the subspace ${\cal H}_{\la-m_1-m_2-1,m_2+1\ra}$ has  a
highest-weight vector $x_{n_1}^{m_1}\zeta_1^{m_2+1}$ of weight
$(m_2+1)(1-\dlt_{n,3})\lmd_{n-3}-m_1\lmd_{n-2}-(m_1+m_2+3)\lmd_{n-1}$.

Assume $n_1=n$. The subspace ${\cal H}_{\la-m_1-m_2,m_2\ra}$ has a
highest-weight vector $x_n^{m_1}\zeta_1^{m_2}$ of weight
$m_2(1-\dlt_{n,2})\lmd_{n-2}+m_1\lmd_{n-1}$.\psp

Now we want to find an explicit expression for ${\cal
H}_{\la\ell_1,\ell_2\ra}$ when it is irreducible. Set
$${\cal
G}'=\sum_{i=1}^{n_1}\sum_{j=n_1+1}^n\mbb{F}E_{j,i},\eqno(4.34)$$
$$\hat{\cal
G}=H+\sum_{r,s\in\ol{1,n_1}\;\mbox{or}\;r,s\in\ol{n_1+1,n};r\neq
s}\mbb{F}E_{r,s}+\sum_{i=1}^{n_1}\sum_{j=n_1+1}^n\mbb{F}E_{i,j}.\eqno(4.35)$$
Then ${\cal G}'$ and $\hat{\cal G}$ are Lie subalgebras of
$sl(n,\mbb{F})$ and $sl(n,\mbb{F})={\cal G}'\oplus \hat{\cal G}$. By
PBW Theorem, $U(sl(n,\mbb{F}))=U({\cal G}')U(\hat{\cal G})$.
According to (1.6)-(1.8),
$$E_{r,s}|_{\cal B}=-x_s\ptl_{x_r}-y_s\ptl_{y_r},\qquad
E_{p,q}|_{\cal B}=x_p\ptl_{x_q}+y_p\ptl_{y_q},\eqno(4.36)$$
$$E_{r,p}|_{\cal B}=\ptl_{x_r}\ptl_{x_p}-\ptl_{y_r}\ptl_{y_p},\qquad
E_{p,r}|_{\cal B}=-x_rx_p+y_ry_p\eqno(4.37)$$ for $r,s\in\ol{1,n_1}$
and $p,q\in\ol{n_1+1,n}$.

First we assume $n_1<n$. For $m_1,m_2\in\mbb{N}$, we have
\begin{eqnarray*}\qquad{\cal H}_{\la
-m_1,-m_2\ra}&=&U(sl(n,\mbb{F}))(x_{n_1}^{m_1}y_{n_1+1}^{m_2})=U({\cal
G}')U(\hat{\cal G})(x_{n_1}^{m_1}y_{n_1+1}^{m_2})
\\ &=&\mbox{Span}\{[\prod_{r=1}^{n_1}x_r^{l_r}]
[\prod_{s=1}^{n-n_1}y_{n_1+s}^{k_s}][\prod_{r=1}^{n_1}\prod_{s=1}^{n-n_1}(x_rx_{n_1+s}-
y_ry_{n_1+s})^{l_{r,s}}]\\& &\qquad\;\;\mid
l_r,k_s,l_{r,s}\in\mbb{N};\sum_{r=1}^{n_1}l_r=m_1;\sum_{s=1}^{n-n_1}k_s=m_2\}\hspace{2.4cm}(4.38)\end{eqnarray*}
by (4.36) and (4.37).
 Furthermore, we assume $n_1>1$. We let
\begin{eqnarray*}& &{\cal H}_{\la
-m_1-m_2,m_2\ra}'\\&=&\mbox{Span}\{[\prod_{r=1}^{n_1}x_r^{l_r}]
[\prod_{1\leq p<q\leq
n_1}(x_py_q-x_qy_p)^{k_{p,q}}][\prod_{r=1}^{n_1}\prod_{s=1}^{n-n_1}(x_rx_{n_1+s}-
y_ry_{n_1+s})^{l_{r,s}}]\\& &\qquad\;\;\mid
l_r,k_{p,q},l_{r,s}\in\mbb{N};\sum_{r=1}^{n_1}l_r=m_1;\sum_{1\leq
p<q\leq n_1}k_{p,q}=m_2\}.\hspace{4cm}(4.39)\end{eqnarray*} By
(3.38), (3.40) and (4.1), we have ${\cal H}_{\la
-m_1-m_2,m_2\ra}'\subset {\cal H}_{\la -m_1-m_2,m_2\ra}$. Moreover,
(4.37) and (4.38) yield
$${\cal H}_{\la -m_1-m_2,m_2\ra}=U(sl(n,\mbb{F}))(x_{n_1}^{m_1}\zeta_1^{m_2})=U({\cal
G}')U(\hat{\cal G})(x_{n_1}^{m_1}\zeta_1^{m_2})\subset {\cal H}_{\la
-m_1-m_2,m_2\ra}'.\eqno(4.40)$$ Thus ${\cal H}_{\la
-m_1-m_2,m_2\ra}'= {\cal H}_{\la -m_1-m_2,m_2\ra}$. Symmetrically,
if $n_1=n_2<n-1$,
\begin{eqnarray*}& &{\cal H}_{\la
m_2,-m_1-m_2\ra}=\mbox{Span}\{ [\prod_{n_1+1\leq p<q\leq
n}(x_py_q-x_qy_p)^{k_{p,q}}][\prod_{r=1}^{n_1}\prod_{s=n_1+1}^n(x_rx_s-
y_ry_s)^{l_{r,s}}]\\& &\times
[\prod_{r=1}^{n-n_1}y_{n_1+r}^{l_r}]\mid
l_r,k_{p,q},l_{r,s}\in\mbb{N};\sum_{r=1}^{n_1}l_r=m_1;\sum_{n_1+1\leq
p<q\leq n}k_{p,q}=m_2\}.\hspace{2.7cm}(4.41)\end{eqnarray*}

When $n_1=n_2=n$, by the arguments between (4.39) and(4.40),
\begin{eqnarray*}{\cal H}_{\la
-m_1-m_2,m_2\ra}&=&\mbox{Span}\{[\prod_{r=1}^nx_r^{l_r}]
[\prod_{1\leq p<q\leq n}(x_py_q-x_qy_p)^{k_{p,q}}]\\&
&\qquad\;\;\mid
l_r,k_{p,q}\in\mbb{N};\sum_{r=1}^nl_r=m_1;\sum_{1\leq p<q\leq
n}k_{p,q}=m_2\},\hspace{2.5cm}(4.42)\end{eqnarray*} which is of
finite-dimensional.

\section{The $o(2n,\mbb{F})$-Variation}

Recall that ${\cal B}=\mbb{F}[x_1,...,x_n,y_1,...,y_n]$ and the
representation of $o(2n,\mbb{F})$ on ${\cal B}$  defined by
(1.14)-(1.16). It is easy to verify
$$T\xi=\xi T\;\;\mbox{on}\;\;{\cal B}\qquad\for\;\;\xi\in o(2n,\mbb{F});T=\flat,\flat',{\cal
D},\eta\eqno(5.1)$$ by (1.9), (1.10), (2.13) and (3.4). Recall the
notions ${\cal B}_{\la k\ra}$ and ${\cal H}_{\la k\ra}$ defined in
(1.17). The ${\cal B}=\bigoplus_{k\in\mbb{Z}}{\cal B}_{\la k\ra}$
forms a $\mbb{Z}$-graded algebra and $${\cal H}_{\la
k\ra}=\bigoplus_{\ell_1,\ell_2\in\mbb{Z};\ell_1+\ell_2=k}{\cal
H}_{\la\ell_1,\ell_2\ra}.\eqno(5.2)$$    Moreover, ${\cal B}_{\la
k\ra}$ and ${\cal H}_{\la k\ra}$ are $o(2n,\mbb{F})$-submodules.
Recall ${\cal K}=\sum_{i,j=1}^n\mbb{F}(E_{i,j}-E_{n+j,n+i})$.\psp

{\bf Theorem 5.1}. {\it For any  $n_1-n_2+1-\dlt_{n_1,n_2}\geq
k\in\mbb{Z}$, ${\cal H}_{\la k\ra}$ is an irreducible
$o(2n,\mbb{F})$-submodule and
$${\cal B}_{\la k\ra}=\bigoplus_{i=0}^\infty\eta^i({\cal H}_{\la k-2i\ra}) \eqno(5.3)$$
is an orthogonal decomposition of irreducible submodules. In
particular, ${\cal B}_{\la k\ra}={\cal H}_{\la k\ra}\oplus
\eta({\cal B}_{\la k-2\ra})$. Moreover, the bilinear form
$(\cdot|\cdot)$ restricted to $\eta^i({\cal H}_{\la k-2i\ra})$ is
nondegenerate. Furthermore, ${\cal H}_{\la k\ra}$ has a basis
\begin{eqnarray*}\qquad& &\big\{\sum_{i=0}^\infty\frac{(x_{n_1+1}y_{n_1+1})^i({\cal
D}-\ptl_{x_{n_1}+1}\ptl_{y_{n_1}+1})^i(x^\al
y^\be)}{\prod_{r=1}^i(\al_{n_1+1}+r)(\be_{n_1+1}+r)}
\mid\al,\be\in\mbb{N}\:^n;\\ &
&\al_{n_1+1}\be_{n_1+1}=0;-\sum_{i=1}^{n_1}\al_i+\sum_{r=n_1+1}^n\al_r+
\sum_{i=1}^{n_2}\be_i-\sum_{r=n_2+1}^n\be_r=k
\big\}\hspace{2.5cm}(5.4)\end{eqnarray*} when $n_1<n_2$. The module
${\cal H}_{\la k\ra}$ under the assumption is of highest-weight type
only if $n_2=n$, in which case $x_{n_1}^{-k}$ is a highest-weight
vector with weight
$-k\lmd_{n_1-1}+(k-1)\lmd_{n_1}+[(k-1)\dlt_{n_1,n-1}-2k\dlt_{n_1,n}]\lmd_n.$
 When $n_1=n_2=n$, all the
irreducible modules ${\cal H}_{\la k\ra}$ with $0\geq k\in\mbb{Z}$
are of $({\cal G},{\cal K})$-type.}

{\it Proof}. Let $n_1-n_2+1\geq k\in\mbb{Z}.$ Note
$sl(n,\mbb{F})|_{\cal B}$ is a subalgebra of $o(2n,\mbb{F})|_{\cal
B}$. Suppose $n_1+1<n_2<n$. By (5.2), Theorem 3.1 and the paragraph
below, the $sl(n,\mbb{F})$-singular vectors in ${\cal H}_{\la k\ra}$
are: for $m_1,m_2\in\mbb{N}$,
$$x_{n_1}^{m_1}y_{n_2+1}^{m_2}\qquad\mbox{with}\;-(m_1+m_2)=k,\eqno(5.5)$$
$$x_{n_1+1}^{m_1}y_{n_2+1}^{m_2}\qquad\mbox{with}\;m_1-m_2=k,\eqno(5.6)$$
$$x_{n_1}^{m_1}y_{n_2}^{m_2}\qquad\mbox{with}\;-m_1+m_2=k.\eqno(5.7)$$

Note
$$(E_{n+n_2+1,n_1}-E_{n+n_1,n_2+1})|_{\cal
B}=-x_{n_1}\ptl_{y_{n_2+1}}-y_{n_1}\ptl_{x_{n_2+1}}\eqno(5.8)$$ by
(1.16). So
$$(E_{n+n_2+1,n_1}-E_{n+n_1,n_2+1})^{m_2}(x_{n_1}^{m_1}y_{n_2+1}^{m_2})=(-1)^{m_2}m_2!x_{n_1}^{-k}\eqno(5.9)$$
for the vectors in (5.5). Moreover,
$$(E_{n+n_2+1,n_1+1}-E_{n+n_1+1,n_2+1})|_{\cal
B}=\ptl_{x_{n_1+1}}\ptl_{y_{n_2+1}}-y_{n_1+1}\ptl_{x_{n_2+1}}\eqno(5.10)$$
again by (1.16), which implies
$$(E_{n+n_2+1,n_1+1}-E_{n+n_1+1,n_2+1})^{m_2}(x_{n_1+1}^{m_1}y_{n_2+1}^{m_2})=m_1![\prod_{r=0}^{m_1-1}(m_2-r)]
y_{n_2+1}^{-k}\eqno(5.11)$$ for the vectors in (5.6). Furthermore,
$$(E_{n_1,n+n_2}-E_{n_2,n+n_1})|_{\cal
B}=\ptl_{x_{n_1}}\ptl_{y_{n_2}}-x_{n_2}\ptl_{y_{n_1}}\eqno(5.12)$$
by (1.15), which implies
$$(E_{n_1,n+n_2}-E_{n_2,n+n_1})^{m_2}(x_{n_1}^{m_1}y_{n_2}^{m_2})=
m_2![\prod_{r=0}^{m_2-1}(m_1-r)] x_{n_1}^{-k}\eqno(5.13)$$ for the
vectors in (5.7).

On the other hand,
$$(E_{n_1,n+n_2+1}-E_{n_2+1,n+n_1})|_{\cal
B}=-y_{n_2+1}\ptl_{x_{n_1}}-x_{n_2+1}\ptl_{y_{n_1}}\eqno(5.14)$$ by
(1.15), which implies
$$(E_{n_1,n+n_2+1}-E_{n_2+1,n+n_1})^{m_2}(x_{n_1}^{-k})=(-1)^{m_2}[\prod_{r=0}^{m_2-1}(-k-r)]
x_{n_1}^{m_1}y_{n_2+1}^{m_2}\eqno(5.15)$$ for the vectors in (5.5).
Moreover,
$$(E_{n_1+1,n+n_2+1}-E_{n_2+1,n+n_1+1})|_{\cal
B}=-x_{n_1+1}y_{n_2+1}-x_{n_2+1}\ptl_{y_{n_1+1}}\eqno(5.16)$$ by
(1.15), which implies
$$(E_{n_1+1,n+n_2+1}-E_{n_2+1,n+n_1+1})^{m_2}(y_{n_2+1}^{-k})=(-1)^{m_2}
x_{n_1}^{m_1}y_{n_2+1}^{m_2} \eqno(5.17)$$ for the vectors in (5.6).
Furthermore,
$$(E_{n+n_2,n_1}-E_{n+n_1,n_2})|_{\cal
B}=-x_{n_1}y_{n_2}-y_{n_1}\ptl_{x_{n_2}}\eqno(5.18)$$ by (1.16),
which implies
$$(E_{n+n_2,n_1}-E_{n+n_1,n_2})^{m_2}(x_{n_1}^{-k})=
(-1)^{m_2}x_{n_1}^{m_1}y_{n_2}^{m_2}\eqno(5.19)$$ for the vectors in
(5.7). Thus for any two vectors in (5.5)-(5.7), there exists an
element in the universal enveloping algebra $U(o(2n,\mbb{F}))$ which
carries one to another. On the other hand, the vectors in
(5.5)-(5.7) have distinct weights (see the paragraph below Theorem
3.1). Thus any nonzero submodule of ${\cal H}_{\la k\ra}$ must
contain one of the vectors in (5.5)-(5.7). Hence all the vectors in
(5.5)-(5.7) are in the submodule by (5.8)-(5.19). Therefore, the
submodule must be equal to ${\cal H}_{\la k\ra}$, that is, ${\cal
H}_{\la k\ra}$ is irreducible. By (5.16) and (5.18), ${\cal H}_{\la
k\ra}$ is not of highest-weight type. The equation (5.3) follows
from Theorem 3.1 and (5.2).

 Assume $n_1+1=n_2<n$. By
Theorem 3.2 and the paragraph below, the $sl(n,\mbb{F})$-singular
vectors in ${\cal H}_{\la k\ra}$ are those in (5.5)-(5.7). So the
theorem holds. Suppose $n_1<n_2=n$. By Theorems 3.1, 3.2 and the
paragraph below them, the $sl(n,\mbb{F})$-singular vectors in ${\cal
H}_{\la k\ra}$ are those in (5.7). Expressions (5.13) and (5.19)
imply the conclusions in the theorem.

 Recall
$$\zeta_1=x_{n_1-1}y_{n_1}-x_{n_1}y_{n_1-1},\;\;
\;\zeta_2=x_{n_2+1}y_{n_2+2}-x_{n_2+2}y_{n_2}. \eqno(5.20)$$ In the
case $n_1=n_2<n-1$, Theorem 4.1 tell us that the
$sl(n,\mbb{F})$-singular vectors in ${\cal H}_{\la k\ra}$ are those
in (5.5) and
$$x_{n_1}^{-k}\zeta_1^{m+1}\qquad\for\;\;m\in\mbb{N},\eqno(5.21)$$
$$y_{n_1+1}^{-k}\zeta_2^{m+1}\qquad\for\;\;m\in\mbb{N}.\eqno(5.22)$$
Again all the singular vectors have distinct weights. If $N$ is a
nonzero submodule of ${\cal H}_{\la k\ra}$, then $N$ must contain
one of the above $sl(n,\mbb{F})$-singular vectors. If $N$ contains a
singular vector in (5.5), then $x_{n_1}^{-k}\in N$ by (5.9). Suppose
$x_{n_1}^{-k}\zeta_1^{m+1}\in N$ for some $m\in\mbb{N}$.
 Note
$$(E_{n_1-1,n+n_1}-E_{n_1,n+n_1-1})|_{\cal
B}=\ptl_{x_{n_1-1}}\ptl_{y_{n_1}} -\ptl_{x_{n_1}}\ptl_{y_{n_1-1}}
\eqno(5.22)$$ by (1.15). Thus
\begin{eqnarray*}\qquad&
&(E_{n_1-1,n+n_1}-E_{n_1,n+n_1-1})^{m+1}(x_{n_1}^{-k}\zeta_1^{m+1})
\\ &=&\left[\sum_{r=0}^{m+1}(-1)^r{m+1\choose
r}(\ptl_{x_{n_1-1}}\ptl_{y_{n_1}})^{m+1-r}(\ptl_{x_{n_1}}\ptl_{y_{n_1-1}})^r\right]
\\ &&\left[\sum_{s=0}^{m+1}(-1)^s{m+1\choose
s}(x_{n_1-1}y_{n_1})^{m+1-s}x_{n_1}^{-k+s}y_{n_1-1}^s\right]
\\&=&\left(\sum_{r=0}^{m+1}{m+1\choose
r}^2[(m+1-r!)]^2r![\prod_{i=1}^r(-k+i)]\right)x_{n_1}^{-k}
\\ &=&[(m+1)!]^2\left(\sum_{r=0}^{m+1}{-k+r\choose
r}\right)x_{n_1}^{-k}\in N.\hspace{6.1cm}(5.23)\end{eqnarray*} So we
have $x_{n_1}^{-k}\in N$ again. Symmetrically, it holds if
 $y_{n_1+1}^{-k}\zeta_2^{m+1}\in N$ for some $m\in\mbb{N}$.
 Therefore, we always have $x_{n_1}^{-k}\in N$.

According to (5.15), $N$ contains all the singular vectors in (5.5).
Observe
$$(E_{n+n_1-1,n_1}-E_{n+n_1,n_1-1})|_{\cal B}=\zeta_1,\;\;
(E_{n_1+2,n+n_1+1}-E_{n_1+1,n+n_1+2})|_{\cal B}=\zeta_2\eqno(5.24)$$
as multiplication operators on ${\cal B}$ by (1.15) and (1.16). Thus
$$(E_{n+n_1-1,n_1}-E_{n+n_1,n_1-1})^{m+1}(x_{n_1}^{-k})=x_{n_1}^{-k}\zeta_1^{m+1},\eqno(5.25)$$
$$(E_{n_1+2,n+n_1+1}-E_{n_1+1,n+n_1+2})^{m+1}(x_{n_1}^{-k})=x_{n_1}^{-k}\zeta_2^{m+1}\in
N.\eqno(5.26)$$ Thus $N$ contains all the $sl(n,\mbb{F})$-singular
vectors in ${\cal H}_{\la k\ra}$, which implies that it contains all
${\cal H}_{\la\ell_1,\ell_2\ra}\subset {\cal H}_{\la k\ra}$. So
$N={\cal H}_{\la k\ra}$, that is, ${\cal H}_{\la k\ra}$ is an
irreducible $o(2n,\mbb{F})$-module, which is of $({\cal G},{\cal
K})$-type if $n_1=n_2=n$ by (5.2). The basis (5.4) is obtained by
(3.72) and (5.2). $\qquad\Box$\psp

Finally, we want to find an expression for ${\cal H}_{\la k\ra}$ for
$0\geq k\in\mbb{Z}$ when $n_1=n_2$.

First we assume $n_1=n_2=1$ and $n\geq 3$. According to (4.26),
(4.38) and (4.41)
\begin{eqnarray*}& &{\cal H}_{\la
-k\ra}\\ &=&\mbox{Span}\{[\prod_{r=2}^ny_r^{\hat l_r}][\prod_{2\leq
p<q\leq n}(x_py_q-x_qy_p)^{\hat k_{p,q}}][\prod_{s=2}^n(x_1x_s-
y_1y_s)^{\hat l_s}],x_1^l [\prod_{s=2}^ny_s^{k_s}]\\
&&\times[\prod_{s=2}^n(x_1x_s- y_1y_s)^{l_s}]\mid l,k_s,l_s,\hat
l,\hat k_{p,q},\hat l_s\in\mbb{N};l+\sum_{s=2}^nk_s=\sum_{r=2}^n\hat
l_r=k\}.\hspace{1.8cm}(5.27)\end{eqnarray*} Next we consider the
subcase $1<n_1=n_2<n-1$.  By (4.23), (4.38), (4.39) (note ${\cal
H}_{\la -m_1-m_2,m_2\ra}'={\cal H}_{\la -m_1-m_2,m_2\ra}$) and
(4.41), we have
\begin{eqnarray*}& &{\cal H}_{\la
-k\ra}\\ &=&\mbox{Span}\{[\prod_{r=1}^{n_1}x_r^{l_r'}] [\prod_{1\leq
p<q\leq
n_1}(x_py_q-x_qy_p)^{k_{p,q}'}][\prod_{r=1}^{n_1}\prod_{s=n_1+1}^n(x_rx_s-
y_ry_s)^{l_{r,s}'}],\\
& &[\prod_{r=1}^{n-n_1}y_{n_1+r}^{\hat l_r}][\prod_{n_1+1\leq
p<q\leq n}(x_py_q-x_qy_p)^{\hat
k_{p,q}}][\prod_{r=1}^{n_1}\prod_{s=n_1+1}^n(x_rx_s- y_ry_s)^{\hat
l_{r,s}}],\\ & &[\prod_{r=1}^{n_1}x_r^{l_r}]
[\prod_{s=1}^{n-n_1}y_{n_1+s}^{k_s}][\prod_{r=1}^{n_1}\prod_{s=1}^{n-n_1}(x_rx_{n_1+s}-
y_ry_{n_1+s})^{l_{r,s}}]\mid l_r,k_s,l_{r,s},l_r',k_{p,q}',\\ &
&l_{r,s}',\hat l_r,\hat k_{p,q},\hat
l_{r,s}\in\mbb{N};\sum_{r=1}^{n_1}l_r+\sum_{s=1}^{n-n_1}k_s=\sum_{r=1}^{n_1}l_r'=\sum_{r=1}^{n-n_1}\hat
l_r=k\}.\hspace{3.6cm}(5.28)\end{eqnarray*}  Consider the subcase
$n_1=n_2=n-1$ and $n\geq 3$. By (4.28), (4.38) and (4.39) (note
${\cal H}_{\la -m_1-m_2,m_2\ra}'={\cal H}_{\la -m_1-m_2,m_2\ra}$),
we obtain
\begin{eqnarray*}& &{\cal H}_{\la
-k\ra}\\ &=&\mbox{Span}\{[\prod_{r=1}^{n-1}x_r^{l_r'}] [\prod_{1\leq
p<q\leq n-1}(x_py_q-x_qy_p)^{k_{p,q}'}][\prod_{r=1}^{n-1}(x_rx_n-
y_ry_n)^{\bar l_r'}],[\prod_{r=1}^{n-1}x_r^{l_r}] y_n^{\hat k}\\
&&\times[\prod_{r=1}^{n-1}(x_rx_n- y_ry_n)^{\bar l_r}]\mid l_r,\hat
k,\bar l_r,l_r',k_{p,q}',\bar
l_r'\in\mbb{N};\sum_{r=1}^{n-1}l_r+\hat
k=\sum_{r=1}^{n-1}l_r'=k\}.\hspace{1.5cm}(5.29)\end{eqnarray*}

Suppose $n_1=n_2=1$ and $n=2$. According to (4.30) and (4.38),
$${\cal H}_{\la
-k\ra}=\mbox{Span}\{[x_1^r y_2^s(x_1x_2- y_1y_2)^l\mid
r,s,l\in\mbb{N};r+s=k\}.\eqno(5.30)$$ Finally we assume $n_1=n_2=n$.
By (4.32) and (4.39) (note ${\cal H}_{\la -m_1-m_2,m_2\ra}'={\cal
H}_{\la -m_1-m_2,m_2\ra}$),
$${\cal H}_{\la -k\ra}=\mbox{Span}\{\prod_{r=1}^nx_r^{l_r}] [\prod_{1\leq
p<q\leq n}(x_py_q-x_qy_p)^{k_{p,q}}]\mid
l_r,k_{p,q}\in\mbb{N};\sum_{r=1}^nl_r=k\},\eqno(5.31)$$ whose
$({\cal G},{\cal K})$-module structure is given by ${\cal H}_{\la
-k\ra}=\bigoplus_{m=0}^\infty {\cal H}_{\la -k-m,m\ra}$ with ${\cal
H}_{\la -k-m,m\ra}$ given in (4.42).

\section{ The $o(2n+1,\mbb{F})$-Variation}

Recall
$$o(2n+1,\mbb{F})=o(2n,\mbb{F})\oplus\bigoplus_{i=1}^n
[\mbb{F}(E_{0,i}-E_{n+i,0})+\mbb{F}(E_{0,n+i}-E_{i,0})]\eqno(6.1)$$
 and ${\cal B}'=\mbb{F}[x_0,x_1,...,x_n,y_1,...,y_n]$.

Fix $n_1,n_2\in\ol{1,n}$ such that $n_1\leq n_2$. The representation
of $o(2n+1,\mbb{F})$ on ${\cal B}'$ by the differential operators in
(1.14)-(1.16), (1.19) and (1.20). Recall ${\cal B}'_{\la
k\ra}=\sum_{i=0}^\infty {\cal B}_{\la k-i\ra}x_0^i.$ Then all ${\cal
B}'_{\la k\ra}$ with $k\in\mbb{Z}$ are $o(2n+1,\mbb{F})$-submodules
and ${\cal B}'=\bigoplus_{k\in\mbb{Z}}{\cal B}'_{\la k\ra}$ forms a
$\mbb{Z}$-graded algebra. Moreover, the variated Laplace operator
${\cal D}'=\ptl_{x_0}^2+2{\cal D}$ by (1.21) and its dual
$\eta'=x_0^2+2\eta$ by (1.22).

A straightforward verification shows
$${\cal D}'\xi=\xi {\cal D}',\;\xi\eta'=\eta'\xi\;\mbox{on}\;{\cal
B}'\qquad\for\;\;\xi\in o(2n+1,\mbb{F}).\eqno(6.2)$$ As in the
introduction, ${\cal H}'_{\la k\ra}=\{f\in{\cal B}'{\la
k\ra}\mid{\cal D}'(f)=0\}.$  According to (6.2), ${\cal H}'_{\la
k\ra}$ is an $o(2n+1,\mbb{F})$-submodule.
 By Lemma 2.1 with $T_1=\ptl_{x_0}^2,\;T_1^-=\int_{(x_0)}^{(2)}$
 (cf. (2.6) and (2.7)) and $T_2=2{\cal D}$, we obtain
$${\cal H}'_{\la k\ra}=\left(\sum_{i=0}^\infty\frac{(-2)^ix_0^{2i}{\cal
D}^i}{(2i)!}\right)({\cal B}_{\la k\ra})\oplus
\left(\sum_{i=0}^\infty\frac{(-2)^ix_0^{2i+1}{\cal
D}^i}{(2i+1)!}\right)({\cal B}_{\la k-1\ra}).\eqno(6.3)$$ Recall
${\cal K}=\sum_{i,j=1}^n\mbb{F}(E_{i,j}-E_{n+j,n+i})$. \psp

{\bf Theorem 6.1}. {\it  For any  $n_1-n_2+1-\dlt_{n_1,n_2}\geq
k\in\mbb{Z}$, ${\cal H}'_{\la k\ra}$ is an irreducible
$o(2n+1,\mbb{F})$-submodule and
$${\cal B}'_{\la k\ra}=\bigoplus_{i=0}^\infty(\eta')^i({\cal H}'_{\la k-2i\ra})
 \eqno(6.4)$$ is an orthogonal decomposition of
irreducible submodules. In particular, ${\cal B}'_{\la k\ra}={\cal
H}'_{\la k\ra}\oplus \eta'({\cal B}'_{\la k-2\ra})$. Moreover, the
bilinear form $(\cdot|\cdot)$ restricted to $(\eta')^i({\cal H}_{\la
k-2i\ra}')$ is nondegenerate. Furthermore,  ${\cal H}_{\la k\ra}$
has a basis
\begin{eqnarray*}\qquad\;\;& &\big\{
\sum_{i=0}^\infty\frac{(-2)^ix_0^{2i+\iota}{\cal D}^i(x^\al
y^\be)}{(2i+\iota)!} \mid\al,\be\in\mbb{N}\:^n;\iota=0,1;\\
& &-\sum_{i=1}^{n_1}\al_i+\sum_{r=n_1+1}^n\al_r+
\sum_{i=1}^{n_2}\be_i-\sum_{r=n_2+1}^n\be_r=k-\iota
\big\}.\hspace{4.4cm}(6.5)\end{eqnarray*}The module ${\cal H}_{\la
k\ra}'$ under the assumption is of highest-weight type only if
$n_2=n$, in which case $x_{n_1}^{-k}$ is a highest-weight vector
with weight
$-k\lmd_{n_1-1}+(k-1)\lmd_{n_1}+[(k-1)\dlt_{n_1,n-1}-2k\dlt_{n_1,n}]\lmd_n.$
 When $n_1=n_2=n$, all the
irreducible modules ${\cal H}_{\la k\ra}$ with $0\geq k\in\mbb{Z}$
are of $({\cal G},{\cal K})$-type. }

{\it Proof}. Observe that
$$(x_0^rx^{\al}y^\be|x_0^sx^{\al_1}y^{\be_1})=\dlt_{r,s}\dlt_{\al,\al_1}
\dlt_{\be,\be_1}(-1)^{\sum_{i=1}^{n_1}\al_i+\sum_{r=n_2+1}^n\be_r}r!\al!\be!
\eqno(6.6)$$ for $r,s\in\mbb{N}$ and
$\al,\be,\al_1,\be_1\in\mbb{N}\:^n.$ By  (1.21) and (1.22),
$$({\cal D}'(f)|g)=(f|\eta'(g))\qquad\for\;\;f,g\in{\cal
B}'.\eqno(6.7)$$

 Let $n_1-n_2+1\geq k\in\mbb{Z}$. First by (5.3) and
(6.3),
$${\cal H}'_{\la
k\ra}=\bigoplus_{r=0}^\infty\left(\sum_{i=0}^\infty\frac{(-2)^ix_0^{2i}{\cal
D}^i}{(2i)!}\right)(\eta^r({\cal H}_{\la k-2r\ra}))\oplus
\bigoplus_{s=0}^\infty\left(\sum_{i=0}^\infty\frac{(-2)^ix_0^{2i+1}{\cal
D}^i}{(2i+1)!}\right)(\eta^s({\cal H}_{\la k-2s-1\ra})).\eqno(6.8)$$
Let $N$ be a nonzero submodule of ${\cal H}'_{\la k\ra}$. By
comparing weights and the arguments in (5.5)-(5.13) and
(5.21)-(5.23),  we have
$$\left(\sum_{i=0}^\infty\frac{(-2)^ix_0^{2i}{\cal
D}^i}{(2i)!}\right)(\eta^{m_1}(x_{n_1}^{-k+2m_1}))\in N\eqno(6.9)$$
for some $m_1\in\mbb{N}$ or
$$\left(\sum_{i=0}^\infty\frac{(-2)^ix_0^{2i+1}{\cal
D}^i}{(2i+1)!}\right)(\eta^{m_2}(x_{n_1}^{-k+2m_2+1}))\in
N\eqno(6.10)$$ for some $m_2\in\mbb{N}$.

Note
$$(E_{n_1,0}-E_{0,n+n_1})=\ptl_{x_0}\ptl_{x_{n_1}}-x_0\ptl_{y_{n_1}}\eqno(6.11)$$
by (1.19) and (1.20). Recall
$${\cal
D}=-\sum_{i=1}^{n_1}x_i\ptl_{y_i}+\sum_{r=n_1+1}^{n_2}\ptl_{x_r}\ptl_{y_r}-\sum_{s=n_2+1}^n
y_s\ptl_{x_s}\eqno(6.12)$$ and
$$\eta=\sum_{i=1}^{n_1}y_i\ptl_{x_i}+\sum_{r=n_1+1}^{n_2}x_ry_r+\sum_{s=n_2+1}^n
x_s\ptl_{y_s}.\eqno(6.13)$$ Then (3.44) gives
\begin{eqnarray*}&&(E_{n_1,0}-E_{0,n+n_1})\left[
\left(\sum_{i=0}^\infty\frac{(-2)^ix_0^{2i}{\cal
D}^i}{(2i)!}\right)(\eta^{m_1}(x_{n_1}^{-k+2m_1}))\right]
\\ &=&\left(\sum_{i=0}^\infty\frac{(i+1)m_1(-k+2m_1)
(-2)^{i+1}x_0^{2i+1}{\cal
D}^i}{(2i+1)!}\right)(\eta^{m_1-1}(x_{n_1}^{-k+2m_1-1}))
\\ & &+\left(\sum_{i=0}^\infty\frac{(-k+2m_1)
(-2)^{i+1}x_0^{2i+1}{\cal
D}^{i+1}}{(2i+1)!}\right)(\eta^{m_1}(x_{n_1}^{-k+2m_1-1}))
\\ & &-m_1(-k+2m_1)\left(\sum_{i=0}^\infty\frac{(-2)^ix_0^{2i+1}{\cal
D}^i}{(2i)!}\right)(\eta^{m_1-1}(x_{n_1}^{-k+2m_1-1}))\\ &=&
m_1(-k+2m_1)\left(\sum_{i=0}^\infty\frac{(-2)^ix_0^{2i+1}{\cal
D}^i}{(2i+1)!}\right)(\eta^{m_1-1}(x_{n_1}^{-k+2m_1-1}))
\\ & &+m_1(-k+2m_1)(m_1-k+n_1-n_2) \left(\sum_{i=0}^\infty\frac{(-2)^{i+1}x_0^{2i+1}{\cal
D}^i}{(2i+1)!}\right)(\eta^{m_1-1}(x_{n_1}^{-k+2m_1-1}))
\\&=&m_1(-k+2m_1)(2m_1-2k+2n_1-2n_2+1)\\ & &\times
\left(\sum_{i=0}^\infty\frac{(-2)^ix_0^{2i+1}{\cal
D}^i}{(2i+1)!}\right)(\eta^{m_1-1}(x_{n_1}^{-k+2m_1-1})).\hspace{5.9cm}(6.14)\end{eqnarray*}
 Moreover,
\begin{eqnarray*}&&(E_{n_1,0}-E_{0,n+n_1})\left[
\left(\sum_{i=0}^\infty\frac{(-2)^ix_0^{2i+1}{\cal
D}^i}{(2i+1)!}\right)(\eta^{m_2}(x_{n_1}^{-k+2m_2+1}))\right]
\\ &=&\left(\sum_{i=0}^\infty\frac{im_2(-k+2m_2+1)
(-2)^ix_0^{2i}{\cal
D}^{i-1}}{(2i)!}\right)(\eta^{m_2-1}(x_{n_1}^{-k+2m_2}))
\\ & &+\left(\sum_{i=0}^\infty\frac{(-k+2m_2+1)
(-2)^ix_0^{2i}{\cal
D}^i}{(2i)!}\right)(\eta^{m_2}(x_{n_1}^{-k+2m_2}))
\hspace{7cm}\end{eqnarray*}\begin{eqnarray*} &
&-m_2(-k+2m_2+1)\left(\sum_{i=0}^\infty\frac{(-2)^ix_0^{2i+2}{\cal
D}^i}{(2i+1)!}\right)(\eta^{m_2-1}(x_{n_1}^{-k+2m_2}))
\\ &=&(-k+2m_2+1)\left(\sum_{i=0}^\infty\frac{
(-2)^ix_0^{2i}{\cal
D}^i}{(2i)!}\right)(\eta^{m_2}(x_{n_1}^{-k+2m_2})).
\hspace{4.7cm}(6.15)\end{eqnarray*} Note $k\leq 0$ by our
assumption. Using (6.9), (6.10), (6.14), (6.15) and induction, we
obtain $x_{n_1}^{-k}\in N$.

Observe
$$(E_{n+n_1,0}-E_{0,n_1})|_{{\cal B}'}=x_0x_{n_1}+y_{n_1}\ptl_{x_0}\eqno(6.16)$$
by (1.19) and (1.20). Then
$$(E_{n+n_1,0}-E_{0,n_1})^m(x_{n_1}^{-k})=x_0^mx_{n_1}^{-k+m}+P_m\in N,\eqno(6.17)$$
where the degree of $P_m$ with respect to $x_0$ is less than $m$.
For any $f\in {\cal H}_{\la k-2m\ra}$ and $g\in{\cal H}_{\la
k-2m-1\ra}$ , (3.44) and (5.2) says that \begin{eqnarray*}\qquad &
&\left(\sum_{i=0}^\infty\frac{(-2)^ix_0^{2i}{\cal
D}^i}{(2i)!}\right)(\eta^m(f))\\
&=&\sum_{i=0}^m\frac{2^ix_0^{2i}\prod_{r=1}^i(m-r)(m-k+n_1-n_2+1+r)}{(2i)!}\eta^{m-i}(f)
\hspace{3cm}(6.18)\end{eqnarray*} and
\begin{eqnarray*}\qquad
& &\left(\sum_{i=0}^\infty\frac{(-2)^ix_0^{2i+1}{\cal
D}^i}{(2i)!}\right)(\eta^m(g))\\
&=&\sum_{i=0}^m\frac{2^ix_0^{2i+1}\prod_{r=1}^i(m-r)(m-k+n_1-n_2+2+r)}{(2i+1)!}\eta^{m-i}(g).
\hspace{2.5cm}(6.19)\end{eqnarray*} This shows that if $x_0^m$ is
the highest $x_0$-power of a nonzero element in ${\cal H}'_{\la
k\ra}$, then its coefficient must be in  ${\cal H}_{\la k-m\ra}$ by
(6.8).

On the other hand, (6.17) implies that
$$\mbox{the coefficients of}\;x_0^m\;\mbox{in}\;U(o(2n,\mbb{F}))[(E_{n+n_1,0}-E_{0,n_1})^m(x_{n_1}^{-k})] ={\cal
H}_{\la k-m\ra},\eqno(6.20)$$ because it is an irreducible
$o(2n,\mbb{F})$-module by Theorem 5.1. By induction on $m$, we can
prove
$${\cal H}'_{\la k\ra}\subset
\sum_{m=0}^\infty
U(o(2n,\mbb{F}))[(E_{n+n_1,0}-E_{0,n_1})^m(x_{n_1}^{-k})]\subset
N.\eqno(6.21)$$ Thus $N={\cal H}'_{\la k\ra}$. This shows that
${\cal H}'_{\la k\ra}$ is irreducible. Since the bilinear form
$(\cdot|\cdot)$ restricted to ${\cal H}_{\la k\ra}\subset {\cal
H}_{\la k\ra}'$ is nondegenerate, the irreducibility of ${\cal
H}'_{\la k\ra}$ implies that the symmetric bilinear form
$(\cdot|\cdot)$ restricted to ${\cal H}_{\la k\ra}'$ is
nondegenerate.

Next want to prove
$$\left({\cal H}_{\la k\ra}'|{\cal H}_{\la k'\ra}'\right)=\{0\}\qquad\for\;\;n_1-n_2+1-\dlt_{n_1,n_2}
\geq k,k'\in\mbb{Z}\;\mbox{such that}\;\;k\neq k'. \eqno(6.22)$$ For
any $f\in {\cal H}_{\la k-2m\ra}$ and $f'\in {\cal H}_{\la
k'-2m'\ra}$, (3.64). (5.2), (6.6) and (6.18) yield
\begin{eqnarray*}
& &\left(\left(\sum_{i=0}^\infty\frac{(-2)^ix_0^{2i}{\cal
D}^i}{(2i)!}\right)(\eta^{2m}(f))|\left(\sum_{r=0}^\infty\frac{(-2)^rx_0^{2r}{\cal
D}^r}{(2r)!}\right)(\eta^{2m'}(f'))\right)
\\
&=&
\sum_{i=0}^m\frac{2^{2i}}{(2i)!}[\prod_{s=1}^i(m-s)(m-k+n_1-n_2+1+s)]\\&
&\times[\prod_{s'=1}^i(m'-s')(m'-k'+n_1-n_2+1+s')]
(\eta^{m-i}(f)|\eta^{m'-i}(f))\\ &=&
0\qquad\mbox{if}\;\;(m,k-2m)\neq
(m',k'-2m').\hspace{7.1cm}(6.23)\end{eqnarray*} Let $g\in{\cal
H}_{\la k-2m-1\ra}$ and $g'\in{\cal H}_{\la k'-2m'-1\ra}$. By
(3.60). (5.2), (6.6) and (6.19), we have
\begin{eqnarray*}
& &\left(\left(\sum_{i=0}^\infty\frac{(-2)^ix_0^{2i+1}{\cal
D}^i}{(2i+1)!}\right)(\eta^{2m+1}(g))|\left(\sum_{r=0}^\infty\frac{(-2)^rx_0^{2r+1}{\cal
D}^r}{(2r+1)!}\right)(\eta^{2m'+1}(g'))\right)
\\&=&
\sum_{i=0}^m\frac{2^{2i}}{(2i+1)!}[\prod_{s=1}^i(m-s)(m-k+n_1-n_2+2+s)]\\&
&\times[\prod_{s'=1}^i(m'-s')(m'-k'+n_1-n_2+2+s')]
(\eta^{m-i}(g)|\eta^{m'-i}(g'))\\ &=&
0\qquad\mbox{if}\;\;(m,k-2m-1)\neq
(m',k'-2m'-1).\hspace{5.6cm}(6.24)\end{eqnarray*} Since
$(x^{2i}_0|x_0^{2i'+1})=0$ for $i,i'\in\mbb{N}$, the elements of the
form (6.18) are orthogonal to those of the form (6.19). Hence (6.22)
holds by (6.8).

For $g\in{\cal H}'_{\la k\ra}$ and $m\in\mbb{N}+1$,
$${\cal
D}'[(\eta')^m(g)]=2m[2(k+n_2-n_1+m-1)+1](\eta')^{m-1}(g)\eqno(6.25)$$
by (3.44) and the facts ${\cal D}'=\ptl_{x_0}^2-2{\cal D}$ and its
dual $\eta'=x_0^2+2\eta$. This shows that
$$((\eta')^m({\cal H}'_{\la k\ra})|(\eta')^{m'}({\cal H}'_{\la
k'\ra}))=\{0\}\qquad\mbox{if}\;\;(m,k)\neq (m',k')\eqno(6.26)$$ for
$n_1-n_2+1-\dlt_{n_1,n_2}\geq k,k'\in\mbb{Z}$ and $m,m'\in\mbb{N}$
by (6.7). Moreover, the symmetric bilinear form $(\cdot|\cdot)$
restricted to $(\eta')^m({\cal H}_{\la k\ra}')$ is nondegenerate.

 Fix $n_1-n_2+1-\dlt_{n_1,n_2}\geq k\in\mbb{Z}$. Denote
$$\hat{\cal B}'_{\la k\ra}=\bigoplus_{i=0}^\infty(\eta')^i({\cal H}'_{\la k-2i\ra})
 \eqno(6.27)$$
Then the symmetric bilinear form $(\cdot|\cdot)$ restricted to
$\hat{\cal B}'_{\la k\ra}$ is nondegenerate. Thus
$${\cal B}'_{\la k\ra}=\hat{\cal B}'_{\la k\ra}\oplus (\hat{\cal B}'_{\la
k\ra})^{\perp}\bigcap {\cal B}'_{\la k\ra}.\eqno(6.28)$$ According
to Lemma 3.2,  $(\hat{\cal B}'_{\la k\ra})^{\perp}\bigcap {\cal
B}'_{\la k\ra}$ is an $o(2n+1,\mbb{F})$-module. Assume $(\hat{\cal
B}'_{\la k\ra})^{\perp}\bigcap {\cal B}'_{\la k\ra} \neq\{0\}$.
 By (5.2), (5.3), (5.8)-(5.13), (5.23)
and (1.23), there exists a nonzero element in $(\hat{\cal B}'_{\la
k\ra})^{\perp}\bigcap {\cal B}'_{\la k\ra}$ of the form:
$$f=\sum_{i=0}^ma_ix_0^{2i}(2\eta)^{m-i}(x_{n_1}^{-k+2m})\eqno(6.29)$$
or
$$g=\sum_{i=0}^mb_ix_0^{2i+1}(2\eta)^{m-i}(x_{n_1}^{-k+2m+1})\eqno(6.30)$$
for some $m\in\mbb{N}+1$. Moreover, we assume that the exponent of
$x_{n_1}$ is minimal.

If (6.29) holds, then (6.11)) and (6.13) give
\begin{eqnarray*}\qquad&&(E_{n_1,0}-E_{0,n+n_1})(f)=(\ptl_{x_0}\ptl_{x_{n_1}}-x_0\ptl_{y_{n_1}})(f)
\\ &=&\sum_{i=1}^m
2i(-k+2m)a_ix_0^{2i-1}(2\eta)^{m-i}(x_{n_1}^{-k+2m-1})\\
& &-\sum_{i=0}^{m-1}2(m-i)(-k+2m)a_ix_0^{2i+1}
(2\eta)^{m-i-1}(x_{n_1}^{-k+2m-1})\\
&=&2(-k+2m)\sum_{i=0}^{m-1}[(i+1)a_{i+1}-(m-i)a_i]x_0^{2i+1}
(2\eta)^{m-i-1}(x_{n_1}^{-k+2m-1})\\
&=&0\hspace{12.8cm}(6.31)\end{eqnarray*}by the minimality of the
exponent of $x_{n_1}$, equivalently
$$(i+1)a_{i+1}=(m-i)a_i\qquad\for\;\;i\in\ol{0,m-1}.\eqno(6.32)$$
Thus
$$a_i=a_0{m\choose i}\qquad\for\;\;i\in\ol{0,m}.\eqno(6.33)$$
So
$$f=\sum_{i=0}^ma_0{m\choose i}x_0^{2i}(2\eta)^{m-i}(x_{n_1}^{-k+2m})
=a_0(\eta')^m(x_{n_1}^{-k+2m})\in \hat{\cal B}'_{\la
k\ra},\eqno(6.34)$$ which contradicts (6.28).

Suppose that (6.30) holds. Note $x_0x_{n_1}^{-k+2m+1}\in{\cal
H}'_{\la k-2m\ra}$
 by (1.23). Expressions (6.11) and (6.13) deduce
 \begin{eqnarray*}\qquad&&(E_{n_1,0}-E_{0,n+n_1})(g)
 =(\ptl_{x_0}\ptl_{x_{n_1}}-x_0\ptl_{y_{n_1}})(g)
\\ &=&\sum_{i=0}^m
(2i+1)(-k+2m+1)b_ix_0^{2i}(2\eta)^{m-i}(x_{n_1}^{-k+2m})\\
& &-\sum_{i=0}^{m-1}2(m-i)(-k+2m+1)b_ix_0^{2i+2}
(2\eta)^{m-i-1}(x_{n_1}^{-k+2m})
\\
&=&(-k+2m+1)\big\{\sum_{i=0}^{m-1}[(2i+3)b_{i+1}-2(m-i)b_i]x_0^{2i+2}
(2\eta)^{m-i-1}(x_{n_1}^{-k+2m})\\
& &+b_0(2\eta)^m(x_{n_1}^{-k+2m})\big\}
=0\hspace{8.7cm}(6.35)\end{eqnarray*}by the minimality of the
exponent of $x_{n_1}$, equivalently
$$b_0=0,\;\;(2i+3)b_{i+1}=2(m-i)b_i\qquad\for\;\;i\in\ol{0,m-1}.\eqno(6.36)$$
Thus $b_i=0$ for $i\in\ol{0,m},$ that is, $g=0$. This contradicts
our choice of nonzero element. Hence $(\hat{\cal B}'_{\la
k\ra})^{\perp}\bigcap {\cal B}'_{\la k\ra}=\{0\}$. Then (6.28) gives
(6.4). Furthermore, (6.5) is obtained by Lemma 3.1 with
$T_1=\ptl_{x_0}^2,\;T_1^-=\int_{(x_0)}^{(2)}$
 (cf. (2.6) and (2.7)) and $T_2=2{\cal D}$.

When $n_1=n_2$, an expression of ${\cal H}_{\la k\ra}'$ can be
obtained via (5.3), (5.27)-(5.31), (6.8), (6.18) and (6.19). In
particular, when $n_1=n_2=n$, the $({\cal G},{\cal K})$-module
structure is given by
 \begin{eqnarray*}\qquad{\cal H}'_{\la
-k\ra}&=&\bigoplus_{m,r=0}^\infty\left(\sum_{i=0}^\infty\frac{(-2)^ix_0^{2i}{\cal
D}^i}{(2i)!}\right)(\eta^r({\cal H}_{\la -k-2r-m,m\ra}))\\ & &\oplus
\bigoplus_{l,s=0}^\infty\left(\sum_{i=0}^\infty\frac{(-2)^ix_0^{2i+1}{\cal
D}^i}{(2i+1)!}\right)(\eta^s({\cal H}_{\la
-k-2s-1-l,l\ra})),\hspace{3.6cm}(6.37)\end{eqnarray*} where ${\cal
H}_{\la -m_1-m_2,m_2\ra}$ given in (4.42). $\qquad\Box$

\section{Noncanonical Representations of $sp(2n,\mbb{F})$}

In this section, we use the results in Sections 3 and 4 to study
noncanonical polynomial representation of $sp(2n,\mbb{F})$.

Recall the symplectic Lie
algebra\begin{eqnarray*}\hspace{1cm}sp(2n,\mbb{F})&=&
\sum_{i,j=1}^n\mbb{F}(E_{i,j}-E_{n+j,n+i})+\sum_{i=1}^n(\mbb{F}E_{i,n+i}+\mbb{F}E_{n+i,i})\\
& &+\sum_{1\leq i<j\leq n
}[\mbb{F}(E_{i,n+j}+E_{n+j,i})+\mbb{F}(E_{n+i,j}+E_{n+j,i})].\hspace{2.4cm}(7.1)\end{eqnarray*}
Again we take the Cartan subalgebra
$H=\sum_{i=1}^n\mbb{F}(E_{i,i}-E_{n+i,n+i})$ and the subspace
spanned by positive root vectors
$$sp(2n,\mbb{F})_+=\sum_{1\leq i<j\leq n}[\mbb{F}(E_{i,j}-E_{n+j,n+i})
+\mbb{F}(E_{i,n+j}+E_{n+j,i})]+\sum_{i=1}^n\mbb{F}E_{i,n+i}.\eqno(7.2)$$
Fix $1\leq n_1\leq n_2\leq n$. The noncanonical oscillator
representation of $sp(2n,\mbb{F})$ on ${\cal
B}=\mbb{F}[x_1,...,x_n,y_1,...,y_n]$ is defined via (1.14)-(1.16).
Recall ${\cal K}=\sum_{i,j=1}^n\mbb{F}(E_{i,j}-E_{n+j,n+i})$. \psp

{\bf Theorem 7.1}. {\it Let $k\in\mbb{Z}$. If $n_1<n_2$ or $k\neq
0$, the subspace ${\cal B}_{\la k\ra}$ (cf. (1.17)) is an
irreducible $sp(2n,\mbb{F})$-module. Moreover, it is a
highest-weight module only if $n_2=n$, in which case for
$m\in\mbb{N}$, $x_{n_1}^{-m}$ is a highest-weight vector of ${\cal
B}_{\la -m\ra}$ with weight $-m\lmd_{n_1-1}+(m-1)\lmd_{n_1}$,
$x_{n_1+1}^{m+1}$ is a highest-weight vector of ${\cal B}_{\la
m+1\ra}$ with weight
$-(m+2)\lmd_{n_1}+(m+1)\lmd_{n_1+1}+(m+1)\dlt_{n_1,n-1}\lmd_n$ if
$n_1<n$ and $y_n^{m+1}$ is a highest-weight vector of ${\cal B}_{\la
m+1\ra}$ with weight $(m+1)\lmd_{n-1}-2(m+1)\lmd_n$ when $n_1=n$.

When $n_1=n_2$, the subspace ${\cal B}_{\la 0\ra}$ is a direct sum
of two irreducible $sp(2n,\mbb{F})$-submodules. If $n_1=n_2=n$, they
are highest-weight modules with a highest-weight vector $1$ of
weight $-2\lmd_n$ and with a highest-weight vector
$x_{n-1}y_n-x_ny_{n-1}$ of weight $\dlt_{n,2}\lmd_{n-2}-4\lmd_n$,
respectively.  If $n_1=n_2=n$, all the irreducible modules are of
$({\cal G},{\cal K})$-type.}

{\it Proof}. Recall that we embed $sl(n,\mbb{F})$ into
$sp(2n,\mbb{F})$ via $E_{i,j}\mapsto E_{i,j}-E_{n+j,n_i}$. Moreover,
${\cal B}$ is nilpotent with respect to $sl(n,\mbb{F})_+$ (cf.
(2.30)) and
$$\eta=\sum_{i=1}^{n_1}y_i\ptl_{x_i}+\sum_{r=n_1+1}^{n_2}x_ry_r+\sum_{s=n_2+1}^n
x_s\ptl_{y_s}.\eqno(7.3)$$ Note
$$(E_{i,n+j}+E_{j,n+i})|_{\cal
B}=\ptl_{x_i}\ptl_{y_j}+\ptl_{x_j}\ptl_{y_i},\;\;(E_{i,n+r}+E_{r,n+i})|_{\cal
B}=\ptl_{x_i}\ptl_{y_r}+x_r\ptl_{y_i},\eqno(7.4)$$
$$(E_{r,n+s}+E_{s,n+r})|_{\cal
B}=x_r\ptl_{y_s}+x_s\ptl_{y_r}\eqno(7.5)$$ for $i,j\in\ol{1,n_1}$
and $r,s\in\ol{n_1+1,n_2}$ by (1.15).
 Moreover,
$$(E_{i,j}-E_{n+j,n+i})|_{\cal
B}=-x_j\ptl_{x_i}-y_j\ptl_{y_i}-\dlt_{i,j}\eqno(7.6)$$ and
$$(E_{i,r}-E_{n+r,n+i})|_{\cal
B}=\ptl_{x_i}\ptl_{x_r}-y_r\ptl_{y_i}\eqno(7.7)$$ for
$i,j\in\ol{1,n_1}$ and $r\in\ol{n_1+1,n_2}$ by (1.7), (1.8) and
(1.14). We will process our arguments in two steps.\psp

{\it Step 1}. $n_2=n$.\psp

Under the assumption, ${\cal B}$ is nilpotent with respect to
$sp(2n,\mbb{F})_+$ by (7.4)-(7.7).

First we assume $n_1+1<n$. According to (3.37), the nonzero weight
vectors in
$$\mbox{Span}\{\eta^{m_3}(x_i^{m_1}y_n^{m_2})\mid
m_r\in\mbb{N};i=n_1,n_1+1\}\eqno(7.8)$$ are all the singular vectors
of $sl(n,\mbb{F})$ in ${\cal B}$. The singular vectors of
$sp(2n,\mbb{F})$ in ${\cal B}$ must be among them. Moreover, the
subalgebra $sp(2n,\mbb{F})_+$ is generated by $sl(n,\mbb{F})_+$ and
$E_{n,2n}$. According to (7.5), $E_{n,2n}|_{\cal B}=x_n\ptl_{y_n}$.
Hence
$$E_{n,2n}(\eta^{m_3}(x_i^{m_1}y_n^{m_2}))=x_n[
m_3x_n\eta^{m_3-1}(x_i^{m_1}y_n^{m_2})+m_2\eta^{m_3}(x_i^{m_1}y_n^{m_2-1})]\eqno(7.9)$$
for $i=n_1,n_1+1$ by (7.3). Considering weights, we conclude that
the vectors $\{x_{n_1}^{m},x_{n_1+1}^{m+1}\mid m\in\mbb{N}$ are all
the singular vectors of $sp(2n,\mbb{F})$ in ${\cal B}$. Furthermore,
$$x_{n_1}^{m}\in {\cal B}_{\la
-m\ra}\;\;\mbox{and}\;\;x_{n_1+1}^{m+1}\in {\cal B}_{\la
m+1\ra}\qquad\for\;\;m\in\mbb{N}.\eqno(7.10)$$ Thus each ${\cal
B}_{\la k\ra}$ has a unique non-isotropic singular vector for
$k\in\mbb{Z}$. By Lemma 3.3, all ${\cal B}_{\la k\ra}$ with
$k\in\mbb{Z}$ are irreducible highest-weight
$sp(2n,\mbb{F})$-submodules.

Consider the case $n_1+1=n$.  According to (3.112), the nonzero
weight vectors in
$$\mbox{Span}\{\eta^{m_2}(x_{n-1}^{m_1}y_n^{m_3}),
x_n^{m_1}y_n^{m_2},\eta^{m_1+m_2}(x_{n-1}^{m_2}y_n^{m_3-m_1})\mid
m_i\in\mbb{N}\}\eqno(7.11)$$ are all the singular vectors of
$sl(n,\mbb{F})$ in ${\cal B}$. Recall $E_{n,2n}|_{\cal
B}=x_n\ptl_{y_n}$. We have
$$E_{n,2n}(x_n^{m_1}y_n^{m_2})=m_2x_n^{m_1+1}y_n^{m_2-1}.\eqno(7.12)$$
By (7.11) and considering weights, we again conclude that the
vectors $\{x_{n-1}^{m},x_n^{m+1}\mid m\in\mbb{N}$ are all the
singular vectors of $sp(2n,\mbb{F})$ in ${\cal B}$. Again all ${\cal
B}_{\la k\ra}$ with $k\in\mbb{Z}$ are irreducible highest-weight
$sp(2n,\mbb{F})$-submodules.

Suppose $n_1=n$. By (7.4), we have $E_{n,2n}=\ptl_{x_n}\ptl_{y_n}$
in this case. According to (4.31), the nonzero weight vectors in
$$\mbox{Span}\{x_n^{m_1}y_n^{m_2}\zeta_1^{m_3}\mid
m_i\in\mbb{N}\}\eqno(7.13)$$
 are all the singular vectors of
$sl(n,\mbb{F})$ in ${\cal B}$,
 where $\zeta_1=x_{n-1}y_n-x_ny_{n-1}$
in this case.
\begin{eqnarray*}
& &E_{n,2n}(x_n^{m_1}y_n^{m_2}\zeta_1^{m_3})\\ &=&
m_1m_2x_n^{m_1-1}y_n^{m_2-1}\zeta_1^{m_3}+m_1m_3x_{n-1}x_n^{m_1-1}y_n^{m_2}\zeta_1^{m_3-1}
\\
&&-m_2m_3y_{n-1}x_n^{m_1}y_n^{m_2-1}\zeta_1^{m_3-1}
-m_3(m_3-1)x_{n-1}y_{n-1}x_n^{m_1}y_n^{m_2}\zeta_1^{m_3-2}
\\
&=&m_1(m_2+m_3)x_n^{m_1-1}y_n^{m_2-1}\zeta_1^{m_3}
+m_3(m_1-m_2-m_3+1)y_{n-1}x_n^{m_1}y_n^{m_2-1}\zeta_1^{m_3-1}
\\& &-m_3(m_3-1)y_{n-1}^2x_n^{m_1+1}y_n^{m_2-1}\zeta_1^{m_3-2}
.\hspace{7.5cm}(7.14)\end{eqnarray*}
 Considering weights, we again conclude that the
vectors $\{x_n^m,y_n^{m+1},\zeta_1\mid m\in\mbb{N}\}$ are all the
singular vectors of $sp(2n,\mbb{F})$ in ${\cal B}$. Moreover,
$$x_n^{m}\in {\cal B}_{\la -m\ra},\;\;\zeta_1\in{\cal B}_{\la 0\ra}\;\;\mbox{and}\;\;y_n^{m+1}
\in{\cal B}_{\la m+1\ra}\qquad\for\;\;m\in\mbb{N}.\eqno(7.15)$$ Thus
each ${\cal B}_{\la k\ra}$ with $k\neq 0$ has a unique non-isotropic
singular vector for $k\in\mbb{Z}$. By Lemma 3.3, all ${\cal B}_{\la
k\ra}$ with $0\neq k\in\mbb{Z}$ are irreducible highest-weight
$sp(2n,\mbb{F})$-submodules.

Set
$${\cal B}_{\la
0,1\ra}=\mbox{Span}\{[\prod_{1\leq r\leq s\leq
n}(x_ry_s+x_sy_r)^{l_{r,s}}]\mid l_{r,s}\in\mbb{N}\}\eqno(7.18)$$and
$${\cal B}_{\la 0,2\ra}=\mbox{Span}\{[\prod_{1\leq r\leq s\leq
n}(x_ry_s+x_sy_r)^{l_{r,s}}](x_py_q-x_qy_p)\mid
l_{r,s}\in\mbb{N};1\leq p<q\leq n\}.\eqno(7.19)$$ Let $${\cal
G}'=\sum_{1\leq r\leq s\leq
n}\mbb{F}(E_{n+s,r}+E_{n+r,s})\eqno(7.20)$$
 and
$$\hat{\cal
G}=\sum_{i,j=1}^n\mbb{F}(E_{i,j}-E_{n+j,n+i})+\sum_{1\leq r\leq
s\leq n}\mbb{F}(E_{r,n+s}+E_{s,n+r}).\eqno(7.21)$$
 Then ${\cal G}'$ and
$\hat{\cal G}$ are Lie subalgebras of $sp(2n,\mbb{F})$ and
$sp(2n,\mbb{F})={\cal G}'\oplus \hat{\cal G}.$ By PBW Theorem
$$U(sp(2n,\mbb{F}))=U({\cal G}')U(\hat{\cal G}).\eqno(7.22)$$
Note
$$(E_{n+s,r}+E_{n+r,s})_{\cal B}=-(x_ry_s+x_sy_r)\qquad\for\;\;r,s\in\ol{1,n}\eqno(7.23)$$
by (1.16). According to (7.4), (7.6) and (7.23),
$${\cal B}_{\la
0,1\ra}=U({\cal G}')(1)=U(sp(2n,\mbb{F}))(1)\eqno(7.24)$$ and
$${\cal B}_{\la
0,2\ra}=\sum_{1\leq p<q\leq n}U({\cal
G}')(x_py_q-x_qy_p)=U(sp(2n,\mbb{F}))(\zeta_1)\eqno(7.25)$$ are
$sp(2n,\mbb{F})$-submodules.

It is obvious, $1\not\in {\cal B}_{\la 0,2\ra}$. On the other hand,
$({\cal B}_{\la 0,1\ra}|x_{n-1}y_n-x_ny_{n-1})=\{0\}$. Hence
$x_{n-1}y_n-x_ny_{n-1}\not \in {\cal B}_{\la 0,1\ra}$. Thus ${\cal
B}_{\la 0,1\ra}$ and ${\cal B}_{\la 0,0\ra}$ have a unique
non-isotropic singular vector. By Lemma 3.3, they are irreducible.
Since $1$ and $x_{n-1}y_n-x_ny_{n-1}$ are the only singular vectors
in ${\cal B}_{\la 0\ra}$ which is nilpotent with respect to
$sp(2n,\mbb{F})_+$, Lemma 2.3 yields
$${\cal B}_{\la 0\ra}={\cal B}_{\la 0,1\ra}\oplus {\cal B}_{\la
0,2\ra}\eqno(7.26)$$ by the similar arguments as those from (3.67)
to (3.69).\psp

{\it Step 2}. $n_2<n$.\psp

We set
\begin{eqnarray*}\hspace{1cm}{\cal G}_1&=&
\sum_{i,j=1}^{n_2}\mbb{F}(E_{i,j}-E_{n+j,n+i})+\sum_{i=1}^{n_2}(\mbb{F}E_{i,n+i}+\mbb{F}E_{n+i,i})\\
& &+\sum_{1\leq i<j\leq n_2
}[\mbb{F}(E_{i,n+j}+E_{n+j,i})+\mbb{F}(E_{n+i,j}+E_{n+j,i})]\hspace{3.6cm}(7.27)\end{eqnarray*}
and
\begin{eqnarray*}\hspace{1cm}{\cal G}_2&=&
\sum_{i,j=n_1+1}^n\mbb{F}(E_{i,j}-E_{n+j,n+i})+\sum_{i=n_1+1}^n(\mbb{F}E_{i,n+i}+\mbb{F}E_{n+i,i})\\
& &+\sum_{n_1+1\leq i<j\leq n
}[\mbb{F}(E_{i,n+j}+E_{n+j,i})+\mbb{F}(E_{n+i,j}+E_{n+j,i})].\hspace{3cm}(7.28)\end{eqnarray*}
Then ${\cal G}_1=sp(2n_2,\mbb{F})$ and ${\cal G}_2\cong
sp(2(n-n_1),\mbb{F})$ are Lie subalgebras of $sp(2n,\mbb{F})$.
Denote
$${\cal M}^1=\mbb{F}[x_1,...,x_{n_2},y_1,...,y_{n_2}],\qquad
{\cal
M}^2=\mbb{F}[x_{n_1+1},...,x_n,y_{n_1+1},...,y_n].\eqno(7.29)$$
Observe that ${\cal M}^1$ is exactly the ${\cal G}_1$-module as
${\cal B}$ in Step 1 with $n\rightarrow n_2$ and ${\cal M}^2$ is
exactly the ${\cal G}_1$-module as ${\cal B}$ in Step 1 with
$n_1=n_2$ and $n\rightarrow n-n_1$. Moreover, we set
$${\cal M}^3=\mbb{F}[x_1,...,x_{n_1},y_1,...,y_{n_1}],\qquad
{\cal
M}^4=\mbb{F}[x_{n_1+1},...,x_{n_2},y_{n_1+1},...,y_{n_2}].\eqno(7.30)$$
Let
$${\cal M}^i_{\la k\ra}={\cal M}^i\bigcap{\cal B}_{\la
k\ra}\qquad\for\;\;i\in\ol{1,4},\;k\in\mbb{Z}.\eqno(7.31)$$ Then
$${\cal M}^1_{\la k\ra}=\bigoplus_{r\in\mbb{Z}}{\cal M}^3_{\la r\ra}{\cal M}^4_{\la
k-r\ra}\qquad\for\;\;k\in\mbb{Z}.\eqno(7.32)$$ Next we prove the
theorem case by case.\psp

{\it Case 1}. $n_1+1<n_2$\psp

 According to (3.36), the
nonzero weight vectors in
$$\mbox{Span}\{\eta^{m_3}(x_i^{m_1}y_j^{m_2})\mid
m_r\in\mbb{N};i=n_1,n_1+1;j=n_2,n_2+1\}\eqno(7.33)$$ are all the
singular vectors of $sl(n,\mbb{F})$ in ${\cal B}$. Fix
$k\in\mbb{N}$. Then the singular vectors of $sl(n,\mbb{F})$  in
${\cal B}_{\la-k\ra}$ are
\begin{eqnarray*}\qquad\qquad& &\{\eta^{m_3}(x_{n_1}^{k+m_2+2m_3}y_{n_2}^{m_2}),
\eta^{m_3}(x_{n_1+1}^{m_1}y_{n_2+1}^{k+m_1+2m_3}),\\
& &\eta^{m_3}(x_{n_1}^{m_4}y_{n_2+1}^{m_5}) \mid
m_i\in\mbb{N};m_4+m_5-2m_3=k\}.\hspace{4.2cm}(7.34)\end{eqnarray*}
Let $M$ be a nonzero $sp(2n,\mbb{F})$-submodule of ${\cal
B}_{\la-k\ra}$. Then $M$ contains a singular of $sl(n,\mbb{F})$.
Suppose some $\eta^{m_3}(x_{n_1}^{k+m_2+2m_3}y_{n_2}^{m_2})\in M$.
We have $E_{n_1,n+n_1}|_{\cal B}=\ptl_{x_{n_1}}\ptl_{y_{n_1}}$ and
$$E_{n_1,n+n_1}^{m_3}[\eta^{m_3}(x_{n_1}^{k+m_2+2m_3}y_{n_2}^{m_2})]=m_3![\prod_{r=1}^{2m_3}(k+m_2+r)]
x_{n_1}^{k+m_2}y_{n_2}^{m_2}\in M\eqno(7.35)$$ by (7.3) and (7.4).
Moreover, $(E_{n_1,n+n_2}+E_{n_2,n+n_1})|_{\cal
B}=\ptl_{x_{n_1}}\ptl_{y_{n_2}}+x_{n_2}\ptl_{y_{n_1}}$ and
$$(E_{n_1,n+n_2}+E_{n_2,n+n_1})^{m_2}(x_{n_1}^{k+m_2}y_{n_2}^{m_2})=m_2!
[\prod_{r=1}^{m_2}(k+r)] x_{n_1}^k\in M\eqno(7.36)$$ by (7.4). Thus
$$x_{n_1}^k\in M.\eqno(7.37)$$

Assume some $\eta^{m_3}(x_{n_1+1}^{m_1}y_{n_2+1}^{k+m_1+2m_3})\in
M$. According to (1.16),
 $$(E_{n+i,j}+E_{n+j,i})|_{\cal B}=\ptl_{x_i}\ptl_{y_j}+\ptl_{x_j}\ptl_{y_i}\qquad\for\;\;i\in\ol{n_2+1,n}.\eqno(7.38)$$
So
$$E_{n+n_2+1,n_2+1}^{m_3}[\eta^{m_3}(x_{n_1+1}^{m_1}y_{n_2+1}^{k+m_1+2m_3})]=
m_3![\prod_{r=1}^{2m_3}(k+m_1+r)]x_{n_1+1}^{m_1}y_{n_2+1}^{k+m_1}\in
M.\eqno(7.39)$$ Moreover,
$$(E_{n+n_2+1,n_1+1}+E_{n+n_1+1,n_2+1})|_{\cal B}=\ptl_{x_{n_1+1}}
\ptl_{y_{n_2+1}}+y_{n_1+1}\ptl_{x_{n_2+1}}\eqno(7.40)$$ by (1.16).
Hence
$$(E_{n+n_2+1,n_1+1}+E_{n+n_1+1,n_2+1})^{m_1}(x_{n_1+1}^{m_1}y_{n_2+1}^{k+m_1})=m_1!
[\prod_{r=1}^{m_1}(k+r)] y_{n_2+1}^k\in M.\eqno(7.41)$$ Furthermore,
$$(E_{n+n_2+1,n_1}+E_{n+n_1,n_2+1})|_{\cal
B}=-x_{n_1}\ptl_{y_{n_2+1}}+y_{n_1}\ptl_{x_{n_2+1}}\eqno(7.42)$$ by
(1.16). Thus
$$(E_{n+n_2+1,n_1}+E_{n+n_1,n_2+1})^k(y_{n_2+1}^k)=(-1)^kk!x_{n_1}^k\in
M.\eqno(7.43)$$ Thus (7.37) holds again.

Consider $\eta^{m_3}(x_{n_1}^{m_4}y_{n_2+1}^{m_5})$ for some
$m_3,m_3,m_4\in\mbb{N}$ such that $m_4+m_5-2m_3=k$.  Note that
$E_{n_1+1,n+n_1+1}|_{\cal B}=x_{n_1+1}\ptl_{y_{n_1+1}}$ by (7.5) and
$$E_{n_1+1,n+n_1+1}^{m_3}[\eta^{m_3}(x_{n_1}^{m_4}y_{n_2+1}^{m_5})]
=m_3!x_{n_1+1}^{2m_3}x_{n_1}^{m_4}y_{n_2+1}^{m_5}\in M.\eqno(7.44)$$
There exists $r_1,r_2\in\mbb{N}$ such that $r_1+r_2=2m_3$ and
$r_1\leq m_4,\;r_2\leq m_5$. Moreover,
$$(E_{n_1,n_1+1}-E_{n+n_1+1,n+n_1})|_{\cal
B}=\ptl_{x_{n_1}}\ptl_{x_{n_1+1}}-y_{n_1+1}\ptl_{y_{n_1}}\eqno(7.45)$$
by (1.7), (1.8) and  (1.14). Moreover, (7.40) and (7.45) yield
\begin{eqnarray*}\qquad &
&(E_{n_1,n_1+1}-E_{n+n_1+1,n+n_1})^{r_1}(E_{n+n_2+1,n_1+1}+E_{n+n_1+1,n_2+1})^{r_2}
(x_{n_1+1}^{2m_3}x_{n_1}^{m_4}y_{n_2+1}^{m_5})\\
&=&(2m_3)![\prod_{s_1=0}^{r_1-1}(m_4-s_1)]
[\prod_{s_2=0}^{r_2-1}(m_5-s_2)]x_{n_1}^{m_4-r_1}y_{n_2+1}^{m_5-r_2}\in
M.\hspace{3.2cm}(7.46)\end{eqnarray*} Furthermore, (7.42) yields
\begin{eqnarray*}\hspace{2cm}& &(E_{n+n_2+1,n_1}+E_{n+n_1,n_2+1})^{m_5-r_2}
(x_{n_1}^{m_4-r_1}y_{n_2+1}^{m_5-r_2})\\
&=&(-1)^{m_5-r_2}(m_5-r_2)!x_{n_1}^k\in
M.\hspace{6.6cm}(7.47)\end{eqnarray*} Thus we always have
$x_{n_1}^k\in M$.

Note that ${\cal M}^1_{\la-k\ra}\ni x_{n_1}^k$ is an irreducible
${\cal G}_1$-module (cf. (7.27) and (7.29)) by Step 1. So
$${\cal M}^1_{\la-k\ra}\subset M.\eqno(7.48)$$
Let $r\in\mbb{Z}$. According to (7.32),
$${\cal M}^3_{\la r\ra}{\cal M}^4_{\la-k- r\ra}\subset {\cal M}^1_{\la-k\ra}\subset M.\eqno(7.49)$$
Moreover, ${\cal M}^2_{\la-k-r\ra}\supset{\cal M}^4_{\la-k- r\ra}$
is an irreducible ${\cal G}_2$-module (cf. (7.28) and (7.29)) by
Step 1. Thus
$${\cal M}^3_{\la r\ra}{\cal M}^2_{\la-k- r\ra}=U({\cal G}_2)({\cal M}^3_{\la r\ra}{\cal M}^4_{\la-k-
r\ra})\subset M.\eqno(7.50)$$ Then
$${\cal B}_{\la -k\ra}=\bigoplus_{r\in\mbb{Z}}{\cal M}^3_{\la r\ra}{\cal M}^2_{\la-k-
r\ra}\subset M\eqno(7.51)$$ by (7.29) and (7.30). Therefore,
$M={\cal B}_{\la -k\ra}$, that is, ${\cal B}_{\la -k\ra}$ is an
irreducible $sp(2n,\mbb{F})$-submodule.

Fix $0<k\in\mbb{N}$.  Then the singular vectors of $sl(n,\mbb{F})$
in ${\cal B}_{\la k\ra}$ are
\begin{eqnarray*}\qquad& &\{\eta^{m_2}(x_{n_1+1}^{k+m_1-2m_2}y_{n_2+1}^{m_1}),
\eta^{m_2}(x_{n_1}^{m_1}y_{n_2}^{k+m_1-2m_2}),\eta^{m_3}(x_{n_1+1}^{m_4}y_{n_2}^{m_5})
\\
& &\mid m_i\in\mbb{N};2m_2\leq
k+m_1;m_4+m_5+2m_3=k\}\hspace{5.1cm}(7.52)\end{eqnarray*} by (7.33).
Let $M$ be a nonzero $sp(2n,\mbb{F})$-submodule of ${\cal B}_{\la
k\ra}$. Then $M$ contains a singular of $sl(n,\mbb{F})$. Suppose
some $\eta^{m_2}(x_{n_1+1}^{k+m_1-2m_2}y_{n_2+1}^{m_1})\in M$ with
$2m_2\leq k+m_1$. We have $E_{n_1+1,n+n_1+1}|_{\cal
B}=x_{n_1+1}\ptl_{y_{n_1+1}}$ and
$$E_{n_1+1,n+n_1+1}^{m_2}[\eta^{m_2}(x_{n_1+1}^{k+m_1-2m_2}y_{n_2}^{m_1})]=m_2!
x_{n_1+1}^{k+m_1}y_{n_2+1}^{m_1}\in M\eqno(7.53)$$ by (7.3) and
(7.5). Moreover, (7.40) gives
$$(E_{n+n_2+1,n_1+1}+E_{n+n_1+1,n_2+1})^{m_1}(x_{n_1+1}^{k+m_1}y_{n_2+1}^{m_1})=m_1!
[\prod_{r=1}^{m_1}(k+r)] x_{n_1+1}^k\in M.\eqno(7.54)$$ Thus
$$x_{n_1+1}^k\in M.\eqno(7.55)$$

Assume some $\eta^{m_2}(x_{n_1}^{m_1}y_{n_2}^{k+m_1-2m_2})\in M$
with $2m_2\leq k+m_1$. Observe $E_{n+n_2,n_2}=y_{n_2}\ptl_{x_{n_2}}$
by (1.16). So
$$E_{n+n_2,n_2}^{m_2}[\eta^{m_2}(x_{n_1}^{m_1}y_{n_2}^{k+m_1-2m_2})]=
m_2!x_{n_1}^{m_1}y_{n_2}^{k+m_1}\in M.\eqno(7.56)$$  Moreover, (7.4)
gives that $(E_{n_1,n+n_2}+E_{n_2,n+n_1})|_{\cal
B}=\ptl_{x_{n_1}}\ptl_{y_{n_2}}+x_{n_2}\ptl_{y_{n_1}}$ and
$$(E_{n_1,n+n_2}+E_{n_2,n+n_1})^{m_1}(x_{n_1}^{m_1}y_{n_2}^{k+m_1})=m_1!
[\prod_{r=1}^{m_1}(k+r)] y_{n_2}^k\in M.\eqno(7.57)$$ Furthermore,
(7.5) yields that $(E_{n_1+1,n+n_2}+E_{n_2,n+n_1+1})|_{\cal
B}=x_{n_1+1}\ptl_{y_{n_2}}+x_{n_2}\ptl_{y_{n_1+1}}$ and
$$(E_{n_1+1,n+n_2}+E_{n_2,n+n_1+1})^k(y_{n_2}^k)=k!x_{n_1+1}^k\in
M.\eqno(7.58)$$ Thus (7.55) holds again.

Consider $\eta^{m_3}(x_{n_1+1}^{m_4}y_{n_2}^{m_5})$ for some
$m_3,m_3,m_4\in\mbb{N}$ such that $m_4+m_5+2m_3=k$.  Note
$E_{n_1+1,n+n_1+1}=x_{n_1+1}\ptl_{y_{n_1+1}}$ by (7.5). So
$$E_{n_1+1,n+n_1+1}^{m_3}
[\eta^{m_3}(x_{n_1+1}^{m_4}y_{n_2}^{m_5})]=m_3!x_{n_1+1}^{m_4+2m_3}y_{n_2}^{m_5}\in
M.\eqno(7.59)$$ According to (7.5),
$$(E_{n_1+1,n+n_2}+E_{n_2,n+n_1+1})^{m_5}(x_{n_1+1}^{m_4+2m_3}y_{n_2}^{m_5})=m_5!x_{n_1+1}^k\in
M.\eqno(7.60)$$ Therefore, we always have $x_{n_1+1}^k\in M$.

 Observe that ${\cal M}^2_{\la k\ra}\ni x_{n_1+1}^k$ is an irreducible ${\cal
G}_2$-module (cf. (7.28) and (7.29)) by Step 1. So
$${\cal M}^2_{\la k \ra}\subset M.\eqno(7.61)$$
Let $r\in\mbb{Z}$. Denote
$${\cal M}^5=\mbb{F}[x_{n_2+1},...,x_n,y_{n_2+1},...,y_n],\qquad{\cal
M}^5_{\la k\ra}={\cal M}^5\bigcap{\cal B}_{\la
k\ra},\;\;k\in\mbb{Z}.\eqno(7.62)$$ Then
$${\cal M}^2_{\la k \ra}=\bigoplus_{r\in\mbb{Z}}{\cal
M}^4_{\la r\ra}{\cal M}^5_{\la k-r\ra}\eqno(7.63)$$ (cf. (7.30)).
Fix $r\in\mbb{Z}$.
$${\cal M}^4_{\la r\ra}{\cal M}^5_{\la k- r\ra}\subset {\cal M}^2_{\la k\ra}\subset M.\eqno(7.64)$$
Moreover, ${\cal M}^1_{\la r\ra}\supset{\cal M}^4_{\la r\ra}$ is an
irreducible ${\cal G}_1$-module (cf. (7.27) and (7.29)) by Step 1.
Thus
$${\cal M}^1_{\la r\ra}{\cal M}^5_{\la k- r\ra}=U({\cal G}_1)({\cal M}^4_{\la r\ra}{\cal M}^5_{\la k-
r\ra})\subset M.\eqno(7.65)$$ Furthermore,
$${\cal B}_{\la k\ra}=\bigoplus_{r\in\mbb{Z}}{\cal M}^1_{\la r\ra}{\cal M}^5_{\la k-
r\ra}\subset M\eqno(7.66)$$ by (7.27) and (7.65). Therefore,
$M={\cal B}_{\la k\ra}$, that is, ${\cal B}_{\la k\ra}$ is an
irreducible $sp(2n,\mbb{F})$-submodule. \psp

{\it Case 2}. $n_2=n_1+1$. \psp

 According to (3.104), the
nonzero weight vectors in
\begin{eqnarray*} & &\mbox{Span}\{\eta^{m_2}(x_i^{m_1}y_j^{m_3}),
x_{n_1+1}^{m_1}y_{n_1+1}^{m_2},\eta^{m_1+m_2}(x_{n_1}^{m_2}y_{n_1+1}^{m_3-m_1}),\eta^{m_1+m_2}(y_{n_1+2}^{m_2}x_{n_1+1}^{m_3-m_1})
\\ & &\qquad\;\;\mid m_r\in\mbb{N};
(i,j)=(n_1,n_1+1),(n_1,n_1+2),(n_1+1,n_1+2)\}.\hspace{2.5cm}(7.67)\end{eqnarray*}
are all the singular vectors of $sl(n,\mbb{F})$ in ${\cal B}$. Fix
$k\in\mbb{N}$. Then the singular vectors of $sl(n,\mbb{F})$  in
${\cal B}_{\la-k\ra}$ are those in (7.34). According to the
arguments in Case 1, ${\cal B}_{\la-k\ra}$ is an irreducible
$sp(2n,\mbb{F})$-submodule. Let $0<k\in\mbb{N}$. Then the singular
vectors of $sl(n,\mbb{F})$ in ${\cal B}_{\la k\ra}$ are
\begin{eqnarray*}\!\!\!\!\!\!\!\! &\{\eta^{m_2}(x_{n_1+1}^{k+m_1-2m_2}y_{n_1+2}^{m_1}),
\eta^{m_2}(x_{n_1}^{m_1}y_{n_1+1}^{k+m_1-2m_2}),
\eta^{m_5+m_6}(x_{n_1}^{m_6}y_{n_1+1}^{m_7-m_5}),\eta^{m_5+m_6}(y_{n_1+2}^{m_6}x_{n_1+1}^{m_7-m_5}),\hspace{0.6cm}
\\
 &\!\!\!\!\!x_{n_1+1}^{m_3}y_{n_1+1}^{m_4} \mid m_i\in\mbb{N};2m_2\leq
k+m_1;m_3+m_4=k=m_5+m_5+m_7\}\hspace{2.1cm}(7.68)\end{eqnarray*} by
(7.67).  Let $M$ be a nonzero $sp(2n,\mbb{F})$-submodule of ${\cal
B}_{\la k\ra}$. As an $sl(n,\mbb{F})$-module, $M$ contains a
singular of $sl(n,\mbb{F})$. If $x_{n_1+1}^{m_3}y_{n_1+1}^{m_4}\in
M$ with $m_3+m_4=k$, then $E_{n_1+1,n+n_1+1}|_{\cal
B}=x_{n_1+1}\ptl_{y_{n_1+1}}$ and
$$E_{n_1+1,n+n_1+1}^{m_4}(x_{n_1+1}^{m_3}y_{n_1+1}^{m_4})=m_4!x_{n_1+1}^k\in M\lra x_{n_1+1}^k\in M\eqno(7.69)$$
by (7.5). Suppose some
$\eta^{m_5+m_6}(x_{n_1}^{m_6}y_{n_1+1}^{m_7-m_5})\in M$ with
$m_5+m_5+m_7=k$. According to (1.16),
$E_{n+n_1+1,n_1+1}=y_{n_1+1}\ptl_{x_{n_1+1}}$. So
$$E_{n+n_1+1,n_1+1}^{m_5+m_6}[\eta^{m_5+m_6}(x_{n_1}^{m_6}y_{n_1+1}^{m_7-m_5})]=(m_5+m_6)!
x_{n_1}^{m_6}y_{n_1+1}^{k+m_6}\in M.\eqno(7.70)$$ Moreover, (7.4)
yields that $(E_{n_1,n+n_1+1}+E_{n_1+1,n+n_1})|_{\cal
B}=\ptl_{x_{n_1}}\ptl_{y_{n_1+1}}+x_{n_1+1}\ptl_{y_{n_1}}$ and
$$(E_{n_1,n+n_1+1}+E_{n_1+1,n+n_1})^{m_6}(x_{n_1}^{m_6}y_{n_1+1}^{k+m_6})=m_6![\prod_{r=1}^{m_6}(k+r)]
y_{n_1+1}^k\in M. \eqno(7.71)$$ Assume some
$\eta^{m_5+m_6}(y_{n_1+2}^{m_6}x_{n_1+1}^{m_7-m_5})\in M$ with
$m_5+m_5+m_7=k$. By (7.3) and (7.5),
$$E_{n_1+1,n+n_1+1}^{m_5+m_6}[\eta^{m_5+m_6}(y_{n_1+2}^{m_6}x_{n_1+1}^{m_7-m_5})]=(m_5+m_6)!
y_{n_1+2}^{m_6}x_{n_1+1}^{k+m_6}\in M.\eqno(7.72)$$ Observe
$$(E_{n+n_1+2,n_1+1}+E_{n+n_1+1,n_1+2})|_{\cal B}=
\ptl_{x_{n_1+1}}\ptl_{y_{n_1+2}}+y_{n_1+1}\ptl_{x_{n_1+2}}\eqno(7.73)$$
by (1.16). Hence
$$(E_{n+n_1+2,n_1+1}+E_{n+n_1+1,n_1+2})^{m_6}(y_{n_1+2}^{m_6}x_{n_1+1}^{k+m_6})
=m_6![\prod_{r=1}^{m_6}(k+r)] x_{n_1+1}^k\in M.\eqno(7.74)$$

Expressions (7.53)-(7.60), (7.69), (7.71) and (7.74) show that we
always have $x_{n_1+1}^k\in M$. Furthermore, (7.61)-(7.66) imply
that ${\cal B}_{\la k\ra}$ is an irreducible
$sp(2n,\mbb{F})$-module. \psp

{\it Case 3}. $n_1=n_2$.\psp

In this case,
$$\eta=\sum_{i=1}^{n_1}y_i\ptl_{x_i}+\sum_{s=n_2+1}^n
x_s\ptl_{y_s}.\eqno(7.75)$$ First we consider the subcase
$1<n_1<n-1$. Expression (4.17) says that the nonzero weight vectors
in
$$\mbox{Span}\{x_{n_1}^{m_1}y_{n_1}^{m_2}\zeta_1^{m_3+1},
x_{n_1+1}^{m_1}y_{n_1+1}^{m_2}\zeta_2^{m_3+1},
\eta^{m_3}(x_{n_1}^{m_1}y_{n_1+1}^{m_2})\mid
m_i\in\mbb{N}\}\eqno(7.76)$$ are all the singular vectors of
$sl(n,\mbb{F})$ in ${\cal B}$, where
$$\zeta_1=x_{n_1-1}y_{n_1}-x_{n_1}y_{n_1-1},\;\;\zeta_2=x_{n_1+1}y_{n_1+2}-x_{n_1+2}y_{n_1+1}.
\eqno(7.77)$$ Fix $k\in\mbb{N}+1$. Then the singular vectors of
$sl(n,\mbb{F})$  in ${\cal B}_{\la-k\ra}$ are
\begin{eqnarray*}\qquad&
&\{x_{n_1}^{k+m_1}y_{n_1}^{m_1}\zeta_1^{m_2+1},
x_{n_1+1}^{m_1}y_{n_1+1}^{k+m_1}\zeta_2^{m_2+1},
\eta^{m_3}(x_{n_1}^{m_4}y_{n_1+1}^{m_5}) \\
& &\mid
m_i\in\mbb{N};m_4+m_5-2m_3=k\}.\hspace{7.7cm}(7.78)\end{eqnarray*}
Let $M$ be a nonzero $sp(2n,\mbb{F})$-submodule of ${\cal
B}_{\la-k\ra}$. As an $sl(n,\mbb{F})$-module, $M$ contains a
singular vector of $sl(n,\mbb{F})$. Suppose some
$x_{n_1}^{k+m_1}y_{n_1}^{m_1}\zeta_1^{m_2+1}\in M$. Note
$E_{n_1,n+n_1}|_{\cal B}=\ptl_{x_{n_1}}\ptl_{y_{n_1}}$ by (7.4), and
so
\begin{eqnarray*} & &
E_{n_1,n+n_1}(x_{n_1}^{k+m_1}y_{n_1}^{m_1}\zeta_1^{m_2}) \\ &=&
(k+m_1)m_1x_{n_1}^{k+m_1-1}y_{n_1}^{m_1-1}\zeta_1^{m_2}-m_2(m_2-1)
x_{n_1}^{k+m_1}y_{n_1}^{m_1}x_{n_1-1}y_{n_1-1}\zeta_1^{m_2-2}
\\ & &+(k+m_1)m_2x_{n_1}^{k+m_1-1}y_{n_1}^{m_1}x_{n_1-1}\zeta_1^{m_2-1}
-m_1m_2x_{n_1}^{k+m_1}y_{n_1}^{m_1-1}y_{n_1-1}\zeta_1^{m_2-1}.
\hspace{1.7cm}(7.79)\end{eqnarray*} Moreover,
$$(E_{n_1-1,n_1}-E_{n+n_1,n+n_1-1})|_{\cal
B}=-(x_{n_1}\ptl_{x_{n_1-1}}+y_{n_1}\ptl_{y_{n_1-1}})\eqno(7.80)$$
by (1.7), (1.8) and (1.14). Thus
\begin{eqnarray*} \qquad\qquad& &
(E_{n_1-1,n_1}-E_{n+n_1,n+n_1-1})^2E_{n_1,n+n_1}(x_{n_1}^{k+m_1}y_{n_1}^{m_1}\zeta_1^{m_2})
\\ &=&
-2m_2(m_2-1) x_{n_1}^{k+m_1+1}y_{n_1}^{m_1+1}\zeta_1^{m_2-2}\in
M.\hspace{5.1cm}(7.81)\end{eqnarray*} Hence
$$x_{n_1}^{k+m_1+1}y_{n_1}^{m_1+1}\zeta_1^{m_2-2}\in
M\qquad\mbox{if}\;\;m_2>1.\eqno(7.82)$$ Furthermore,
$$(E_{n_1-1,n_1}-E_{n+n_1,n+n_1-1})E_{n_1,n+n_1}(x_{n_1}^{k+m_1}y_{n_1}^{m_1}\zeta_1)
= -kx_{n_1}^{k+m_1}y_{n_1}^{m_1}\in M.\eqno(7.83)$$  So we always
have $x_{n_1}^{k+m}y_{n_1}^m\in M$ for some $m\in\mbb{N}$ by
induction on $m_2$.

 Observe
 $$E_{n_1,n+n_1}(x_{n_1}^{k+m}y_{n_1}^m)=\ptl_{x_{n_1}}\ptl_{y_{n_1}}(x_{n_1}^{k+m}y_{n_1}^m)=m![\prod_{r=1}^m
 (k+r)]x_{n_1}^k\eqno(7.84)$$
 by (7.4). Thus
 $$x_{n_1}^k\in M.\eqno(7.85)$$
Symmetrically, if some
$x_{n_1+1}^{m_1}y_{n_1+1}^{k+m_1}\zeta_2^{m_2+1}\in M$, we have
$y_{n_1+1}^k\in M$. But
$$(E_{n+n_1+1,n_1}+E_{n+n_1,n_1+1})|_{\cal
B}=-x_{n_1}\ptl_{y_{n_1+1}}+y_{n_1}\ptl_{x_{n_1+1}}\eqno(7.86)$$ by
(1.16), which gives
$$(E_{n+n_1+1,n_1}+E_{n+n_1,n_1+1})^k(y_{n_1+1}^k)=(-1)^kk!x_{n_1}^k\in
M.\eqno(7.87)$$ Thus (7.85) holds again.

Assume that some $\eta^{m_3}(x_{n_1}^{m_4}y_{n_1+1}^{m_5})\in M$
with $m_4+m_5-2m_3=k$. Note there exists $r_1,r_2\in\mbb{N}$ such
that $r_1+r_2=m_3$ and $2r_1\leq m_4,\;2r_2\leq m_5$. Moreover, $
E_{n_1,n+n_1}|_{\cal B}=\ptl_{x_{n_1}}\ptl_{y_{n_1}}$ by (7.4) and $
E_{n+n_1+1,n_1+1}|_{\cal B}=\ptl_{x_{n_1+1}}\ptl_{y_{n_1+1}}$ by
(1.16). Thus
\begin{eqnarray*}\qquad & &E_{n_1,n+n_1}^{r_1}E_{n+n_1+1,n_1+1}^{r_2}
[\eta^{m_3}(x_{n_1}^{m_4}y_{n_1+1}^{m_5})]\\
&=&m_3![\prod_{s_1=0}^{2r_1-1}(m_4-s_1)]
[\prod_{s_2=0}^{2r_2-1}(m_5-s_2)]x_{n_1}^{m_4-2r_1}y_{n_1+1}^{m_5-2r_2}\in
M. \hspace{3.2cm}(7.88)\end{eqnarray*} Furthermore, (1.16) gives
$(E_{n+n_1+1,n_1}+E_{n+n_1,n_1+1})|_{\cal
B}=-x_{n_1}\ptl_{y_{n_1+1}}+y_{n_1}\ptl_{x_{n_1+1}}$ , and so
\begin{eqnarray*}\hspace{2cm}& &(E_{n+n_1+1,n_1}+E_{n+n_1,n_1+1})^{m_5-2r_2}
(x_{n_1}^{m_4-2r_1}y_{n_1+1}^{m_5-2r_2})\\
&=&(-1)^{m_5-2r_2}(m_5-2r_2)!x_{n_1}^k\in
M.\hspace{6.2cm}(7.89)\end{eqnarray*} Thus we always have
$x_{n_1}^k\in M$.

Now
$$(E_{n_1,n+n_1+1}+E_{n_1+1,n+n_1})|_{\cal
B}=-y_{n_1+1}\ptl_{x_{n_1}}+x_{n_1+1}\ptl_{y_{n_1}}\eqno(7.90)$$ by
(7.4). For any $r\in\mbb{N}+1$,
$$\frac{(-1)^r}{\prod_{s=0}^{r-1}(k-r)}(E_{n_1,n+n_1+1}+E_{n_1+1,n+n_1})^r(x_{n_1}^k)
=x_{n_1}^{k-r}y_{n_1+1}^r\in M.\eqno(7.91)$$ If $k\geq 2$ and
$r\in\ol{1,k-1}$, then
$${\cal M}^1_{\la -k+r\ra}{\cal M}^2_{\la -r\ra}=U({\cal G}_1)U({\cal
G}_2)(x_{n_1}^{k-r}y_{n_1+1}^r)\subset M\eqno(7.92)$$ because ${\cal
M}^1_{\la -k+r\ra}$ is an irreducible ${\cal G}_1$-module and ${\cal
M}^2_{\la -r\ra}$ is an irreducible ${\cal G}_2$-module by Step 1.
Moreover,
$${\cal M}^1_{\la -k\ra}=U({\cal G}_1)(x_{n_1}^k),\;{\cal M}^2_{\la -k\ra}=U({\cal G}_2)
(y_{n_1+1}^k)\subset M.\eqno(7.93)$$ Furthermore,
$${\cal M}^1_{\la -k\ra}{\cal M}^2_{\la 0\ra}=U({\cal G}_1)U({\cal
G}_2)(x_{n_1}^k)\subset M \;\;\mbox{if}\;\;n_1=n-1\eqno(7.94)$$ and
$${\cal M}^1_{\la 0\ra}{\cal M}^2_{\la -k\ra}=U({\cal G}_1)U({\cal G}_2)
(y_{n_1+1}^k)\subset M \;\;\mbox{if}\;\;n_1=1.\eqno(7.95)$$

Note
$$(E_{r,i}-E_{n+i,n+r})|_{\cal
B}=y_iy_r-x_ix_r\qquad\for\;\;i\in\ol{1,n_1},\;r\in\ol{n_1+1,n}\eqno(7.96)$$
by (1.7), (1.8) and (1.14). In particular, if $k>1$ or $n_1=1$, we
have
$$(E_{n_1+1,n_1}-E_{n+n_1,n+n_1+1})(x_{n_1}^k)=y_{n_1}x_{n_1}^ky_{n_1+1}-x_{n_1}^{k+1}x_{n_1+1}\in
M.\eqno(7.97)$$ Since
$$y_{n_1}x_{n_1}^ky_{n_1+1}\in {\cal M}^1_{\la -k+1\ra}{\cal M}^2_{\la
-1\ra}\subset M,\eqno(7.98)$$ we get $$ x_{n_1}^{k+1}x_{n_1+1}\in
M.\eqno(7.99)$$ Suppose $k=1$ and $n_1>1$. By (7.93),
$$\zeta_1x_{n_1}=(x_{n_1-1}y_{n_1}-x_{n_1}y_{n_1-1})x_{n_1}\in
M.\eqno(7.100)$$ Observe
$$(E_{n_1+1,n+n_1-1}+E_{n_1-1,n+n_1+1})|_{\cal
B}=x_{n_1+1}\ptl_{y_{n_1-1}}-y_{n_1+1}\ptl_{x_{n_1-1}}\eqno(7.101)$$
by (1.15). So
$$-(E_{n_1+1,n+n_1-1}+E_{n_1-1,n+n_1+1})(\zeta_1x_{n_1})=x_{n_1}^2x_{n_1+1}
-x_{n_1}y_{n_1}y_{n_1+1}\in M.\eqno(7.102)$$ On the other hand,
(1.16) gives
$$(E_{n+i,j}+E_{n+j,i})|_{\cal
B}=-(x_iy_j+x_jy_i)\qquad\for\;\;i,j\in\ol{1,n_1},\eqno(7.103)$$
which implies
$$-E_{n+n_1,n_1}(y_{n_1+1})=x_{n_1}y_{n_1}y_{n_1+1}\in
M.\eqno(7.104)$$ By (7.102), we have $x_{n_1}^2x_{n_1+1}\in M.$ So
(7.99) always holds.

By Step 1,
$${\cal M}^1_{\la -k-1\ra}{\cal M}^2_{\la 1\ra}=U({\cal G}_1)U({\cal
G}_2)(x_{n_1}^{k+1}x_{n_1+1})\subset M.\eqno(7.105)$$ Suppose
$${\cal M}^1_{\la -k-i\ra}{\cal M}^2_{\la i\ra}\subset
M\eqno(7.106)$$ for $1\leq i\leq m$. Then
\begin{eqnarray*}\qquad\qquad& &(E_{n_1+1,n_1}-E_{n+n_1,n+n_1+1})(x_{n_1}^{k+m}x_{n_1+1}^m)
\\ &=&y_{n_1}x_{n_1}^{k+m}x_{n_1+1}^my_{n_1+1}-x_{n_1}^{k+m+1}x_{n_1+1}^{m+1}\in
M\hspace{5cm}(7.107)\end{eqnarray*} by (7.96). If $m>1$, we have
$$y_{n_1}x_{n_1}^{k+m}x_{n_1+1}^my_{n_1+1}\in
{\cal M}^1_{\la -k-(m-1)\ra}{\cal M}^2_{\la m-1\ra}\subset
M.\eqno(7.108)$$ Note
$$(E_{r,n+s}+E_{s,n+r})|_{\cal
B}=-(x_ry_s+x_sy_r)\qquad\for\;\;r,s\in\ol{n_1+1,n}\eqno(7.109)$$ by
(1.15). If $m=1$, we have
$$y_{n_1}x_{n_1}^{k+1}x_{n_1+1}y_{n_1+1}
=-E_{n_1+1,n+n_1+1}(y_{n_1}x_{n_1}^{k+1})\subset
E_{n_1+1,n+n_1+1}({\cal M}^1_{\la -k\ra})\subset M.\eqno(7.110)$$
Then (7.107), (7.108) and (7.110) give
$$x_{n_1}^{k+m+1}x_{n_1+1}^{m+1}\in
M.\eqno(7.111)$$ Furthermore,
$${\cal M}^1_{\la -k-m-1\ra}{\cal M}^2_{\la m+1\ra}=U({\cal G}_1)U({\cal
G}_2)(x_{n_1}^{k+m+1}x_{n_1+m+1})\subset M.\eqno(7.112)$$ Thus
(7.106) holds for any $i\in\mbb{N}+1$. Symmetrically, we have
$${\cal M}^1_{\la i\ra}{\cal M}^2_{\la -k-i\ra}\subset M\qquad\for\;\;i\in\mbb{N}+1.\eqno(7.113)$$

Suppose $n_1<n-1$. Then $x_{n_1}^{k+1}x_{n_1+1}\zeta_2\in M$ by
(7.105). Moreover,
$$(k+1)y_{n_1}x_{n_1}^ky_{n_1+1}\zeta_2=-(k+1)E_{n+n_1,n_1}(x_{n_1}^{k-1}y_{n_1+1}\zeta_2)\in
M\eqno(7.114)$$ by (7.92) and (7.103). According (1.7), (1.8) and
(1.14),
$$(E_{n_1,n_1+1}-E_{n+n_1+1,n+n_1})|_{\cal
B}=\ptl_{x_{n_1}}\ptl_{x_{n_1+1}}-\ptl_{y_{n_1}}\ptl_{y_{n_1+1}}.\eqno(7.115)$$
Thus
\begin{eqnarray*} \qquad& &(E_{n_1,n_1+1}-E_{n+n_1+1,n+n_1})
[(x_{n_1}^{k+1}x_{n_1+1}-(k+1)y_{n_1}x_{n_1}^ky_{n_1+1})\zeta_2]
\\ &=& 3(k+1)x_{n_1}^k\zeta_2\in
M\hspace{9.6cm}(7.116)\end{eqnarray*} by (7.77). Hence
$${\cal M}^1_{\la-k\ra}{\cal M}^2_{\la 0\ra}=U({\cal G}_1)U({\cal
G}_2)(x_{n_1}^k)+U({\cal G}_1)U({\cal G}_2)(x_{n_1}^k\zeta_2)\subset
M\eqno(7.117)$$ by (7.26) and(7.85). Symmetrically,
$${\cal M}^1_{\la 0 \ra}{\cal M}^2_{\la -k\ra}\subset
M.\eqno(7.118)$$ By (7.92)-(7.95), (7.106), (7.112), (7.113),
(7.117) and (7.118),
$${\cal M}^1_{\la -k-r \ra}{\cal M}^2_{\la r\ra}\subset
M\qquad\for\;\;r\in\mbb{Z}.\eqno(7.119)$$ Therefore,
$${\cal B}_{\la -k\ra}=\bigoplus_{r\in\mbb{Z}}
{\cal M}^1_{\la -k-r \ra}{\cal M}^2_{\la r\ra}\subset
M.\eqno(7.120)$$ We get $M={\cal B}_{\la -k\ra}$, that is, ${\cal
B}_{\la -k\ra}$ is an irreducible $sp(2n,\mbb{F})$-module. We can
similarly prove that ${\cal B}_{\la k\ra}$ is an irreducible
$sp(2n,\mbb{F})$-module.

Finally, we study ${\cal B}_{\la 0\ra}$. We first consider the
generic case $1<n_1<n-1$. Set
\begin{eqnarray*} \qquad
{\cal B}_{\la 0,1\ra}&=&\mbox{Span}\{[\prod_{1\leq r\leq s\leq
n_1\;\mbox{or}\;n_1+1\leq r\leq s\leq n}(x_ry_s+x_sy_r)^{l_{r,s}}]\\
& &\qquad\;\;\times
[\prod_{p=1}^{n_1}\prod_{q=n_1+1}^n(x_px_q-y_py_q)^{k_{p,q}}] \mid
l_{r,s},k_{p,q}\in\mbb{N}\}\hspace{3cm}(7.121)\end{eqnarray*} and
$${\cal B}_{\la
0,2\ra}=\sum_{1\leq r< s\leq n_1\;\mbox{or}\;n_1+1\leq r< s\leq
n}{\cal B}_{\la
0,1\ra}(x_ry_s-x_sy_r)+\sum_{p=1}^{n_1}\sum_{q=n_1+1}^n{\cal B}_{\la
0,1\ra}(x_px_q+y_py_q).\eqno(7.122)$$ We want to prove that ${\cal
B}_{\la 0,1\ra}$ and ${\cal B}_{\la 0,2\ra}$ forms
$sp(2n,\mbb{F})$-submodules.

Let \begin{eqnarray*}\qquad {\cal G}'&=&\sum_{1\leq r\leq s\leq
n_1}\mbb{F}(E_{n+s,r}+E_{n+r,s})+\sum_{n_1+1\leq p\leq q\leq
n}\mbb{F}(E_{p,n+q}+E_{q,n+p})\\ &
&+\sum_{r=1}^{n_1}\sum_{p=n_1+1}^n
\mbb{F}(E_{p,r}-E_{n+r,n+p})\hspace{7cm}(7.123)\end{eqnarray*} and
\begin{eqnarray*}\hat{\cal
G}&=&\sum_{i,j=1}^{n_1}\mbb{F}(E_{i,j}-E_{n+j,n+i})+\sum_{r,s=n_1+1}^n\mbb{F}
(E_{r,s}-E_{n+s,n+r})+\sum_{1\leq r\leq s\leq
n_1}\mbb{F}(E_{r,n+s}+E_{s,n+r})\\ & &+\sum_{n_1+1\leq p\leq q\leq
n}\mbb{F}(E_{n+q,p}+E_{n+p,q})+\sum_{r=1}^{n_1}\sum_{p=n_1+1}^n
[\mbb{F}(E_{r,p}-E_{n+p,n+r})\\
&
&+\mbb{F}(E_{r,n+p}+E_{p,n+r})+\mbb{F}(E_{n+r,p}-E_{n+p,r})].\hspace{5.8cm}(7.124)\end{eqnarray*}
Then ${\cal G}'$ and $\hat{\cal G}$ are Lie subalgebras of
$sp(2n,\mbb{F})$ and $sp(2n,\mbb{F})={\cal G}'\oplus \hat{\cal G}.$
By PBW Theorem $U(sp(2n,\mbb{F}))=U({\cal G}')U(\hat{\cal G}).$

By (7.96), (7.103) and (7.109),
$$U({\cal G}')|_{\cal B}={\cal B}_{\la
0,1\ra}\;\;\mbox{as multiplication operators on}\;\;{\cal
B}.\eqno(7.125)$$ Moreover,
$$(E_{r,s}-E_{n+s,n+r})|_{\cal
B}=x_r\ptl_{x_s}+y_r\ptl_{y_s}+\dlt_{r,s},\eqno(7.126)$$
$$(E_{n+r,s}+E_{n+s,r})|_{\cal
B}=\ptl_{x_r}\ptl{y_s}+\ptl_{x_s}\ptl{y_r},\eqno(7.127)$$
$$(E_{n+r,i}+E_{n+i,r})|_{\cal
B}=-x_i\ptl{y_r}+y_i\ptl{x_r},\eqno(7.128)$$
$$(E_{i,n+r}+E_{r,n+i})|_{\cal
B}=-y_r\ptl{x_i}+x_r\ptl{y_i},\eqno(7.129)$$
$$(E_{i,r}-E_{n+r,n+i})|_{\cal
B}=\ptl_{x_i}\ptl{x_r}-\ptl_{y_i}\ptl{y_r}\eqno(7.130)$$ for
$i\in\ol{1,n_1}$ and $r,s\in\ol{n_1+1,n}$. According to (7.4),
(7.6), (7.124) and (7.126)-(7.130), $U(\hat{\cal G})(1)=\mbb{F}.$
Thus
$${\cal B}_{\la 0,1\ra}=U({\cal
G}')(1)=U(sp(2n,\mbb{F}))(1)\eqno(7.131)$$ forms an
$sp(2n,\mbb{F})$-submodule.

Let
$$W=\sum_{1\leq r< s\leq n_1\;\mbox{or}\;n_1+1\leq r< s\leq
n}\mbb{F}(x_ry_s-x_sy_r)+
\sum_{p=1}^{n_1}\sum_{q=n_1+1}^n\mbb{F}(x_px_q+y_py_q).\eqno(7.132)$$
By  (7.4), (7.6) and (7.126)-(7.130), we can verify that $W$ forms
an irreducible $\hat{\cal G}$-submodule. Hence
$${\cal B}_{\la 0,2\ra}=U({\cal
G}')(W)= U(sp(2n,\mbb{F}))(W)\eqno(7.133)$$ forms an
$sp(2n,\mbb{F})$-submodule. Moreover,
$${\cal B}_{\la 0,1\ra}\bigcap W=\{0\}.\eqno(7.134)$$

Next we want to prove that ${\cal B}_{\la 0,1\ra}$ and ${\cal
B}_{\la 0,2\ra}$ are irreducible $sp(2n,\mbb{F})$-submodules.
According to (7.78), the singular vectors of $sl(n,\mbb{F})$  in
${\cal B}_{\la 0 \ra}$ are
\begin{eqnarray*}\qquad&
&\{x_{n_1}^{m_1}y_{n_1}^{m_1}\zeta_1^{m_2+1},
x_{n_1+1}^{m_1}y_{n_1+1}^{m_1}\zeta_2^{m_2+1},
\eta^{m_3}(x_{n_1}^{m_4}y_{n_1+1}^{m_5}) \\
& &\mid
m_i\in\mbb{N};m_4+m_5=2m_3\}.\hspace{8.2cm}(7.135)\end{eqnarray*}
Let $M$ be a nonzero submodule of ${\cal B}_{\la 0,1\ra}$. Then $M$
contains a singular vector of $sl(n,\mbb{F})$. Suppose some
$x_{n_1}^{m_1}y_{n_1}^{m_1}\zeta_1^{m_2}\in M$. By (7.79)-(7.82), we
can assume $m_2=0,1$. If $m_2=0$, (7.84) yields $1\in M$. Then
$M={\cal B}_{\la 0,1\ra}$ by (7.131). Suppose $m_2=1$. We have
$E_{n_1,n+n_1}|_{\cal B}=\ptl_{x_{n_1}}\ptl_{y_{n_1}}$ by (7.4), and
$$E_{n_1,n+n_1}[x_{n_1}^{m_1}y_{n_1}^{m_1}\zeta_1]
=m_1(m_1+1)x_{n_1}^{m_1-1}y_{n_1}^{m_1-1}\zeta_1\eqno(7.136)$$ by
(7.79). By induction on $m_1$, we have $\zeta_1\in M\subset {\cal
B}_{\la 0,1\ra}$, which contradicts (7.134). Similarly, if some
$x_{n_1+1}^{m_1}y_{n_1+1}^{m_1}\zeta_2^{m_2+1}\in M$, we have $M=
{\cal B}_{\la 0,1\ra}$. Assume some
$\eta^{m_3}(x_{n_1}^{m_4}y_{n_1+1}^{m_5})\in M$ with $m_4+m_5=2m_3$.
Note $m_4$ and $m_5$ are both even or odd. If $m_4=2r_1$ and
$m_5=2r_2$ are even, then (7.88) gives $1\in M$, equivalently
$M={\cal B}_{\la 0,1\ra}$. Suppose that $m_4=2r_1+1$ and
$m_5=2r_2+1$ are odd. Expression (7.75) yields
$$\eta(x_{n_1}y_{n_1+1})=x_{n_1}x_{n_1+1}
+y_{n_1}y_{n_1+1}\in M\subset {\cal B}_{\la 0,1\ra},\eqno(7.137)$$
 which contradicts (7.134) again. Thus we always have
$M={\cal B}_{\la 0,1\ra}$, that is, ${\cal B}_{\la 0,1\ra}$ is
irreducible. Similarly, we can prove that ${\cal B}_{\la 0,2\ra}$ is
irreducible.

If $n_1=1$ and $n=2$, we let $${\cal B}_{\la
0,1\ra}=\mbox{Span}\{[\prod_{i=1}^n(x_iy_i)^{m_i}]
(x_1x_2-y_1y_2)^{m_3}] \mid m_i\in\mbb{N}\}\eqno(7.138)$$
 and
${\cal B}_{\la 0,2\ra}={\cal B}_{\la 0,1\ra}(x_1x_2+y_1y_2).$ When
$n_1=1$ and $n>2$, we set
$${\cal B}_{\la 0,1\ra}=\mbox{Span}\{[(x_1y_1)^l\prod_{2\leq r\leq s\leq
n}(x_ry_s+x_sy_r)^{l_{r,s}}] [\prod_{q=2}^n(x_1x_q-y_1y_q)^{k_q}]
\mid l, l_{r,s},k_{q}\in\mbb{N}\}\eqno(7.139)$$ and
$${\cal B}_{\la
0,2\ra}=\sum_{2\leq r< s\leq n}{\cal B}_{\la
0,1\ra}(x_ry_s-x_sy_r)+\sum_{q=2}^n{\cal B}_{\la
0,1\ra}(x_1x_q+y_1y_q).\eqno(7.140)$$ In the case $1<n_1=n-1$, we
put
\begin{eqnarray*} \qquad
{\cal B}_{\la 0,1\ra}&=&\mbox{Span}\{(x_ny_n)^l[\prod_{1\leq r\leq
s\leq n-1}(x_ry_s+x_sy_r)^{l_{r,s}}]\\ & &\qquad\;\;\times
[\prod_{p=1}^{n-1}(x_px_n-y_py_n)^{k_{p}}] \mid
l,l_{r,s},k_{p}\in\mbb{N}\}\hspace{4.2cm}(7.141)\end{eqnarray*} and
$${\cal B}_{\la
0,2\ra}=\sum_{1\leq r< s\leq n_1}{\cal B}_{\la
0,1\ra}(x_ry_s-x_sy_r)+\sum_{p=1}^{n_1}{\cal B}_{\la
0,1\ra}(x_px_n+y_py_n).\eqno(7.142)$$ The above corresponding
partial arguments show that ${\cal B}_{\la 0,1\ra}$ and ${\cal
B}_{\la 0,2\ra}$ are irreducible in the corresponding case.

Now $1$ is a non-isotropic element in
 ${\cal B}_{\la 0,1\ra}$ and $x_{n_1}x_{n_1+1}+y_{n_1}y_{n_1+1}$ a non-isotropic element in ${\cal B}_{\la 0,2\ra}$
by (3.54). By Lemma 2.3, the symmetric bilinear form $(\cdot|\cdot)$
restricted to them are nondegenrate. Since $(1|{\cal B}_{\la
0,2\ra})=\{0\}$, ${\cal B}_{\la 0,1\ra}$ is orthogonal to ${\cal
B}_{\la 0,2\ra}$. Thus the symmetric bilinear form $(\cdot|\cdot)$
restricted ${\cal B}_{\la 0,1\ra}+{\cal B}_{\la 0,2\ra}$ is
nondegenerate. Then
$${\cal B}_{\la 0\ra}=({\cal B}_{\la 0,1\ra}+{\cal B}_{\la
0,2\ra})\oplus ({\cal B}_{\la 0,1\ra}+{\cal B}_{\la
0,2\ra})^\perp\bigcap {\cal B}_{\la 0\ra}.\eqno(7.143)$$ If $({\cal
B}_{\la 0,1\ra}+{\cal B}_{\la 0,2\ra})^\perp\bigcap {\cal B}_{\la
0\ra}\neq\{0\}$, then it contains a singular vector of
$sl(n,\mbb{F})$. Our above arguments in proving the irreducibility
of ${\cal B}_{\la 0,1\ra}$ show that it contains either ${\cal
B}_{\la 0,1\ra}$ or ${\cal B}_{\la 0,2\ra}$, which is absurd.
Therefore, ${\cal B}_{\la 0\ra}={\cal B}_{\la 0,1\ra}\oplus{\cal
B}_{\la 0,2\ra}$ is an orthogonal decomposition of irreducible
$sp(2n,\mbb{F})$-submodules.

Suppose $n_1=n_2=n$. For $k\in\mbb{N}+1$, (4.21) and (4.31) imply
that
$${\cal B}_{\la k\ra}=\bigoplus_{m=0}^{\infty}\bigoplus_{r=\llbracket (k+1)/2
\rrbracket}^\infty\eta^r({\cal H}_{\la k-2r-m,m\ra})\eqno(7.144)$$
and
$${\cal B}_{\la -k\ra}=\bigoplus_{m,r=0}^{\infty}\eta^r({\cal H}_{\la -k-2r-m,m\ra})\eqno(7.145)$$
are $({\cal G},{\cal K})$-structures, where ${\cal H}_{\la
-m_1-m_2,m_2\ra}$ is given in (4.42). Moreover,
$${\cal B}_{\la 0,1\ra}=\bigoplus_{m,r=0}^{\infty}\eta^r({\cal H}_{\la -2r-2m,2m\ra})\eqno(7.146)$$
and
$${\cal B}_{\la 0,2\ra}=\bigoplus_{m,r=0}^{\infty}\eta^r({\cal H}_{\la -2r-2m-1,2m+1\ra})\eqno(7.147)$$
are $({\cal G},{\cal K})$-structures by the arguments in
(7.79)-(7.82), (7.84) and (7.136) (cf. (7.24), (7.25)). $\qquad\Box$

\vspace{1cm}

\noindent{\Large \bf References}

\hspace{0.5cm}

\begin{description}

\item[{[C]}] B. Cao, Solutions of Navier Equations and Their
Representation Structure, {\it Advances in Applied Mathematics,}
\textbf{43} (2009), 331-374.

\item[{[DES]}] M. Davidson, T.
Enright, and R. Stanke, {\it Differential Operators and Highest
Weight Representations}, Memoirs of American Mathematical Society
{\bf 94}, no. 455, 1991.

\item[{[FC]}] F. M. Fern\'{a}ndez and E. A. Castro, {\it Algebraic
Methods in Quantum Chemistry and Physics}, CRC Press, Inc., 1996.

\item[{[FSS]}] L.Frappat, A.Sciarrino and P.sorba, Dictionary on Lie
Algebras and Superalgebras, Academic Press,2000

\item[{[G]}] H. Georgi, {\it Lie Algebras in Particle Physics},
Second Edition, Perseus Books Group, 1999.

\item[{[Ho]}] R. Howe, Perspectives on invariant theory: Schur
duality, multiplicity-free actions and beyond, {\it The Schur
lectures} (1992) ({\it Tel Aviv}), 1-182, {\it Israel Math. Conf.
Proc.,} 8, {\it Bar-Ilan Univ., Ramat Gan,} 1995.

\item[{[Hu]}] J. E. Humphreys, {\it Introduction to Lie Algebras and Representation Theory},
 Springer-Verlag New York Inc., 1972.

\item[{[K]}] V. G. Kac, {\it Infinite Dimensional Lie Algebras},
Third edition, Cambridge University Press, 1990.

\item[{[I1]}] N. H. Ibragimov, {\it Transformation groups applied to
mathematical physics}, Nauka, 1983.

\item[{[I2]}] N. H. Ibragimov, {\it Lie Group Analysis of
Differential Equations}, Volume 2, CRC Handbook, CRC Press, 1995.

\item[{[L]}] F. S. Levin, {\it An Introduction to Quantum Theory},
Cambridge University Press, 2002.

 \item[{[LF]}] W. Ludwig and C. Falter, {\it
Symmetries in Physics}, Second Edition, Springer-Verlag,
Berlin/Heidelberg, 1996.

\item[{[X]}] X.Xu, Flag partial differential equations and
representations of Lie algebras, {\it Acta Appl Math} {\bf 102}
(2008), 249--280.

\end{description}

\end{document}